%% file: main.tex
\documentclass[fleqn, preprint,review,11pt]{elsarticle}

\input{preamble}

\begin{document}

\begin{frontmatter}

\title{Non-intrusive Hybrid Scheme for Multiscale Heat Transfer: Thermal Runaway in a Battery Pack}                      

\author[1]{Yinuo Noah Yao\corref{contrib}}
\cortext[cor1]{Corresponding author}
\cortext[contrib]{Current address: Department of Civil and Environmental Engineering, Texas A\&M University, College Station, TX 77840, USA}
\author[2]{Perry Harabin}
\author[2]{Morad Behandish} 
\author[1]{Ilenia Battiato\corref{cor1}}

\address[1]{Department of Energy Science and Engineering, Stanford University, Stanford, CA 94305, USA}
\address[2]{Palo Alto Research Center (PARC), Palo Alto, CA 94305, USA}
% \address[3]{Department of Civil and Environmental Engineering, Texas A\&M University, College Station, TX 77840, USA}

\input{abstract}

\end{frontmatter}

\input{intro}

\input{metrics} % Morad is working on this

\input{equations}

\input{thermal}
\input{conclusion}

\section*{Acknowledgments}
This material is based upon work supported by the Defense Advanced Research Projects Agency (DARPA) under Agreement No. HR00112090061. The views, opinions and/or findings expressed are those of the authors and should not be interpreted as representing the official views or policies of the Department of Defense or U.S. Government. We also would like to thank Kyle Pietrzyk for the valuable inputs and feedbacks provided throughout all the stages of this work.

\newpage
\appendix

\input{BCs_pack_temp}

\input{BCs_cell_temp}

\input{closure}

\input{accuracy}

\input{coupling}

\newpage
\bibliographystyle{elsarticle-num}
%\bibliography{library_new}

\end{document}

%% file: preamble.tex
\usepackage[margin=1in]{geometry}

\usepackage{physics}
\usepackage{amsmath}
\usepackage{amssymb}
\usepackage{esint}
\usepackage{siunitx}
\usepackage{multicol}
\usepackage{pdflscape}
\usepackage{multirow}
\usepackage{multicol}
\usepackage{nicefrac}
\usepackage{booktabs}
\usepackage{float}
\usepackage{svg}
\usepackage{epstopdf}
\usepackage{xspace}
\usepackage[explicit]{titlesec}
\usepackage{xcolor}
\usepackage{verbatim}
\usepackage{mathrsfs}
\usepackage{hhline}
\usepackage{longtable}
\usepackage{tabu}
\usepackage{booktabs}
\usepackage{bm}

\usepackage{lineno}
\usepackage{caption}
\usepackage{subcaption}
\makeatletter
\let\LN@align\align
\let\LN@endalign\endalign
\renewcommand{\align}{\linenomath\LN@align}
\renewcommand{\endalign}{\LN@endalign\endlinenomath}
\let\LN@gather\gather
\let\LN@endgather\endgather
\renewcommand{\gather}{\linenomath\LN@gather}
\renewcommand{\endgather}{\LN@endgather\endlinenomath}
\makeatother

\newcommand{\ensmean}[1]{\left\langle #1 \right\rangle}

\let\oldfrac\frac% Store \frac
\renewcommand{\frac}[2]{%
  \mathchoice
    {\oldfrac{#1}{#2}}% display style
    {#1/#2}% text style
    {#1/#2}% script style
    {#1/#2}% script-script style
}

\newcommand{\ddfrac}[2]{\frac{\displaystyle #1}{\displaystyle #2}}
\def\K{\mathbf{K}}

\def\n{\mathbf{n}}
\def\R{\mathbf{R}}

\def\U{\mathbf{U}}

\def\V{\mathbf{V}}

\def\x{\mathbf{x}}

\def\chiv{\bm{\chi}}
\def\xiv{\bm{\xi}}

\def\I{\mathbf{I}}

\def\derive#1#2{\frac{\partial #1}{\partial #2}}

\def\neweq#1#2{\begin{equation} \label{#1}\begin{aligned} #2 \end{aligned}\end{equation}}
\def\neweqgat#1#2{\begin{equation} \label{#1}\begin{gathered} #2 \end{gathered}\end{equation}}

\def\xloop#1{\ifx\relax#1\else[#1]\expandafter\xloop\fi}

\def\lookfori#1{\ifx i#1 \expandafter\lookfors \else \expandafter\lookfori\fi}
\def\lookfors#1{\ifx s#1 \expandafter\lookforunderscore \else \expandafter\lookfori\fi}
\def\lookforunderscore#1{\ifx \_#1 \expandafter\savefirstletter \else \expandafter\lookfori\fi}
\def\savefirstletter#1{\def\letterone{#1}\expandafter\stopatnext}
\def\stopatnext#1{\ifx \_#1 \break\fi \savesecondletter#1}
\def\savesecondletter#1{\def\lettertwo{#1}\expandafter\stopatnext}

\DeclareMathOperator*{\argmax}{arg\,max}

\newcommand{\symbolica}{\textsf{Symbolica}\xspace}

\newcommand{\fenics}{\textsf{FEniCS}\xspace}

\journal{}

%% file: abstract.tex
\begin{abstract}
Accurate analytical and numerical modeling of multiscale systems is a daunting task. The need to properly resolve spatial and temporal scales spanning multiple orders of magnitude pushes the limits of both our theoretical models as well as our computational capabilities. Rigorous upscaling techniques enable efficient computation while bounding/tracking errors and helping to make informed cost-accuracy tradeoffs. The biggest challenges arise when the applicability conditions of upscaled models break down. Here, we present a non-intrusive two-way (iterative bottom-up top-down) coupled  hybrid model, applied to thermal runaway in battery packs, that combines fine-scale and upscaled equations in the same numerical simulation to achieve predictive accuracy while limiting computational costs. First, we develop two methods with different orders of accuracy to enforce continuity at the coupling boundary. Then, we derive weak (i.e., variational) formulations of the fine-scale and upscaled governing equations for finite element (FE) discretization and numerical implementation in \fenics. We demonstrate that hybrid simulations can accurately predict the average temperature fields within error bounds determined \emph{a priori} by homogenization theory. Finally, we demonstrate the computational efficiency of the hybrid algorithm against fine-scale simulations.
\end{abstract}

\begin{keyword}



Homogenization \sep Upscaling \sep Multiscale Modeling \sep Hybrid Models \sep Algorithmic Refinement

\end{keyword}

%% file: intro.tex
\section{Introduction} \label{sec_intro}

The premise of homogenization theory is to derive upcaled models of physical systems by exploiting the principle of separation of scales. Rigorous homogenization/coarse-graining methods \cite{Langlo1994-wq, Vasilyeva2019-wh,Das2005-gm,Frippiat2008-cl,Battiato2019-xk, Battiato2011-ad} provide \emph{a priori} guarantees, such as asymptotic bounds, on the deviation of the upscaled model's solution from that of the averaged fine-scale model (Figure~\ref{fig:homogenization-illu}(a) and Figure~\ref{fig:homogenization-illu}(b)). These error bounds remain valid over regions of space and durations of time where/when a set of {\it applicability conditions} hold. However, deriving the upscaled governing equations is tedious and often limits the adoption of the method. Recently, Pietryzk \emph{et al.} \cite{Pietrzyk2021-lu} developed a symbolic upscaling engine, \symbolica, capable of automatically deriving  homogenized partial differential equations (PDEs), as well as their applicability conditions, from  fine-scale PDEs, initial conditions (ICs) and boundary conditions (BCs). \symbolica is capable of handling complex multiscale, multi-physics, heterogeneous, and nonlinear systems of PDEs that surpass humans' ability to handle by cumbersome pen-and-paper analysis. 

Although \emph{a priori} error bounds provide a solid ground for identifying cost-accuracy tradeoffs where/when the applicability conditions are satisfied, the challenge is to deal with situations in which such conditions are violated. This has been shown to happen in presence of, e.g., large spatial gradients of the quantities of interest \cite{Battiato2011-ad, Pietrzyk2021-lu, Boso2013-ag}. % In fact, the most interesting phenomena occur when the applicability conditions are violated due to large gradients, material heterogeneity and aging.
To overcome these challenges, hybrid or algorithmic refinement formulations have been developed \cite{Battiato2011-zo,Yousefzadeh2017-yc}: fine-scale equations are solved only in a sub-domain in which applicability conditions are violated whereas the upscaled equations are solved everywhere else (Figure~\ref{fig:homogenization-illu}(c)). Hybrid simulations can be built through either intrusive or nonintrusive coupling conditions between the two sub-domains (fine-scale and {macro-scale}). Intrusive coupling methods are characterized by the existence of an ``overlapping" or ``handshake" region where both fine-scale and upscaled governing equations are concurrently solved \cite{Battiato2011-zo,Roubinet2013-fg,Pettersson2013-ou}. Intrusive coupling methods are typically more expensive due to the existence of the ``overlapping" or ``handshake" region, and more difficult to implement in legacy codes. In nonintrusive coupling methods, each subdomain is only solved with one set of governing equations, either fine-scale or upscaled, and the coupling between the subdomains is formulated exclusively in terms of boundary conditions \cite{Yousefzadeh2017-yc,Kadeethum2022-ag}. 

In this paper, we develop a predictive nonintrusive two-way coupled hybrid formulation. The computational domain is represented by two different materials, each of which is governed by a set of governing equations.  We apply this approach to the use-case of thermal runaway simulation in Li-ion batteries (LIBs) because they provide a sufficiently complex example with heterogeneity and nonlinearity in the source term. In Section~\ref{sec:govern-eqs}, we describe the fine-scale and upscaled governing equations and formulations of the coupling boundary conditions with two different approaches. In Section~\ref{sec:therm-res}, we evaluate the efficiency and accuracy of the proposed hybrid formulation for heat transfer in a battery pack as a demonstration.

%% file: metrics.tex
\subsection{DARPA Program Metrics} \label{sec_metrics}

The methods in this paper, developed with support and oversight of the DARPA Computable Models Disruption Opportunity \cite{DARPACOMPMods}, demonstrate several measurable advancements over the state-of-the-art. Here, we summarize them in terms of the relevant program metrics, i.e., modeling accuracy and numerical efficiency. As previously discussed, upscaling theory by multiple scale expansions ensures that the modeling error of coarse-grained approximations is \emph{a priori} bounded under appropriate dynamic conditions expressed {in terms} of dimensionless numbers. When such conditions are locally (in space and/or time violated), it is therefore important that any further strategy (numerical or analytical) that aims at coupling fine-scale models with their continuum-scale counterpart in the same simulation domain be bounded by the {aforementioned} upscaling error. In this regard, the accuracy of any proposed hybrid scheme can be directly assessed against such an a priori error. In Sections~\ref{sec:acc-hc} and~\ref{sec:xhc-acc}, we show that both coupling schemes satisfy the requested accuracy. An additional important metric is that the computation cost associated with the iterative coupling between fine- and coarse-scale models does not overcome the cost of full fine-scale simulations over the microscopic domain (here considered the benchmark for both accuracy and cost). In Section~\ref{sec:efficiency}, we provide both an extensive analysis of the cost-accuracy tradeoffs as well as guidelines for the efficient adoption of hybrid algorithms in large-scale domains.

\subsubsection{\bf Direct Computability:}
{Governing PDEs of heterogeneous multiscale systems are derived from first principles (e.g., conservation laws) and constitutive/material laws, applicable at length/time scales that are not ideal for computation, optimization or design, since they require resolving geometric features at scales much finer than the device scale. Upscaled (e.g., homogenized) models, on the other hand, are computationally more tractable  and, generally, more easily parametrizable since  their effective properties can be measured at the device scales. However, they are not directly computable when their applicability conditions  are violated because they loose accuracy and error guarantees can not be satisfied. Our approach enables bridging the gap between formulation and computation of multiscale systems where diffent models may be needed in the same simulation domain to ensure predictive accuracy while keeping computational costs in check: we achieve this by using the proper model  at the proper scales (fine-scale or upscaled), in different regions of space or intervals of time, and coupling them at the interfaces, while respecting the error bounds. Moreover, the non-intrusive nature of the coupling conditions allows one to use existing (e.g., off-the-shelf commercial or open-source) solvers, specialized for each regime/scale, without having to rewrite the code. By providing the information on the dimensionless numbers, the method can automatically determine the location of the breakdown region and compute the coupling locations.} 

\subsubsection{Accuracy and Efficiency:}
{As a result of the development of the hybrid non-intrusive two-way coupling approach, the errors are still bounded by the upscaling errors which can be determined \emph{a priori}. Additionally, the computational cost of hybrid simulations is significantly lower than that of fine-scale simulations 
%\deleted{with the existence of a breakdown region} 
(Section~\ref{sec:efficiency}). The performance metrics on computatibility and accuracy are summarized in Table~\ref{tab:darpa-metrics}.}

\begin{table}[H]
\centering
\caption{
  \label{tab:darpa-metrics} Summary of performance metrics.}
\begin{tabular}{p{0.15\linewidth}|p{0.25\linewidth}|p{0.25\linewidth}|p{0.25\linewidth}}
\hline \hline
Metric &  SOA & Initial & Ultimate Target  \\ \hline
Computability  & Fine-scale simulations over macroscopic domains & Hybrid coupling

For 80 battery cells and 2.5\% of volume solved by fine-scale model: speedup=3 
 & Hybrid coupling

For 800 battery cells and 2.5\% of volume solved by fine-scale model: speedup=15 \\ \hline
Conservation \& Accuracy  & Fine-scale simulations over macroscopic domains & Not Applicable.

Cannot guarantee error bounds, unless by fully resolving fine-scale models.

 & A priori asymptotic error guarantees of order $\epsilon$. \\ \hline \hline
\end{tabular}%
\end{table}

%% file: equations.tex
\section{Thermal Runaway in Li-Ion Batteries}
\label{sec:govern-eqs}

\begin{figure}
    \centering
    \includegraphics[width=0.9\textwidth]{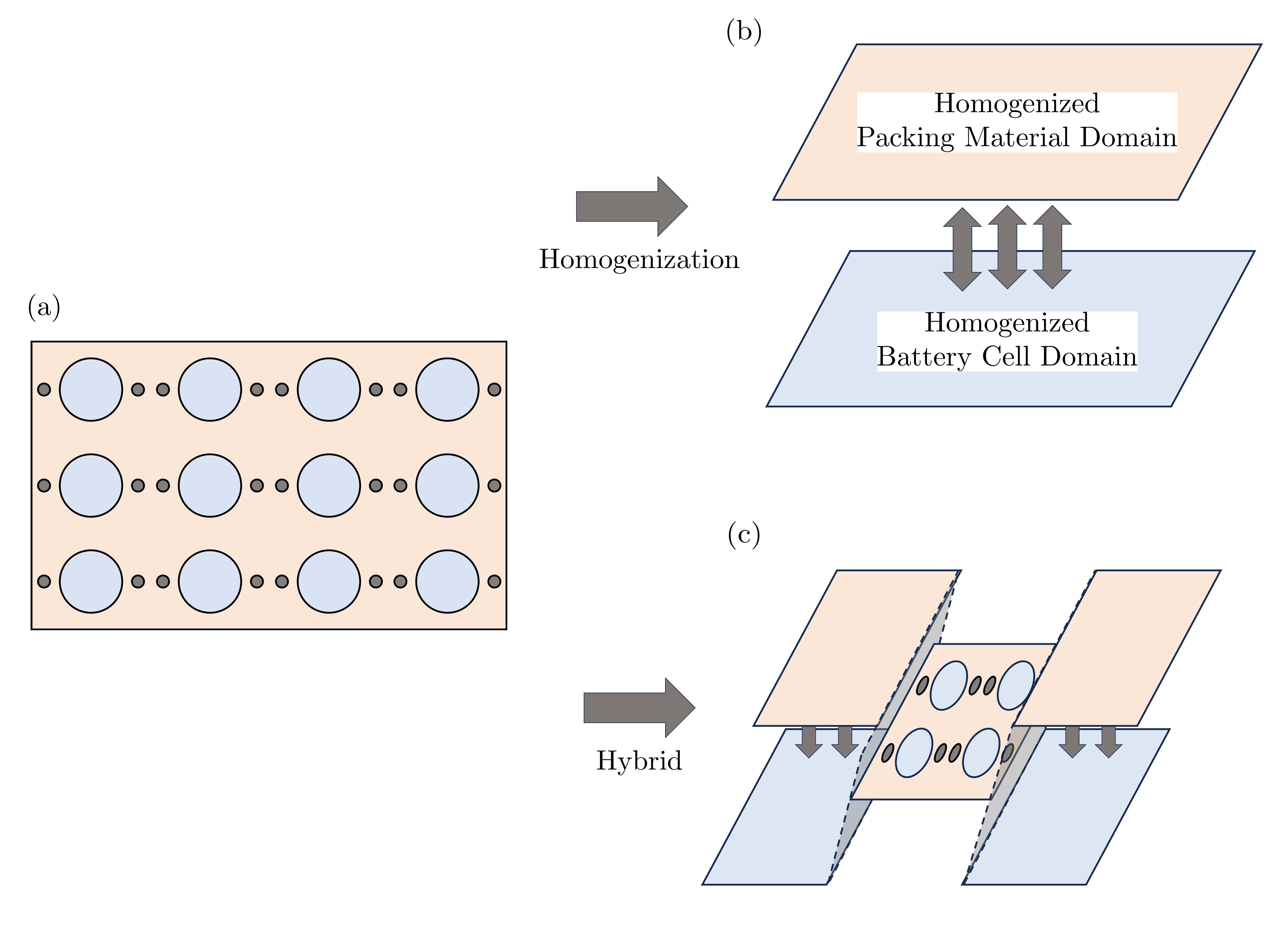}
    \caption{{Illustration of (a) a fine-scale battery pack domain, (b) homogeization of the fine-scale battery pack domain into two upscaled packing material and battery cell domains, and (c) a hybrid domain consists of both fine-scale and homogenized domains.}}
    \label{fig:homogenization-illu}
\end{figure}

The onset of thermal runaway in a cell due to mechanical, thermal, and electric abuse can compromise the entire battery pack and lead to an explosion \cite{Feng2018-ix}. Understanding heat transfer in these systems at the relevant scales is critical to optimize design and operation. Nevertheless, the development of accurate heat transfer models in battery systems, ranging from the sub-electrode to the battery pack scale, represents a formidable multiscale task because of the complex interactions between the processes at each scale, such as heat generation at the electrode scale and thermal runaway at the packing/module scale. 

One common approach to modeling thermal runaway in battery cells is to develop spatially independent models with experimentally determined and calibrated parameters \cite{Ren2018-zj}. However, such parameters, calibrated on the basis of lab-scale systems, can have significant uncertainties and tend to vary from system to system. Therefore, the applicability of these models to large-scale systems is questionable. 

The alternative is to perform fine-scale simulations that resolve the governing partial differential equations in a defined computational domain \cite{Kong2021-sb, Guo2019-td}. With appropriate mesh resolution and computational methods, heat transfer in the battery pack can be accurately simulated. However, these high-fidelity simulations generally require significant computational costs \cite{Yousefzadeh2023-rg}, and are unlikely to be used as predictive tools for large battery packs that contain thousands of cells.

Another approach is to develop upscaled models from fine-scale equations based on homogenization/coarse-graining theory. Recently, Pietrzyk \emph{et al.} \cite{Pietrzyk2023-ou} have generated upscaled equations using the automated upscaling engine, \symbolica, to model heat transfer in battery packs. Yet, when certain applicability conditions are invalidated due to, e.g., manufacturing defects and cell aging, the prediction from upscaled models can deviate significantly from the average fine-scale behavior. In the following, we develop a hybrid formulation that couples fine-scale equations with upscaled equations when applicatibility conditions of the latter are violated within a small portion of the computational domain. In Section~\ref{sec:sub-govens}, we present the governing equations at the fine- and continuum-scales for the thermal runaway problem.

\subsection{Problem description and governing equations}
\label{sec:sub-govens}
\begin{figure}
\centerline{
 {\includegraphics[width=0.6\textwidth]{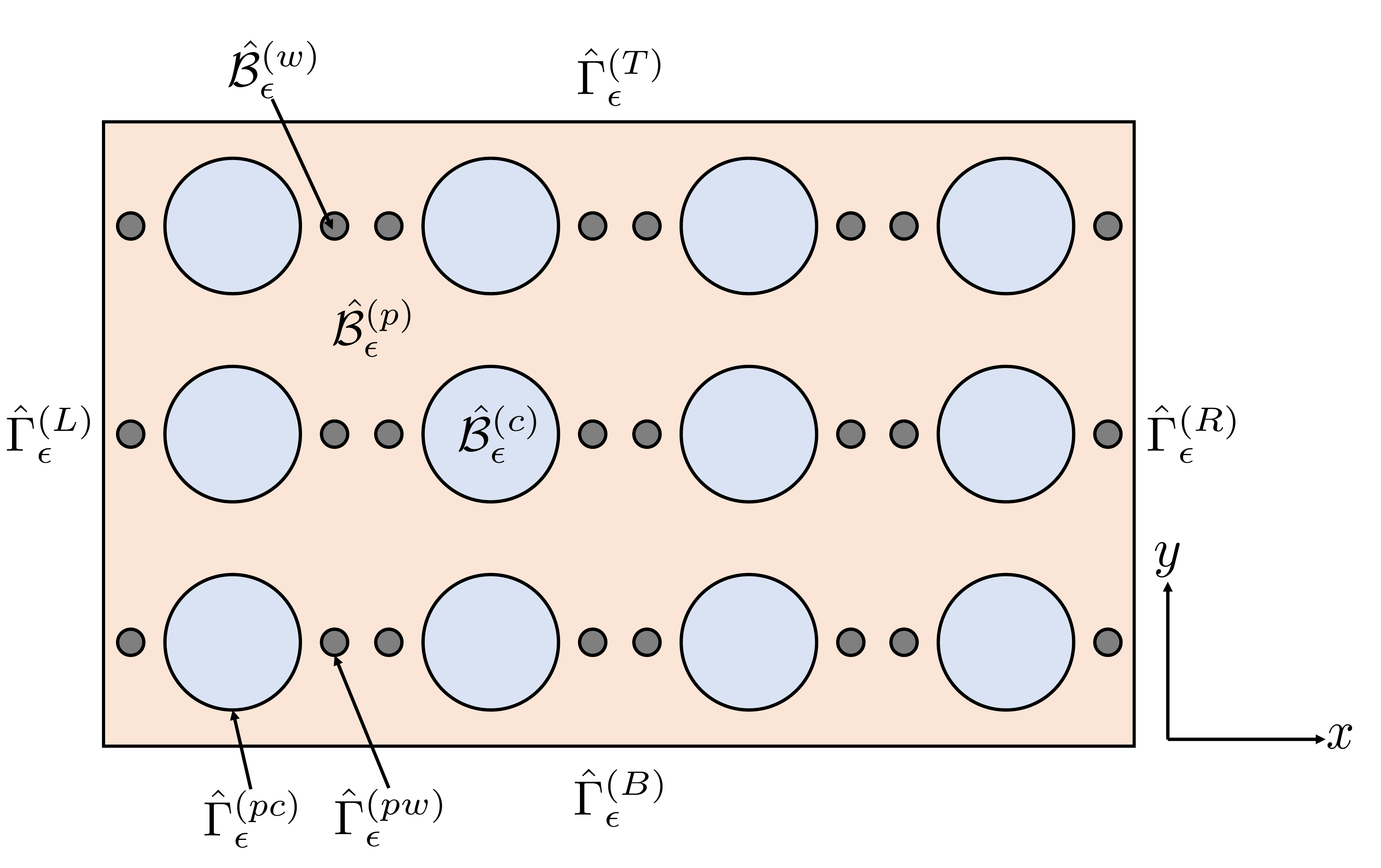}}}
\caption{Example of a fine-scale two-dimensional battery pack.}
\label{fig:domain_eg}
\end{figure}

We consider heat transfer within a two-dimensional battery pack $\hat{\Omega}_\epsilon \subset \mathbb{R}^2$ (Figure~\ref{fig:domain_eg}) that is bounded by packing edges $\hat{\Gamma}_{\epsilon} \subset \mathbb{R}^2$. The edges of the pack include the top, bottom, left, and right and are referred to as $\hat{\Gamma}_{\epsilon}^{\left(T\right)} \subset \hat{\Gamma}_{\epsilon}$, $\hat{\Gamma}_{\epsilon}^{\left(B\right)} \subset \hat{\Gamma}_{\epsilon}$, $\hat{\Gamma}_{\epsilon}^{\left(L\right)} \subset \hat{\Gamma}_{\epsilon}$ and $\hat{\Gamma}_{\epsilon}^{\left(R\right)} \subset \hat{\Gamma}_{\epsilon}$, respectively. The battery pack consists of three distinct domains: battery cells $\hat{\mathcal{B}}_{\epsilon}^{\left(c\right)} \subset \hat{\Omega}_{\epsilon}$, packing material $\hat{\mathcal{B}}_{\epsilon}^{\left(p\right)} \subset \hat{\Omega}_{\epsilon}$ and cooling water pipes $\hat{\mathcal{B}}_{\epsilon}^{\left(w\right)} \subset \hat{\Omega}_{\epsilon}$ (Figure~\ref{fig:domain_eg}); $\hat{\Gamma}_\epsilon^{\left(pc\right)}$ denotes the interface between packing material and battery cell, while $\hat{\Gamma}_\epsilon^{\left(pw\right)}$ denotes the interface between packing material and cooling water pipes. Both packing materials and cooling water pipes are used as heat sinks to absorb the heat generated by the battery cells. Following the notation of Pietrzyk et al.~\cite{Pietrzyk2023-ou}, we use the ``hat'' and subscript $\epsilon$ to denote dimensional and fine-scale quantities, respectively. The governing equations of heat transfer in the packing material and battery cells are modeled as 
\begin{subequations}
\label{eq:Heat_Transfer_Eq}
\begin{align}
&\derive{\left(\hat{\rho}^{\left(p\right)}\hat{C}^{\left(p\right)}\hat{T}_{\epsilon}^{\left(p\right)}\right)}{\hat{t}} = \hat{\nabla}\bm{\cdot}\left(\hat{k}^{\left(p\right)}\hat{\nabla}\hat{T}_{\epsilon}^{\left(p\right)}\right) \quad \text{for } \hat{\mathbf{x}} \in \hat{\mathcal{B}}_{\epsilon}^{\left(p\right)}, \label{eq:Heat_Transfer_Eq_Packing} \\
% &\derive{\left(\hat{\rho}^{\left(c\right)}\hat{C}^{\left(c\right)}\hat{T}_{\epsilon}^{\left(c\right)}\right)}{\hat{t}} = \hat{\nabla}\bm{\cdot}\left(\hat{k}^{\left(c\right)}\hat{\nabla}\hat{T}_{\epsilon}^{\left(c\right)}\right) + \hat{\Pi}\left(\hat{t}, \hat{\mathbf{x}}\right) \quad \text{for } \hat{\mathbf{x}} \in \hat{\mathcal{B}}_{\epsilon}^{\left(c\right)}, 
&\derive{\left(\hat{\rho}^{\left(c\right)}\hat{C}^{\left(c\right)}\hat{T}_{\epsilon}^{\left(c\right)}\right)}{\hat{t}} = \hat{\nabla}\bm{\cdot}\left(\hat{k}^{\left(c\right)}\hat{\nabla}\hat{T}_{\epsilon}^{\left(c\right)}\right) + {\hat{\Pi}} \quad \text{for } \hat{\mathbf{x}} \in \hat{\mathcal{B}}_{\epsilon}^{\left(c\right)}, \label{eq:Heat_Transfer_Eq_Cell} 
\end{align}

\noindent respectively, subject to boundary conditions
\begin{align}
 &-\n_{\epsilon}^{\left(p\right)} \bm{\cdot} \hat{k}^{\left(p\right)}\hat{\nabla}\hat{T}_{\epsilon}^{\left(p\right)} = \hat{U}^{\left(pc\right)}\left(\hat{T}_{\epsilon}^{\left(p\right)} - \hat{T}_{\epsilon}^{\left(c\right)}\right) \quad \text{for } \hat{\mathbf{x}} \in \hat{\Gamma}_{\epsilon}^{\left(pc\right)}, \label{eq:BC_pc_packing}\\
% &-\n_{\epsilon}^{\left(p\right)} \bm{\cdot} \hat{k}^{\left(p\right)}\hat{\nabla}\hat{T}_{\epsilon}^{\left(p\right)} = \hat{q}_{\epsilon}^{\left(pw\right)}\left(\hat{t}, \hat{\mathbf{x}}\right) \quad \text{for } \hat{\mathbf{x}} \in \hat{\Gamma}_{\epsilon}^{\left(pw\right)}, \label{eq:BC_pw}\\
&-\n_{\epsilon}^{\left(p\right)} \bm{\cdot} \hat{k}^{\left(p\right)}\hat{\nabla}\hat{T}_{\epsilon}^{\left(p\right)} = {\hat{q}_{\epsilon}^{\left(pw\right)}} \quad \text{for } \hat{\mathbf{x}} \in \hat{\Gamma}_{\epsilon}^{\left(pw\right)}, \label{eq:BC_pw}\\
&-\n_{\epsilon}^{\left(c\right)} \bm{\cdot} \hat{k}^{\left(c\right)}\hat{\nabla}\hat{T}_{\epsilon}^{\left(c\right)} = \hat{U}^{\left(pc\right)}\left(\hat{T}_{\epsilon}^{\left(c\right)} - \hat{T}_{\epsilon}^{\left(p\right)}\right) \quad \text{for } \hat{\mathbf{x}} \in \hat{\Gamma}_{\epsilon}^{\left(pc\right)}, \label{eq:BC_pc_cell}  
\end{align}
\end{subequations}
\noindent where $i=p$ or $c$ refers to packing material or the battery cell, respectively,  $\hat{\rho}^{\left(i\right)}$ [\si{ML\tothe{-3}}] is the density, $\hat{C}^{\left(i\right)}$ [\si{L\tothe{2}T\tothe{-2}\Theta\tothe{-1}}] is the heat capacity, {$\hat{T}_{\epsilon}^{\left(i\right)}  \equiv \hat{T}_{\epsilon}^{\left(i\right)}\left(\hat{t}, \hat{\mathbf{x}}\right) [\si{\Theta}]$}  is the temperature at time $\hat{t} > 0$ and location $\hat{\mathbf{x}} \in \hat{\mathcal{B}}_{\epsilon}^{\left(i\right)}$, $\hat{k}^{\left(i\right)}$ [\si{MLT\tothe{-3}\Theta\tothe{-1}}] is the thermal conductivity, $\n_{\epsilon}^{\left(i\right)} \equiv \n_{\epsilon}^{\left(i\right)}\left(\hat{\mathbf{x}}\right)$ is the normal vector to the interfaces pointing away from the domain, $\hat{U}^{\left(pc\right)}$ [\si{MT\tothe{-3}\Theta\tothe{-1}}] is the total heat transfer coefficient between the packing material and battery cells,  $\hat{q}_{\epsilon}^{\left(pw\right)}\left(\hat{t}, \hat{\mathbf{x}}\right)$ [\si{MT\tothe{-3}}]  is a power flux between the packing material and the cooling water pipes, and $\hat{\Pi}(\hat{t}, \hat{\mathbf{x}})$ [\si{ML\tothe{-1}T\tothe{-3}}] is a power flux source term. {Equation~\eqref{eq:Heat_Transfer_Eq} still holds with temperature-dependent material properties. For this study, we assume that the material properties are constant.}

To model the battery cell heat generation, the power flux source term $\hat{\Pi}(\hat{t}, \hat{\mathbf{x}})$ is approximated using polynomial of $\hat{T}_{\epsilon}^{\left(c\right)}${~\cite{Pietrzyk2023-ou} where}
\begin{subequations}\label{PI}
\begin{align}
    & {\hat{\Pi}\left(\hat{t}, \hat{\x}\right) \equiv \hat{\Pi}\left(\hat{T}_{\epsilon}^{\left(c\right)}, \hat{\x}\right) \equiv \sum_{k = 0}^{\infty} \hat{a}_k\left(\hat{\x}\right) \left(\hat{T}_{\epsilon}^{\left(c\right)}\right)^k,}
\end{align}
such that 
\begin{align}
\hat{\Pi}\left(\hat{T}_{\epsilon}^{\left(c\right)}, \hat{\mathbf{x}}\right) &= \hat{\Pi}_{\text{base}}\left(\hat{\mathbf{x}}\right) + \frac{1}{2}\left\{\text{Erf}\left[C_1\left(2\frac{\hat{T}_a - \hat{T}_{\epsilon}^{\left(c\right)} + \hat{T}_{ref}}{\hat{T}_{s1}} + 1\right)\right] + 1\right\}\left(\hat{\Pi}_{\text{burn}} - \hat{\Pi}_{\text{base}}\left(\hat{\mathbf{x}}\right)\right) \nonumber \\   
&- \frac{1}{2}\left\{\text{Erf}\left[C_2\left(2\frac{\hat{T}_{max} - \hat{T}_{\epsilon}^{\left(c\right)} + \hat{T}_{ref}}{\hat{T}_{s2}} - 1\right)\right] + 1\right\}\hat{\Pi}_{\text{burn}}, \label{eq:new_pi_term1}
\end{align}
\noindent and
\begin{align}
&\hat{T}_{max} = \hat{T}_{Max} - \hat{T}_{ref}, \label{eq:new_pi_term2}\\ 
&\hat{T}_{Max} = \hat{T}_{ref} + \hat{T}_{a} + \hat{T}_{s1} + \hat{T}_{b} + \hat{T}_{s2}, \label{eq:new_pi_term2_1}\\
&{C_1 = \text{Erf}^{-1}\left(2\epsilon_{s1} - 1\right),} \label{eq:new_pi_term3}\\
&{C_2 = \text{Erf}^{-1}\left(2\epsilon_{s2} - 1\right)}. \label{eq:new_pi_term3_2}  
\end{align}
\end{subequations}
\noindent {In Eq.~\eqref{PI}, $\hat{a}_k(\hat{\x})$ are the coefficients of the polynomial}, $\hat{\Pi}_{\text{base}}(\hat{\mathbf{x}})$ and $\hat{\Pi}_{\text{burn}}$ are the \textit{base} and \textit{burn} power flux values, respectively, $\hat{T}_{ref}$~[\si{\Theta}] is the reference temperature, $\hat{T}_a$~[\si{\Theta}] and $\hat{T}_b$~[\si{\Theta}] are temperature ranges {over which} $\hat{\Pi}(\hat{T}_{\epsilon}^{\left(c\right)},\hat{\mathbf{x}}) = \hat{\Pi}_{\text{base}}(\hat{\mathbf{x}})$ and $\hat{\Pi}(\hat{T}_{\epsilon}^{\left(c\right)}, \hat{\mathbf{x}}) = \hat{\Pi}_{\text{burn}}$, respectively, $\hat{T}_{s1}$ and $\hat{T}_{s2}$ are the temperature ranges for which $\hat{\Pi}(\hat{T}_{\epsilon}^{\left(c\right)}, \hat{\mathbf{x}})$ transitions from $\hat{\Pi}_{\text{base}}(\hat{\mathbf{x}})$ to $\hat{\Pi}_{\text{burn}}$ and from $\hat{\Pi}_{\text{burn}}$ to $0$, respectively, and {$\epsilon_{s1} = 0.0005$ and $\epsilon_{s2} = 0.0005$ are dimensionless parameters associated with} smoothness of the error functions. {A plot for $\hat \Pi$ can be found in Figure 2 of Pietrzyk et al.~\cite{Pietrzyk2023-ou}. Additionally,} the detailed formulation and validation of the source term can be found in Pietrzyk et al.~\cite{Pietrzyk2023-ou}.

\subsection{Dimensionless fine-scale governing equations}

{In order to apply homogenization theory, we followed the nondimensionalization procedures in Pietrzyk et al.~\cite{Pietrzyk2023-ou} to derive a set of dimensionless numbers that provides insight into the physical mechanisms within a system. We define the reference scales as} 

%\textcolor{red}{[It should be $\hat{\mathbf{x}} = \hat{\mathcal{L}}\mathbf{x}$. Also the comma was left in the fraction]}

\begin{gather}
\hat{T}_{\epsilon}^{\left(p\right)} = \hat{T}_{max}T_{\epsilon}^{\left(p\right)} + \hat{T}_{ref}, \quad \hat{T}_{\epsilon}^{\left(c\right)} = \hat{T}_{max}T_{\epsilon}^{\left(c\right)} + \hat{T}_{ref}, \quad \hat{k}^{\left(p\right)} = \hat{\mathcal{K}}^{\left(p\right)} k^{\left(p\right)}, \nonumber \\ 
\hat{k}^{\left(c\right)} = \hat{\mathcal{K}}^{\left(c\right)} k^{\left(c\right)}, \quad \hat{\nabla} = \frac{1}{\hat{\mathcal{L}}} \nabla, \quad \hat{t} = \frac{\hat{\rho}^{\left(p\right)}\hat{C}^{\left(p\right)}\hat{\mathcal{L}}^2}{\hat{\mathcal{K}}^{\left(p\right)}} t, \quad {\hat{\mathbf{x}} = \hat{\mathcal{L}}\mathbf{x}},  \nonumber \\ \hat{\Pi}\left(\hat{T}_{\epsilon}^{\left(c\right)}, \hat{\mathbf{x}}\right) = \hat{\Pi}_{\text{burn}} \Pi\left(T_{\epsilon}^{\left(c\right)}, \mathbf{x}\right), \quad \hat{q}_{\epsilon}^{\left(pw\right)}\left(\hat{t}, \hat{\mathbf{x}}\right) = \hat{Q}^{\left(pw\right)} q_{\epsilon}^{\left(pw\right)}\left(t, \mathbf{x}\right), \label{eq:Scales}
\end{gather}
\noindent where {$\hat{\mathcal{L}}$ is the size of the macroscopic domain,} $\hat{\mathcal{K}}^{\left(i\right)}$ is the reference scale of the thermal conductivity, and $\hat{Q}^{\left(pw\right)}$ is the reference scale of the power flux sink term at $\hat{\Gamma}_{\epsilon}^{\left(pw\right)}$. By scaling the equations with the reference scales, the nondimensional fine-scale governing equations are obtained as 
\begin{subequations}
\begin{align}
&\derive{T_{\epsilon}^{\left(p\right)}}{t} = \nabla \bm{\cdot} \left(k^{\left(p\right)} \nabla T_{\epsilon}^{\left(p\right)}\right) \quad \text{for } \mathbf{x} \in \mathcal{B}_{\epsilon}^{\left(p\right)}, \label{eq:Heat_Transfer_Eq_Packing_dimless} \\
&\derive{T_{\epsilon}^{\left(c\right)}}{t} = \left(\varrho \bm{\cdot} \varsigma\right) \nabla \bm{\cdot} \left(k^{\left(c\right)} \nabla T_{\epsilon}^{\left(c\right)}\right) + \left(\varrho \bm{\cdot} \mathcal{R}\right) {\Pi} \quad \text{for } \mathbf{x} \in \mathcal{B}_{\epsilon}^{\left(c\right)}, \label{eq:Heat_Transfer_Eq_Cell_dimless}
\end{align}
\label{eq:pore_goven_eqs}
\noindent subject to boundary conditions 
\begin{align}
&-\n_{\epsilon}^{\left(p\right)} \bm{\cdot} k^{\left(p\right)} \nabla T_{\epsilon}^{\left(p\right)} = \text{Bi}^{\left(p\right)}\left(T_{\epsilon}^{\left(p\right)} - T_{\epsilon}^{\left(c\right)}\right) \quad \text{for } \mathbf{x} \in \Gamma_{\epsilon}^{\left(pc\right)}, \label{eq:BC_pc_packing_dimless} \\
&-\n_{\epsilon}^{\left(p\right)} \bm{\cdot} k^{\left(p\right)}\nabla T_{\epsilon}^{\left(p\right)} = \mathcal{Q} \bm{\cdot} {q_{\epsilon}^{\left(pw\right)}} \quad \text{for } \mathbf{x} \in \Gamma_{\epsilon}^{\left(pw\right)}, \label{eq:BC_pw_dimless} \\
&-\n_{\epsilon}^{\left(c\right)} \bm{\cdot} k^{\left(c\right)} \nabla T_{\epsilon}^{\left(c\right)} = \text{Bi}^{\left(c\right)}\left(T_{\epsilon}^{\left(c\right)} - T_{\epsilon}^{\left(p\right)}\right) \quad \text{for } \mathbf{x} \in \Gamma_{\epsilon}^{\left(pc\right)}. \label{eq:BC_pc_cell_dimless} 
\end{align}
\end{subequations}
The dimensionless power flux source term is defined as 
\begin{subequations}
\begin{align}
\Pi\left(T_{\epsilon}^{\left(c\right)}, \mathbf{x}\right) &= \Pi_{\text{base}}\left(\mathbf{x}\right) + \frac{1}{2}\left\{\text{Erf}\left[A_1 T_{\epsilon}^{\left(c\right)} + B_1\right] + 1\right\}\left(1 - \Pi_{\text{base}}\left(\mathbf{x}\right)\right) \nonumber \\ 
& - \frac{1}{2}\left\{\text{Erf}\left[A_2 T_{\epsilon}^{\left(c\right)} + B_2\right] + 1\right\},\label{eq:new_pi_term_dimless}
\end{align}

\noindent where

\neweq{new_pi_term_dimless_coefs}{A_1 = -2C_1\frac{\hat{T}_{max}}{\hat{T}_{s1}}, \quad B_1 = 2C_1\frac{\hat{T}_{a}}{\hat{T}_{s1}} + C_1, \quad A_2 = -2C_2\frac{\hat{T}_{max}}{\hat{T}_{s2}}, \quad B_2 = 2C_2\frac{\hat{T}_{max}}{\hat{T}_{s2}} - C_2.}
\end{subequations}

\noindent Scaling the governing equations results in six dimensionless numbers $\text{Bi}^{\left(p\right)}$, $\mathcal{Q}$, $\varrho$, $\varsigma$, $\text{Bi}^{\left(c\right)}$, and $\mathcal{R}$ such that

\neweqgat{eq:dimless_groups}{\text{Bi}^{\left(p\right)} = \frac{\hat{U}^{\left(pc\right)}\hat{\mathcal{L}}}{\hat{\mathcal{K}}^{\left(p\right)}}, \quad \mathcal{Q} = \frac{\hat{Q}^{\left(pw\right)}\hat{\mathcal{L}}}{\hat{T}_{max}\hat{\mathcal{K}}^{\left(p\right)}}, \quad \varrho = \frac{\hat{\rho}^{\left(p\right)}\hat{C}^{\left(p\right)}}{\hat{\rho}^{\left(c\right)}\hat{C}^{\left(c\right)}}, \\
\varsigma = \frac{\hat{\mathcal{K}}^{\left(c\right)}}{\hat{\mathcal{K}}^{\left(p\right)}},  \quad \text{Bi}^{\left(c\right)} = \frac{\text{Bi}^{\left(p\right)}}{\varsigma}, \quad \mathcal{R} = \frac{\hat{\Pi}_{\text{burn}}\hat{\mathcal{L}}^2}{\hat{T}_{max}\hat{\mathcal{K}}^{\left(p\right)}}.}

\subsection{Unit-cell and domain formulation}
\begin{figure}
\centerline{
 {\includegraphics[width=.5\linewidth]{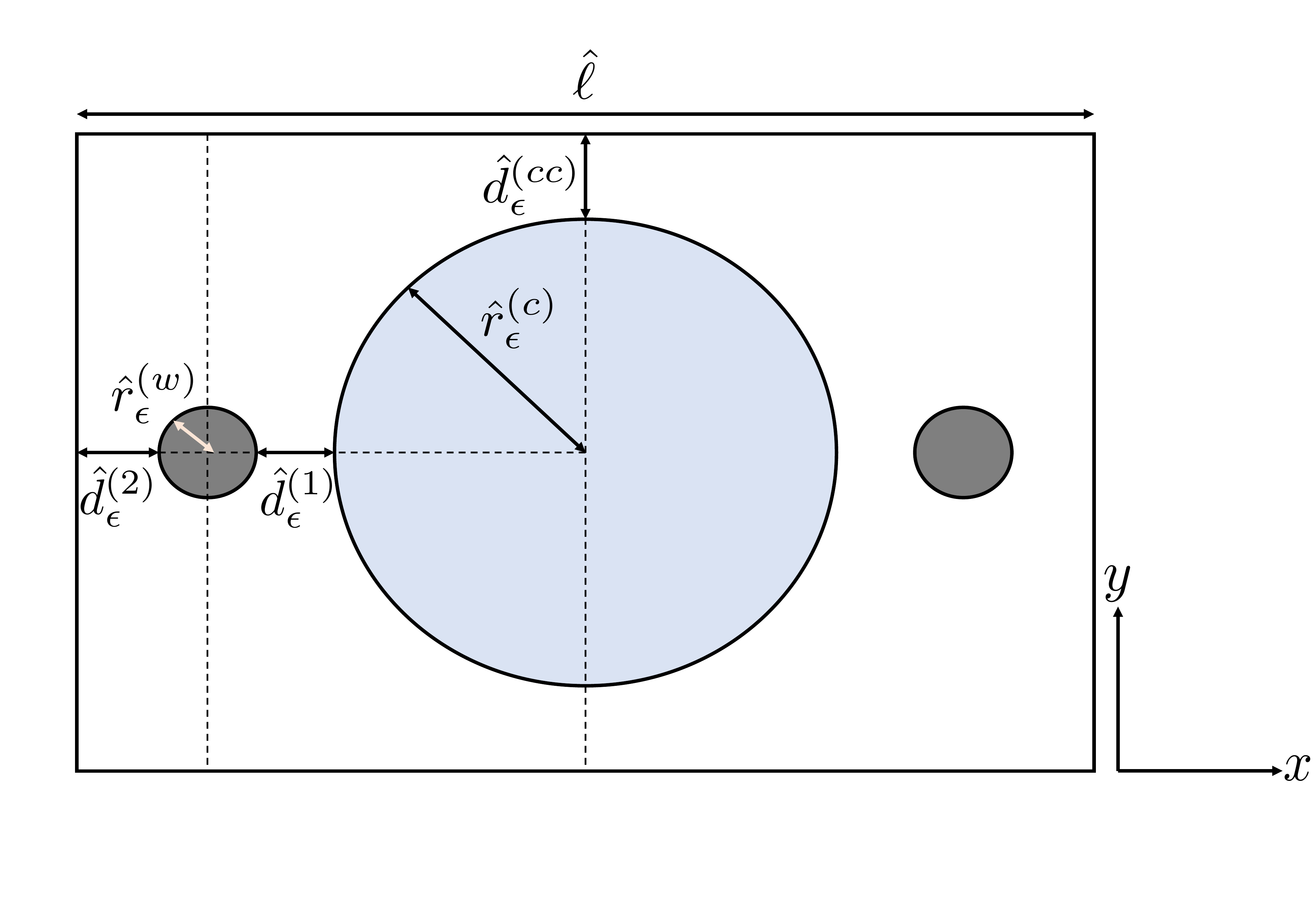}}}
\caption{Unit cell geometry with defined dimensional parameters {in the $x$ and $y$ directions.}}
\label{fig:unit-cell-geom}
\end{figure}

In each battery pack, the battery cells can be divided into $n$ regions such that $\mathcal{B}_\epsilon^{\left(c,n\right)} \subseteq \mathcal{B}_\epsilon^{\left(c\right)}$ $n \in \mathbb{Z}^{+}$, $n \leq N^{\left(c\right)}$, and $N^{\left(c\right)}$ is the number of battery cells in the battery pack. The geometry of the two-dimensional battery pack is defined by the unit cell, the number of unit cells in the $x$ and $y$ directions $N^{\left(c\right)}_{x}$ and $N^{\left(c\right)}_{y}$, respectively. The unit cell geometry (Figure~\ref{fig:unit-cell-geom}) is then defined by the distance between the battery cell and the boundary of the unit cell $\hat{d}_\epsilon^{\left(cc\right)}$, the radius of the battery cell $\hat{r}_\epsilon^{\left(c\right)}$, the radius of the cooling water pipe $\hat{r}_\epsilon^{\left(w\right)}$ and the distances between the battery cell and the cooling pipe $\hat{d}_\epsilon^{\left(1\right)}$ and $\hat{d}_\epsilon^{\left(2\right)}$. The length of the unit cell is computed as 
\begin{align}
\hat{\ell} &= 2\left( \hat{d}_\epsilon^{\left(1\right)} + \hat{d}_\epsilon^{\left(2\right)} + \hat{r}_\epsilon^{\left(c\right)} + \hat{r}_\epsilon^{\left(w\right)} \right),
\end{align}
the aspect ratio of the unit cell $a$ is defined as 
\begin{align}
a &= \ddfrac{2\left(\hat{d}_\epsilon^{\left(cc\right)} + \hat{r}_\epsilon^{\left(c\right)}\right)}{\hat{\ell}}.
\end{align}
With the defined unit cell geometry, {the length $\hat{\mathcal{L}}_{x}$ and the width $\hat{\mathcal{L}}_{y}$ of the battery pack are defined as} 
\begin{align}
\hat{\mathcal{L}}_{x} &= N_{x}^{\left(c\right)}\hat{\ell},\\
\hat{\mathcal{L}}_{y} &= N_{y}^{\left(c\right)}\hat{\ell}a.
\end{align}
{The main objective of this study is to demonstrate the capability and accuracy of the approach. As a result, the arrangement of cooling pipes and battery cells may not be realisitc since their location is not optimized.}

\subsection{Dimensionless upscaled governing equations}

Based on homogenization theory, the upscaled governing equations are derived from the fine-scale governing equations using the automated pipeline, \symbolica~\cite{Pietrzyk2023-ou} by applying averaging operators as

\neweqgat{eq:ave-op}{\langle \;\; \bm{\cdot}^{\left(i\right)} \rangle_{Y} \equiv \frac{1}{|Y|} \int_{\mathcal{B}^{\left(i\right)}}\left( \;\; \bm{\cdot}^{\left(i\right)}\right)\;d\xiv, \quad \langle \bm{\cdot} \rangle_{\mathcal{B}^{\left(i\right)}} \equiv \frac{1}{|\mathcal{B}^{\left(i\right)}|} \int_{\mathcal{B}^{\left(i\right)}}\left(\bm{\cdot}\right)\;d\xiv, \quad \langle \bm{\cdot} \rangle_{\Gamma^{\left(j\right)}} \equiv \frac{1}{|\Gamma^{\left(j\right)}|} \int_{\Gamma^{\left(j\right)}}\left(\bm{\cdot}\right)\;d\xiv,}

\noindent respectively, where $i = p$, $c$, or $w$ refer to packing material, battery cell and cooling water pipe, respectively; $j = pc$ or $pw$ refers to the interface between packing material and battery cell, and the interface between packing material and cooling water pipes, respectively; $|Y|$ is the area of the unit-cell domain; $|\mathcal{B}^{\left(i\right)}|$ is the area of the region $\mathcal{B}^{\left(i\right)}$ in the unit-cell, and $|\Gamma^{\left(j\right)}|$ is the length of interface $\Gamma^{\left(j\right)}$ in the unit-cell. {While the second averaging operator in Equation~\eqref{eq:ave-op} is typically used in homogenization, the first averaging operator ensures a better comparison between fine-scale and upscaled results. Hence, for completeness, we report both.}
% $\langle \bm{\cdot} \rangle_Y$ such that 
% \begin{align}
% \label{eq:ave-op}
% \langle \bm{\cdot} \rangle_Y \equiv \ddfrac{1}{\abs{Y}} \int_Y \left( \bm{\cdot} \right) \ \dd\boldsymbol{\xi} = \ddfrac{1}{\abs{Y}} \left[ \int_{\mathcal{B}^p} \left(\cdot\right) \ \dd\boldsymbol{\xi} + \int_{\mathcal{B}^c} \left(\cdot\right) \ \dd\boldsymbol{\xi}\right],
% \end{align}
% where $\abs{Y}$ refers to the volume of the unit cell $Y$.

{

The main objective of deriving upscaled equations is to obtain average packing material and cell temperatures.} The accuracy of the derived upscaled equations has been thoroughly validated by Pietrzyk et al.~\cite{Pietrzyk2023-ou}. The upscaled governing equations of the average temperatures $\ensmean{T^{p}}_Y$ and $\ensmean{T^{c}}_Y$ are {derived by applying the averaging operator 
 (equation~\eqref{eq:ave-op}) such that}
 
\begin{subequations}\label{upscaled}
\label{eq:upscaled_goven_eqs}
\neweq{Nc_2_homo_eq_1}{\phi^{\left(p\right)}\derive{\langle T^{\left(p\right)} \rangle_{Y}}{t} + \U^{\left(p\right)} \bm{\cdot} \nabla_{\mathbf{x}} \langle T^{\left(p\right)} \rangle_{Y} - \V^{\left(p\right)} \bm{\cdot} \nabla_{\mathbf{x}} \langle T^{\left(c\right)} \rangle_{Y} - \nabla_{\mathbf{x}} \bm{\cdot} \left(\K^{\left(p\right)} \bm{\cdot} \nabla_{\mathbf{x}} \langle T^{\left(p\right)} \rangle_{Y}\right)                   \\
= -R_1^{\left(p\right)}\langle T^{\left(p\right)} \rangle_{Y} + R_2^{\left(p\right)}\langle T^{\left(c\right)} \rangle_{Y} - R_3^{\left(p\right)}q^{\left(pw\right)}\left(t, \mathbf{x}\right) + \R_4^{\left(p\right)} \bm{\cdot} \nabla_{\mathbf{x}} {q^{\left(pw\right)}},}

\neweq{Nc_2_homo_eq_2}{\phi^{\left(c\right)}\derive{\langle T^{\left(c\right)} \rangle_{Y}}{t} + \U^{\left(c\right)} \bm{\cdot} \nabla_{\mathbf{x}} \langle T^{\left(c\right)} \rangle_{Y} - \V^{\left(c\right)} \bm{\cdot} \nabla_{\mathbf{x}} \langle T^{\left(p\right)} \rangle_{Y} - \nabla_{\mathbf{x}} \bm{\cdot} \left(\K^{\left(c\right)} \bm{\cdot} \nabla_{\mathbf{x}} \langle T^{\left(c\right)} \rangle_{Y}\right)               \\
= R_1^{\left(c\right)} \langle T^{\left(p\right)} \rangle_{Y} - R_2^{\left(c\right)} \langle T^{\left(c\right)} \rangle_{Y} + R_3^{\left(c\right)} q^{\left(pw\right)}\left(t, \mathbf{x}\right) + R_4^{\left(c\right)} {\overline{\Pi}},}

\neweq{eq:upscale_power_flux_source}{{\overline{\Pi}} = \Pi_{\text{base}} + \frac{1}{2}\left\{\text{Erf}\left[\frac{A_1}{\phi^{\left(c\right)}} \langle T^{\left(c\right)} \rangle_{Y} + B_1\right] + 1\right\}\left(1 - \Pi_{\text{base}}\right)                   \\
- \frac{1}{2}\left\{\text{Erf}\left[\frac{A_2}{\phi^{\left(c\right)}} \langle T^{\left(c\right)} \rangle_{Y} + B_2\right] + 1\right\}.}
\end{subequations}

\noindent where

\begin{subequations}\label{effectivecoefficients}
\neweq{}{\U^{\left(p\right)} = \phi^{\left(p\right)}\frac{\text{Bi}^{\left(p\right)}}{|\mathcal{B}^{\left(p\right)}|} |\Gamma^{\left(pc\right)}| \langle \chiv^{\left(p\right)\left[3\right]} \rangle_{\Gamma^{\left(pc\right)}} - k^{\left(p\right)}\langle\nabla_{\xiv}\chi^{\left(p\right)\left[2\right]}\rangle_{Y},}
\neweq{}{\V^{\left(p\right)} = \frac{\phi^{\left(p\right)}}{\phi^{\left(c\right)}}\left[\phi^{\left(p\right)}\frac{\text{Bi}^{\left(p\right)}}{|\mathcal{B}^{\left(p\right)}|} |\Gamma^{\left(pc\right)}| \langle \chiv^{\left(c\right)\left[2\right]} \rangle_{\Gamma^{\left(pc\right)}} - k^{\left(p\right)}\langle \nabla_{\xiv}\chi^{\left(p\right)\left[2\right]} \rangle_{Y}\right],}
\neweq{}{\K^{\left(p\right)} = k^{\left(p\right)}\left[\phi^{\left(p\right)}\I + \langle \nabla_{\xiv}\chiv^{\left(p\right)\left[3\right]} \rangle_{Y}\right],}
\neweq{}{R_1^{\left(p\right)} = \phi^{\left(p\right)}\frac{\text{Bi}^{\left(p\right)}}{|\mathcal{B}^{\left(p\right)}|}|\Gamma^{\left(pc\right)}|\left(\frac{1}{\epsilon} - \langle\chi^{\left(c\right)\left[1\right]}\rangle_{\Gamma^{\left(pc\right)}} + \langle\chi^{\left(p\right)\left[2\right]}\rangle_{\Gamma^{\left(pc\right)}}\right),}
\neweq{}{R_2^{\left(p\right)} = \frac{\phi^{\left(p\right)}}{\phi^{\left(c\right)}}R_1^{\left(p\right)},}
\neweq{}{R_3^{\left(p\right)} = \phi^{\left(p\right)^2}\left[\frac{\mathcal{Q}|\Gamma^{\left(pw\right)}|}{|\mathcal{B}^{\left(p\right)}|\epsilon} + \frac{\text{Bi}^{\left(p\right)}}{|\mathcal{B}^{\left(p\right)}|} |\Gamma^{\left(pc\right)}| \langle \chi^{\left(p\right)\left[1\right]} \rangle_{\Gamma^{\left(pc\right)}}\right],}
\neweq{}{\R_4^{\left(p\right)} = \phi^{\left(p\right)}k^{\left(p\right)}\langle\nabla_{\xiv}\chi^{\left(p\right)\left[1\right]}\rangle_{Y}.}
\neweq{}{\U^{\left(c\right)} = \varrho \varsigma \left[\phi^{\left(c\right)}\frac{\text{Bi}^{\left(c\right)}}{|\mathcal{B}^{\left(c\right)}|} |\Gamma^{\left(pc\right)}| \langle \chiv^{\left(c\right)\left[2\right]} \rangle_{\Gamma^{\left(pc\right)}} + k^{\left(c\right)}\langle\nabla_{\xiv}\chi^{\left(c\right)\left[1\right]}\rangle_{Y}\right],}
\neweq{}{\V^{\left(c\right)} = \frac{\phi^{\left(c\right)}}{\phi^{\left(p\right)}}\varrho \varsigma \left[\phi^{\left(c\right)}\frac{\text{Bi}^{\left(c\right)}}{|\mathcal{B}^{\left(c\right)}|} |\Gamma^{\left(pc\right)}| \langle \chiv^{\left(p\right)\left[3\right]} \rangle_{\Gamma^{\left(pc\right)}} + k^{\left(c\right)}\langle\nabla_{\xiv}\chi^{\left(c\right)\left[1\right]}\rangle_{Y}\right],}
\neweq{}{\K^{\left(c\right)} = \varrho \varsigma k^{\left(c\right)}\left[\phi^{\left(c\right)}\I + \langle\nabla_{\xiv}\chiv^{\left(c\right)\left[2\right]}\rangle_{Y}\right],}
\neweq{}{R_1^{\left(c\right)} = \frac{\phi^{\left(c\right)}}{\phi^{\left(p\right)}} R_2^{\left(c\right)},}
\neweq{}{R_2^{\left(c\right)} = \phi^{\left(c\right)}\frac{\left(\text{Bi}^{\left(c\right)} \varrho \varsigma\right)}{|\mathcal{B}^{\left(c\right)}|} |\Gamma^{\left(pc\right)}| \left(\frac{1}{\epsilon} - \langle \chi^{\left(c\right)\left[1\right]} \rangle_{\Gamma^{\left(pc\right)}} + \langle \chi^{\left(p\right)\left[2\right]} \rangle_{\Gamma^{\left(pc\right)}}\right),}
\neweq{}{R_3^{\left(c\right)} = \phi^{\left(c\right)^2}\frac{\left(\text{Bi}^{\left(c\right)} \varrho  \varsigma\right)}{|\mathcal{B}^{\left(c\right)}|} |\Gamma^{\left(pc\right)}| \langle \chi^{\left(p\right)\left[1\right]} \rangle_{\Gamma^{\left(pc\right)}},}
\neweq{}{R_4^{\left(c\right)} = \phi^{\left(c\right)^2} \varrho \mathcal{R},}
\end{subequations}
{In equations \eqref{upscaled}~--~\eqref{effectivecoefficients}, $\langle T^{\left(p\right)} \rangle_{Y}\left(t, \mathbf{x}\right)$ and $\langle T^{\left(c\right)} \rangle_{Y}\left(t, \mathbf{x}\right)$ are the averaged packing and battery cell temperatures, $\overline{\Pi}(\langle T^{\left(c\right)} \rangle_{Y}, \mathbf{x})$ is the homogenized power flux source term and ${\Pi}_{\text{base}}(\hat{\mathbf{x}})$ is the dimensionless \textit{base} power flux values, respectively.}  $\U^{\left(i\right)}$ are effective velocities, $\V_{i}^{\left(\bm{\cdot}\right)}$ are effective parameters corresponding to the emergent terms, $\K^{\left(i\right)}$ are effective thermal conductivities, $R_{\bm{\cdot}}^{\left(i\right)}$ and $\R_{5}^{\left(p\right)}$ are effective reaction rates, and ${\chi}_0^{p}$, ${\chi}_1^{p}$, $\boldsymbol{\chi}_2^{p}$ and ${\chi}_0^{c}$ are closure problems that must be solved in the unit cell with periodic boundary conditions. Since these closure problems are only solved once, the increase in computation cost is negligible. A detailed formulation of closure problems can be found in~\ref{subsection:Appendix_E_Closure_Problems} or in Pietrzyk et al.~\cite{Pietrzyk2023-ou}.

\subsection{Hybrid governing equations}
\label{sec:coupling-conditions}
\begin{figure}
\centerline{
 {\includegraphics[width=.7\linewidth]{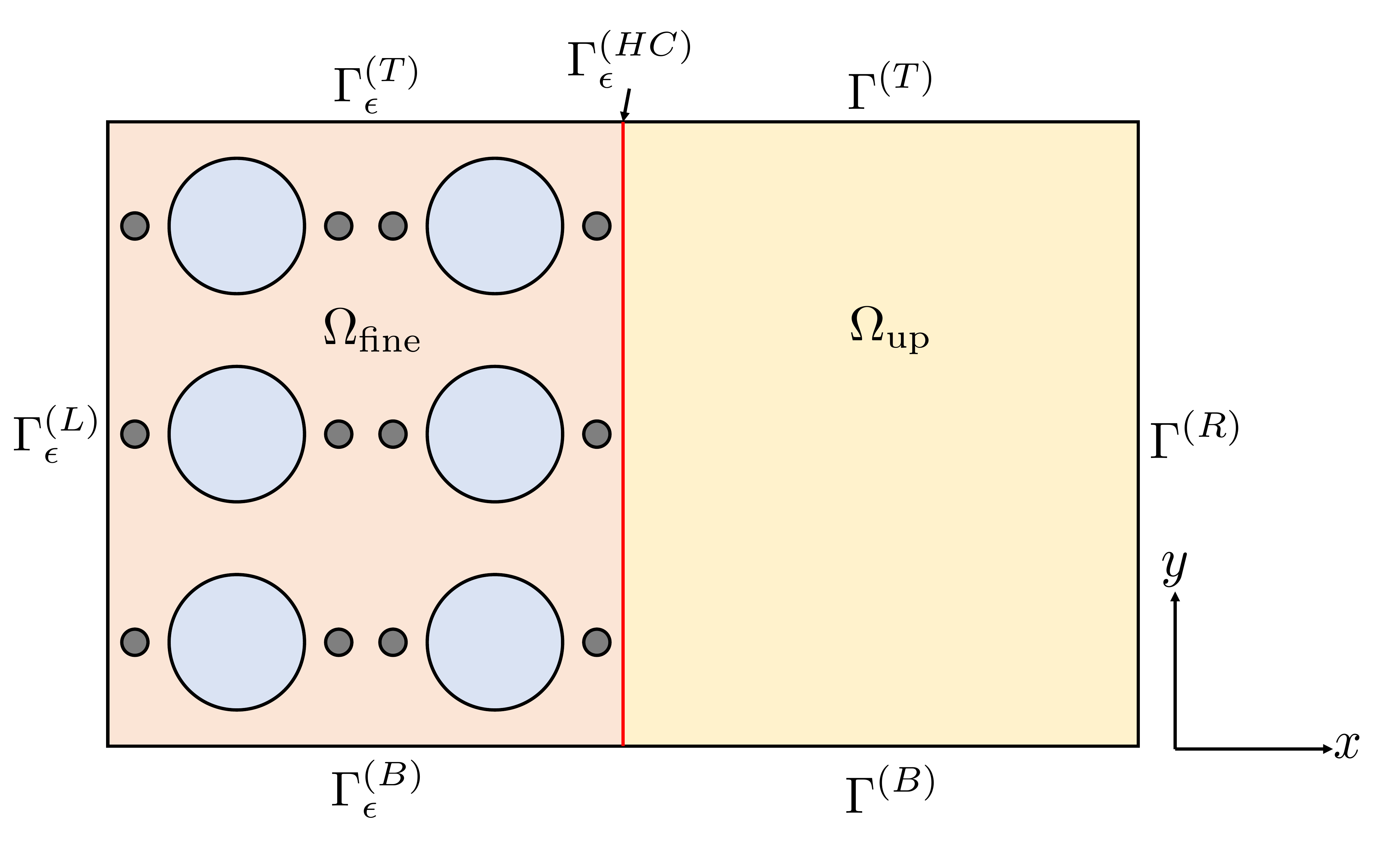}}}
\caption{{Schematic diagram of a hypothetical hybrid domain with fine-scale and upscaled subdomains. The sizes of the fine-scale and upscaled subdomains are arbitrary.}}
\label{fig:het-domain-eg}
\end{figure}

Assuming the existence of breakdown regions ${\Omega}_{\text{break}}$ in a dimensionless battery pack ${\Omega}_{\epsilon}$, we consider a two-dimensional heterogeneous battery pack ${\Omega}_{\epsilon,\text{het}} \in \mathbb{R}^2$ where the homogeneous assumption in the upscaled governing equations~\eqref{eq:upscaled_goven_eqs} is violated in the breakdown regions. Instead of solving the entire ${\Omega}_{\epsilon,\text{het}}$ with fine-scale equations~\eqref{eq:pore_goven_eqs} which is computationally expensive, a hybrid approach that combines the advantages of fine-scale and upscaled approaches is preferred. In hybrid simulations, the fine-scale subdomain {${\Omega}_{\text{fine}}$} is defined as ${\Omega}_{\text{break}} \cap {\Omega}_{\epsilon,\text{het}}$ while the upscaled subdomain ${\Omega}_{\text{up}}$ is defined as ${\Omega}_{\text{break}} \backslash {\Omega}_{\epsilon,\text{het}}$. There,  {${\Omega}_{\text{fine}}$} and ${\Omega}_{\text{up}}$ are solved by fine-scale and upscaled equations, respectively. These two subdomains are coupled by a coupling boundary ${\Gamma}^{\left(HC\right)}_{\epsilon} \subset \mathbb{R}^2$. In this study, we simplify the problem by assuming that ${\Gamma}^{\left(HC\right)}_{\epsilon}$ is defined in $\mathcal{B}_\epsilon^{\left(p\right)}$ and only dependent on $x$ instead of both $x$ and $y$. Therefore, ${\Omega}_{\text{up}}$ is bounded by packing edges ${\Gamma}^{\left(R\right)}$,  ${\Gamma}^{\left(T\right)}$ and ${\Gamma}^{\left(B\right)}$ and the coupling boundary ${\Gamma}^{\left(HC\right)}_{\epsilon}$, while the fine-scale subdomain {${\Omega}_{\text{fine}}$} is bounded by packing edges ${\Gamma}_\epsilon^{\left(L\right)}$,  ${\Gamma}_\epsilon^{\left(T\right)}$ and ${\Gamma}_\epsilon^{\left(B\right)}$ and the coupling boundary ${\Gamma}^{\left(HC\right)}_{\epsilon}$ (Figure~\ref{fig:het-domain-eg}). 

% Let's focus on a  region across the coupling boundary with $\mathbf{x}_{HC} \in \Gamma_\epsilon^{\left(HC\right)}$ (Figure~\ref{fig:hybrid_illu_1}). {We now define two arbitrary volumes, which live in} {${\Omega}_{\text{fine}}$} and ${\Omega}_{\text{up}}$, with centroids $\mathbf{x}^{+}$ and $\mathbf{x}^{-}$. {As $\mathbf{x}^{+}$ and $\mathbf{x}^{-}$ approach $\mathbf{x}_{HC}$ from the right and the left, respectively, =======

% move toward and eventually coincide with $\mathbf{x}_{HC} \in \Gamma_\epsilon^{\left(HC\right)}$ from their respective subdomains where} $\mathbf{x}_{HC}$ is the centroid of the coupling volume $Y_{HC}$, heat and flux conservation will be guaranteed (Figure~\ref{fig:hybrid_illu_2}) as 

\begin{figure*}
    \centering
    \begin{subfigure}[ht]{0.7\textwidth}
    \caption{}
    \centerline{
     {\includegraphics[width=\textwidth]{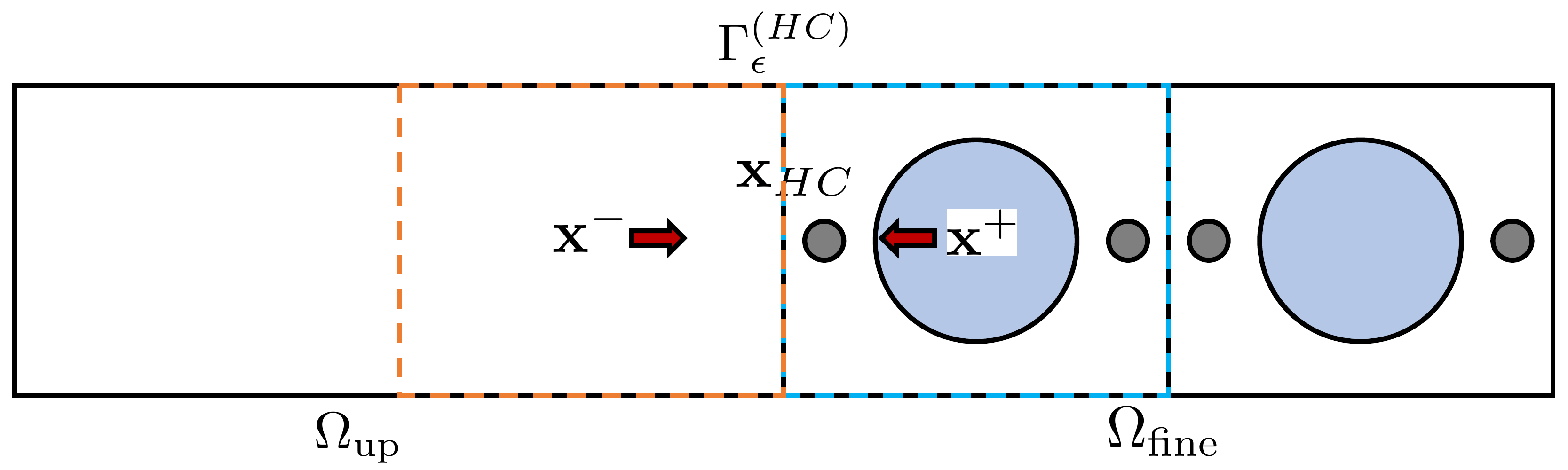}}}
    \label{fig:hybrid_illu_1}
    \end{subfigure}
    
    \begin{subfigure}[ht]{0.7\textwidth}
    \caption{}
    \centerline{
     {\includegraphics[width=\textwidth]{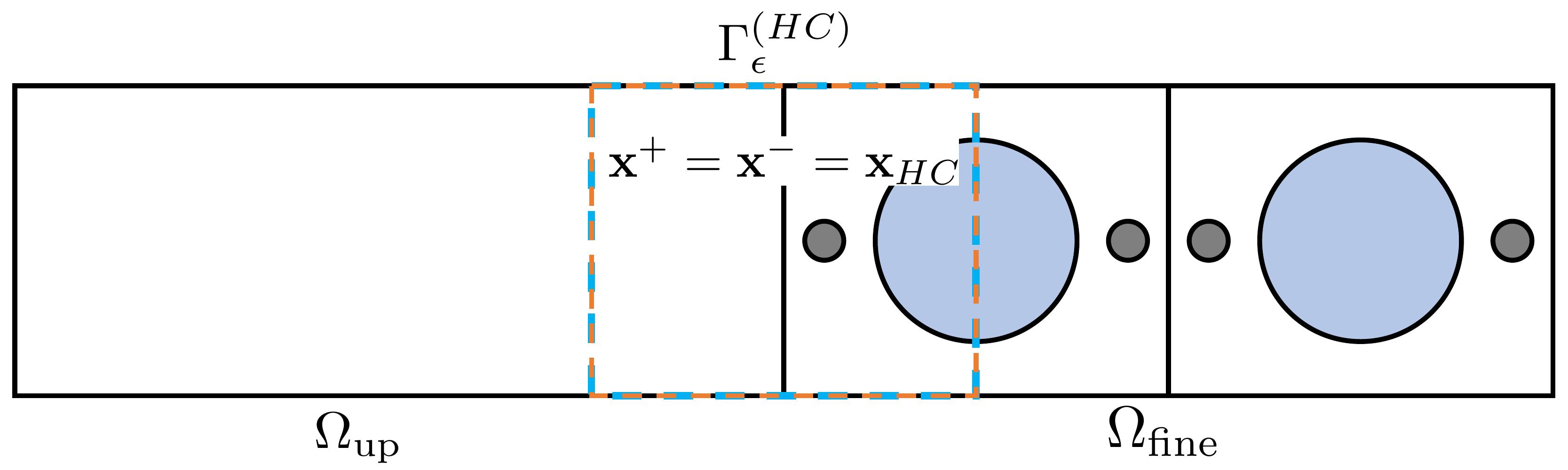}}}
    \label{fig:hybrid_illu_2}
    \end{subfigure}
    
    \begin{subfigure}[ht]{0.7\textwidth}
    \caption{}
    \centerline{
     {\includegraphics[width=\textwidth]{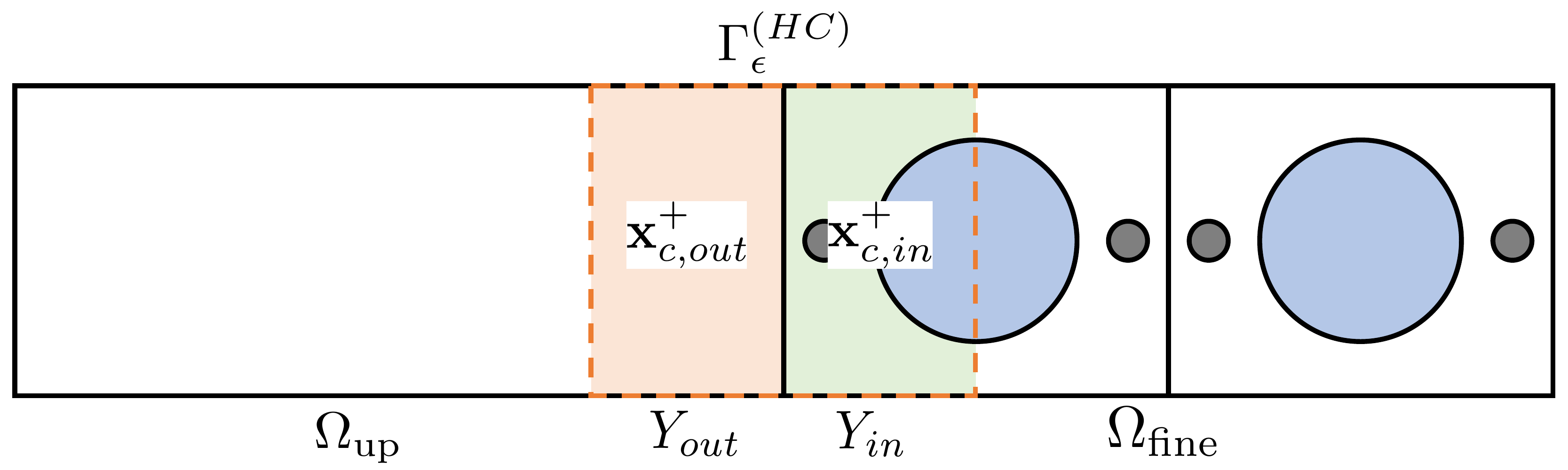}}}
    \label{fig:hybrid_illu_3}
    \end{subfigure}
     \caption{{Schematic illustration of the formulation of the hybrid approach. (a) Two volumes with centroids $\mathbf{x}^{+}$ in the fine-scale domain (blue) and $\mathbf{x}^{-}$ in the upscaled domain (orange) move toward each other, (b) two volumes coincide where $\mathbf{x}^{+} = \mathbf{x}^{-} = \mathbf{x}_{HC}$, and (c) the coupling volume from the fine-scale domain is divided into two regions $Y_{in}$ and $Y_{out}$ with centroids $\mathbf{x}_{c,in}^{+}$ and $\mathbf{x}_{c,out}^{+}$, respectively, where $Y_{in}$ intersects with the fine-scale domain and $Y_{out}$ intersects with the upscaled domain.}}
     \label{fig:hybrid_illu_all}
\end{figure*}

{

Let's focus on a region across the coupling boundary with $\mathbf{x}_{HC} \in \Gamma_\epsilon^{\left(HC\right)}$ (Figure~\ref{fig:hybrid_illu_1}). %\textcolor{red}{[Here, it looks like you use $\langle \cdot \rangle_{Y}$ on Pore-scale values. I know what you are trying to say, but it might be worth defining a separate average for this, since you technically can only use that average to denote upscaled quantities. If you look at the other battery paper, we define an average over moving domain $\mathcal{W}\left(\mathbf{x}\right)$ (which is rigorously defined) for what you are trying to do here. It looks like the reviewer has some confusion on similar points as well (i.e., doesn't understand the double integrals in the appendix, etc.]}
We now define two average temperatures $\ensmean{T_\epsilon^{\left(p\right)}}_{Y}(\mathbf{x}^{+})$ and $\ensmean{T^{\left(p\right)}}_{Y}(\mathbf{x}^{-})$, with centroids $\mathbf{x}^{+} \in \Omega_{\text{fine}}$ and $\mathbf{x}^{-} \in \Omega_{\text{up}}$ , respectively. As $\mathbf{x}^{+}$ and $\mathbf{x}^{-}$ approach $\mathbf{x}_{HC}$ from the right and the left, %\textcolor{red}{[should be left and right? not right and left]}, 
respectively, the limiting average temperatures, $\ensmean{T_\epsilon^{\left(p\right)}}_{Y}(\mathbf{x}_{HC}) :=\lim_{\mathbf{x}^{+} \rightarrow \mathbf{x}_{HC}} \ensmean{T_\epsilon^{\left(p\right)}}_{Y}(\mathbf{x}^{+})$ and  $\ensmean{T^{\left(p\right)}}_{Y}(\mathbf{x}_{HC}) :=\lim_{\mathbf{x}^{-} \rightarrow \mathbf{x}_{HC}} \ensmean{T^{\left(p\right)}}_{Y}(\mathbf{x}^{-})$, will be equivalent, satisfying the continuity of heat. By applying a similar concept to the heat flux across the coupling boundary, the continuities of heat and flux are defined as (Figure~\ref{fig:hybrid_illu_2})} %\textcolor{red}{[I think there might be implications to the $\mathbf{n}$ being outside the integral (i.e., boundary HC must be straight). I don't know if you want to address this or not]}

% volumes, which live in{${\Omega}_{\text{fine}}$} and ${\Omega}_{\text{up}}$, with centroids $\mathbf{x}^{+}$ and $\mathbf{x}^{-}$. {As $\mathbf{x}^{+}$ and $\mathbf{x}^{-}$ approach $\mathbf{x}_{HC}$ from the right and the left, respectively, =======

% move toward and eventually coincide with $\mathbf{x}_{HC} \in \Gamma_\epsilon^{\left(HC\right)}$ from their respective subdomains where} $\mathbf{x}_{HC}$ is the centroid of the coupling volume $Y_{HC}$, heat and flux conservation will be guaranteed (Figure~\ref{fig:hybrid_illu_2}) as 

\begin{subequations}
\label{eq:coupling-conds}
\begin{align}
&{\ensmean{T^{\left(p\right)}}_{Y}(\mathbf{x}^{-}) = \ensmean{T_\epsilon^{\left(p\right)}}_{Y}(\mathbf{x}^{+}),  \quad \text{for} \quad \abs{\mathbf{x}^{+} - \mathbf{x}^{-}} \rightarrow 0, }\\
& {\ensmean{\mathbf{J}^{\left(p\right)}}_{Y}(\mathbf{x}^{-})  \bm{\cdot} \mathbf{n}^{\left(p\right)}_\epsilon = \phi^{\left(p\right)}\ensmean{\mathbf{J}_\epsilon^{\left(p\right)}}_{Y}(\mathbf{x}^{+})  \bm{\cdot} \mathbf{n}^{\left(p\right)}_\epsilon, \quad \text{for} \quad \abs{\mathbf{x}^{+} - \mathbf{x}^{-}} \rightarrow 0,}\label{eq:packing_flux_hc}
\end{align}
\end{subequations}
where  ${\mathbf{J}_\epsilon^{\left(p\right)}}(\mathbf{x}^{+})$ and $\ensmean{\mathbf{J}^{\left(p\right)}}_{Y}(\mathbf{x}^{-})$ are the packing material fluxes in the fine-scale and upscaled equations, respectively. Although there are packing and cell temperatures, heat and flux continuities are only required for the packing temperature. The derivations of the flux conditions in equation~\eqref{eq:packing_flux_hc} can be found in~\ref{app:packing-flux-proof} and~\ref{app:cell-flux-proof}. One challenge associated with equation~\eqref{eq:coupling-conds} is that the average fine-scale quantities {(RHS of equation~\eqref{eq:coupling-conds})} are unknown. Therefore, appropriate approximation methods such as the Taylor (Section~\ref{sec:taylor}) or Series (Section~\ref{sec:series}) expansion approach are needed to compute the average fine-scale quantities.

% In addition, as demonstrated in Section X, since equation~\eqref{eq:cell_flux_hc} is valid, no boundary conditions on the cell temperatures in fine-scale and upscaled simulations have been imposed, resulting in a jump or discontinuity at the coupling boundary due to errors between fine-scale and upscaled simulations. To eliminate the discontinuity and enforce equations~\ref{eq:coupling-conds} simultaneously, we added a flux to the battery cell nearest to the coupling boundary (Section~\ref{sec:taylor} and~\ref{sec:series}).  
\subsubsection{Taylor expansion approach}
\label{sec:taylor}
At the coupling boundary, we separate the coupling volume $Y(\mathbf{x}_{HC})$ into $Y_{in}(\mathbf{x}^{+}_{c,in}) = Y \cap {\Omega}_{\text{fine}}$ and $Y_{out}(\mathbf{x}^{+}_{c,out}) = Y \cap {\Omega}_{\text{up}}$ where $\mathbf{x}^{+}_{c,in}$ and $\mathbf{x}^{+}_{c,out}$ are the centroids of $Y_{in}$ and $Y_{out}$, respectively (Figure~\ref{fig:hybrid_illu_3}). Using this definition, $\ensmean{T_\epsilon^{\left(p\right)}}_{Y}(\mathbf{x}^{+})$ can be expressed as %\textcolor{red}{[can maybe get rid of the $"\in Y_{in}"$ and $"\in Y_{out}"$]}
\begin{align}
\label{eq:Tp}
&\ensmean{T_\epsilon^{\left(p\right)}}_{Y}(\mathbf{x}^{+}) = \ddfrac{1}{\abs{Y}} \int_{\mathcal{B}^{\left(p\right)}_{in}} T_\epsilon^{\left(p\right)}(\mathbf{y}) \ \dd\mathbf{y} + \ddfrac{1}{\abs{Y}} \int_{\mathcal{B}^{\left(p\right)}_{out} } T_\epsilon^{\left(p\right)}(\mathbf{y}) \ \dd\mathbf{y},
\end{align}
where $\mathcal{B}^{\left(p\right)}_{in}$ and $\mathcal{B}^{\left(p\right)}_{out}$ are the packing material domains in $Y_{in}$ and $Y_{out}$, respectively.

Since the quantities in $Y_{out}$ are solved with upscaled equations, the fine-scale packing temperature $T_\epsilon^{\left(p\right)}$ is unknown. To approximate the second integral in equation~\eqref{eq:Tp}, one approach is use a  Taylor series at the nearest defined location (the coupling boundary $\mathbf{x}_{HC}$). 
% However, $T_\epsilon^{\left(p\right)}(\mathbf{x}_{HC})$ is not defined unless $\mathbf{x}_{HC} \in \mathcal{B}^p$ is satisfied. Therefore, instead of expanding at $\mathbf{x}_{HC}$, the integral is expanded at $\mathbf{x}^{*}_{HC}$ and retains the first two terms where $\mathbf{x}^{*}_{HC}$ is defined as
% \begin{align}
%     &\mathbf{x}^{*}_{HC} = \argmin_{y} \left(y-\mathbf{x}_{HC} \right),
% \end{align}
% where $y \in \mathcal{B}^p$ and $\abs{\mathbf{x}^{*}_{HC} - \mathbf{x}_{HC}} < \epsilon$ are satisfied to ensure that $\mathbf{x}^{*}_{HC}$ is in the same unit cell as $\mathbf{x}_{HC}$. 
As such, the integral can be expressed as
\begin{align}
&\int_{\mathcal{B}^{\left(p\right)}_{out}} T_\epsilon^{\left(p\right)}(\mathbf{y}) \ \dd\mathbf{y} = \int_{\mathcal{B}^{\left(p\right)}_{out}} \left[ T_\epsilon^{\left(p\right)}(\mathbf{x}_{HC}) + \pdv{T_\epsilon^{\left(p\right)}(\mathbf{x}_{HC})}{\mathbf{x}} \left(\mathbf{y} - \mathbf{x}_{HC}\right) + \mathcal{O}\left((\mathbf{y}-\mathbf{x}_{HC})^2\right) \right] \ \dd\mathbf{y}.
\end{align}
By defining the centroid of $Y_{out}$ as
\begin{align}
&\mathbf{x}^{+}_{c,out} = \ddfrac{1}{\abs{\mathcal{B}^{\left(p\right)}_{out}} } \int_{\mathcal{B}^{\left(p\right)}_{out}} \mathbf{y} \ \dd\mathbf{y},
\end{align}
and $\abs{\mathcal{B}^{\left(p\right)}_{out}} = \abs{Y_{out}}\phi^{\left(p\right)}_{out}$, the approximated integral is then expressed as
\begin{align}
\int_{\mathcal{B}^{\left(p\right)}_{out}} T_\epsilon^{\left(p\right)}(\mathbf{y}) \ \dd\mathbf{y} &= \abs{Y_{out}}\phi^{\left(p\right)}_{out} \left[ T_\epsilon^{\left(p\right)}(\mathbf{x}_{HC}) + \pdv{T_\epsilon^{\left(p\right)}(\mathbf{x}_{HC})}{\mathbf{x}} \left(\mathbf{x}^+_{c,out} - \mathbf{x}_{HC}\right) \right] \\
&+  \int_{\mathcal{B}^{\left(p\right)}_{out}} \mathcal{O}\left((\mathbf{y}-\mathbf{x}_{HC})^2\right) \ \dd\mathbf{y}
\end{align}
where $\phi^{\left(p\right)}_{out}$ is the volume fraction of the packing materials in $Y_{out}$ such that
\begin{align}
    \phi^{\left(p\right)}_{out} = \frac{\abs{Y}\phi^{(p)} -  \int_{\mathcal{B}^{\left(p\right)}_{in}} \mathbf{1} \ \dd\mathbf{y}}{\abs{Y_{out}}}.
\end{align}
Finally, by truncating at the second-order term, the integral is approximated as 
% However, $\phi^{\left(p\right)}_{out}$ cannot be computed because $Y_{out}$ is only defined in the upscaled domain, which does not contain information on battery cells and cooling water pipes. By ensuring $\abs{Y_{out}} \ge \abs{Y}$, $\phi^{\left(p\right)}_{out} = \phi^{\left(p\right)}$ is guaranteed where $\phi^{\left(p\right)}$ is the known volume fraction of packing material in a unit cell that can be computed with the known unit cell geometry. Finally, the second-order accurate approximate integral is simplified to 

\begin{align}
\label{eq:taylor-tp-approx}
&\int_{\mathcal{B}^{\left(p\right)}_{out}} T_\epsilon^{\left(p\right)}(\mathbf{y}) \ \dd\mathbf{y} \approx \abs{Y_{out}}\phi^{\left(p\right)}_{out} \left[ T_\epsilon^{\left(p\right)}(\mathbf{x}_{HC}) + \pdv{T_\epsilon^{\left(p\right)}(\mathbf{x}_{HC})}{\mathbf{x}} \left(\mathbf{x}^+_{c,out} - \mathbf{x}_{HC}\right) \right].
\end{align}
By substituting equation~\eqref{eq:taylor-tp-approx} into equation~\eqref{eq:Tp}, the average packing temperature at the coupling boundary can be approximated as 
\begin{align}
\label{eq:approx-tp}
&\ensmean{T_\epsilon^{\left(p\right)}}_{Y}(\mathbf{x}^{+}) \approx \ensmean{T_\epsilon^{\left(p\right)}}_{Y_{in}}(\mathbf{x}^{+}) + \alpha^{-1} \phi^{\left(p\right)}_{out} \left[ T_\epsilon^{\left(p\right)}(\mathbf{x}_{HC}) + \pdv{T_\epsilon^{\left(p\right)}(\mathbf{x}_{HC})}{\mathbf{x}} \left(\mathbf{x}^+_{c,out} - \mathbf{x}_{HC}\right) \right],
\end{align}
where $\ensmean{T_\epsilon^{\left(p\right)}}_{Y_{in}}(\mathbf{x}^{+}) = \ddfrac{1}{\abs{Y}} \int_{\mathcal{B}^{\left(p\right)}_{in}} T_\epsilon^{\left(p\right)}(\mathbf{y}) \ \dd\mathbf{y}$ and $\alpha = \abs{Y}/\abs{Y_{out}}$.

% By applying a similar approach to approximate the average cell temperature $\ensmean{T_\epsilon^{c}}_Y_{Y}(\mathbf{x}^{+})$, we derive the second-order accurate approximation as 
% \begin{align}
% \label{eq:approx-tc}
% &\ensmean{T_\epsilon^{c}}_Y_{Y}(\mathbf{x}^{+}) \approx \ensmean{T_\epsilon^{c}}_Y_{Y_{in}}(\mathbf{x}^{+}) + \alpha^{-1} \phi^{c} \left[ T_\epsilon^{c}(\mathbf{x}^{*}_{HC}) + \pdv{T_\epsilon^{c}(\mathbf{x}^{*}_{HC})}{\mathbf{x}} \left(\mathbf{x}^+_{c,out} - \mathbf{x}^{*}_{HC}\right) \right],
% \end{align}
% where $\phi^{c}$ is the known volume fraction of battery cell in a unit cell. A point to note is that $\mathbf{x}^{*}_{HC}$ in equations~\ref{eq:approx-tp} and~\ref{eq:approx-tc} might not be equivalent since the approximate centroids are defined on the respective domain.

To derive the average flux of packing material $\ensmean{\mathbf{J}_\epsilon^{(p)}}_{Y}(\mathbf{x}^{+})$, a similar approach as that used in equation~\eqref{eq:Tp} is employed, such that 
\begin{align}
\phi^{\left(p\right)}\ensmean{\mathbf{J}_\epsilon^{\left(p\right)}}_{Y}(\mathbf{x}^{+})  \bm{\cdot} \mathbf{n}^{\left(p\right)}_\epsilon &= \ddfrac{\phi^{\left(p\right)}}{\abs{Y}} \int_{\mathcal{B}^{\left(p\right)}_{in}} \mathbf{J}_\epsilon^{\left(p\right)}(\mathbf{y}) \ \dd\mathbf{y}  \bm{\cdot} \mathbf{n}_{\epsilon}^{\left(p\right)} + \ddfrac{\phi^{\left(p\right)}}{\abs{Y}} \int_{\mathcal{B}^{\left(p\right)}_{out}} \mathbf{J}_\epsilon^{\left(p\right)}(\mathbf{y}) \ \dd\mathbf{y} \bm{\cdot} \mathbf{n}_{\epsilon}^{\left(p\right)} \nonumber \\
&\approx \phi^{\left(p\right)}\ensmean{\mathbf{J}_\epsilon^{\left(p\right)}}_{Y_{in}}(\mathbf{x}^{+}) \bm{\cdot} \mathbf{n}_{\epsilon}^{\left(p\right)} + q^{\left(p,n\right)},
\end{align}
where {$\ensmean{\mathbf{J}_\epsilon^{\left(p\right)}}_{Y_{in}}(\mathbf{x}^{+}) =  \ddfrac{1}{\abs{Y}} \int_{\mathcal{B}^{\left(p\right)}_{in}} \mathbf{J}_\epsilon^{\left(p\right)}(\mathbf{y}) \ \dd\mathbf{y}$} and $q^{\left(p,n\right)} = {\ddfrac{\phi^{\left(p\right)}}{\abs{Y}} \int_{\mathcal{B}^{\left(p\right)}_{out}} \mathbf{J}_\epsilon^{\left(p\right)}(\mathbf{y}) \ \dd\mathbf{y} \bm{\cdot} \mathbf{n}_{\epsilon}^{\left(p\right)} }$ is the {unresolved} %\text{color}{is it unknown or unresolved? I am not sure we talked about guessing yet. Guessing is how we handle this unknown.} 
flux of the packing material at iteration $n$. {To derive the flux boundary conditions of the fine-scale equation, we expand the unresolved flux at $\mathbf{x}_{HC}$ and the integral can be approximated as }% \textcolor{red}{[It might just be me, but I got a little lost here; I see what's happening after rereading it a few times, but having a formula for a ``guess flux" was a little weird at first.]}
\begin{align}
\label{eq:approx-flux-long}
   q^{\left(p,n\right)} &= \alpha^{-1}\phi^{\left(p\right)} \phi^{\left(p\right)}_{out}\left[ \mathbf{J}_\epsilon^{\left(p\right)}(\mathbf{x}_{HC}) + \pdv{\mathbf{J}_\epsilon^{\left(p\right)}(\mathbf{x}_{HC})}{\mathbf{x}} \left(\mathbf{x}^+_{c,out} - \mathbf{x}_{HC}\right) \right] \bm{\cdot} \mathbf{n}_{\epsilon}^{\left(p\right)} \nonumber \\
  &+  \int_{\mathcal{B}^{\left(p\right)}_{out}} \mathcal{O}\left((\mathbf{y}-\mathbf{x}_{HC})^2\right) \ \dd\mathbf{y}.
\end{align}
To relate $q^{\left(p,n\right)}$ to the flux boundary conditions of the fine-scale equations $\mathbf{J}_\epsilon^{(p)}(\mathbf{x}_{HC})$, equation~\eqref{eq:approx-flux-long} is truncated to the first-order  and approximated as 
\begin{align}
\label{eq:approx-flux}
  &  q^{\left(p,n\right)} \approx \alpha^{-1}\phi^{\left(p\right)} \phi^{\left(p\right)}_{out} \mathbf{J}_\epsilon^{\left(p\right)}(\mathbf{x}_{HC}) \bm{\cdot} \mathbf{n}_{\epsilon}^{\left(p\right)} .
\end{align}
By rearranging the equation~\eqref{eq:approx-flux}, the flux boundary condition for the fine-scale equations is defined as 
\begin{align}
\mathbf{J}_\epsilon^{\left(p\right)}(\mathbf{x}_{HC}) \bm{\cdot} \mathbf{n}_{\epsilon}^{\left(p\right)}  \approx {\alpha}{\left(\phi^{\left(p\right)} \phi^{\left(p\right)}_{out}\right)}^{-1} q^{\left(p,n\right)}. 
\end{align}
In summary, by using Taylor series expansion, we derived the first-order accurate coupling boundary conditions for both fine-scale and upscaled equations as 
\begin{subequations}
\label{eq:taylor-hc}
\begin{align}
&\ensmean{T_\epsilon^{\left(p\right)}}_{Y}(\mathbf{x}^{+}) \approx \ensmean{T_\epsilon^{\left(p\right)}}_{Y_{in}}(\mathbf{x}^{+}) + \alpha^{-1} \phi^{\left(p\right)}_{out} \left[ T_\epsilon^{\left(p\right)}(\mathbf{x}_{HC}) + \pdv{T_\epsilon^{\left(p\right)}(\mathbf{x}_{HC})}{\mathbf{x}} \left(\mathbf{x}^+_{c,out} - \mathbf{x}_{HC}\right) \right], \label{eq:taylor-hc-temp}\\
&\phi^{\left(p\right)}\ensmean{\mathbf{J}_\epsilon^{\left(p\right)}}_{Y}(\mathbf{x}^{+})  \bm{\cdot} \mathbf{n}^{\left(p\right)}_\epsilon \approx \phi^{\left(p\right)}\ensmean{\mathbf{J}_\epsilon^{\left(p\right)}}_{Y_{in}}(\mathbf{x}^{+}) \bm{\cdot} \mathbf{n}_{\epsilon}^{\left(p\right)} + q^{\left(p,n\right)}, \label{eq:taylor-hc-flux}\\
& \mathbf{J}_\epsilon^{\left(p\right)}(\mathbf{x}_{HC}) \bm{\cdot} \mathbf{n}_{\epsilon}^{\left(p\right)} \approx {\alpha}{\left(\phi^{\left(p\right)} \phi^{\left(p\right)}_{out}\right)}^{-1} q^{\left(p,n\right)}.
\end{align}
\end{subequations}
{For equation~\eqref{eq:taylor-hc} to be valid, two conditions on $Y$ must be satisfied: (1) $\Gamma^{\left(HC\right)}$ only intersects with the packing materials, and (2) $\mathbf{x}_{HC} \in \Gamma^{\left(HC\right)}$ are the centroids of the coupling volume $Y$.}

\subsubsection{Series expansion approach}
\label{sec:series}
% In Section~\ref{sec:taylor}, we show that the Taylor series approach is first-order accurate due to the relationship between the guess flux and the flux boundary conditions of the fine-scale equations. 
{In Section~\ref{sec:taylor}, we show that the Taylor series approach is first-order accurate for the approximation of coupling fluxes (Equation~\eqref{eq:taylor-hc-flux}) because the Taylor expansion for flux is truncated at the first-order.} According to Pietrzyk et al.~\cite{Pietrzyk2021-lu}, \symbolica is capable of achieving higher-order accurate upscaled equations; therefore, limiting the hybrid coupling method to be first-order accurate narrows its applicability. Therefore, we propose a one-sided scheme that uses Series expansion to derive high-order accurate coupling conditions. The main advantage of this approach is that the information in the upscaled subdomain will not be {necessary}. For any quantity of interest $\langle \bm{\cdot} \rangle_{Y} (\mathbf{x}^+)$, in a homogeneous domain, we apply a second-order accurate approximation such that %\textcolor{red}{[It took me awhile to figure out equation (28) came from (i.e., Taylor expansions). As such, I think this paragraph might lean-in too hard on ``This is not Taylor expansion; it is something else". Maybe just call it "An approach for second order accuracy"?]}
\begin{align}
\label{eq:general_series}
\langle \bm{\cdot} \rangle_Y \left(\mathbf{x}^++\ddfrac{\epsilon}{2}\right) \approx \ddfrac{ \langle \bm{\cdot} \rangle_Y (\mathbf{x}^+) + \langle \bm{\cdot} \rangle_Y (\mathbf{x}^++\epsilon)}{2},
\end{align}
where %\textcolor{red}{[shouldn't it be $\mathcal{B}^{\left(i\right)}\left(\mathbf{x}\right)$? Otherwise I think $\langle \bm{\cdot} \rangle_{Y}\left(\mathbf{x}\right)$ should just be $\langle \bm{\cdot} \rangle_{Y}$. See previous comment on defining an average over $\mathcal{W}$.]}
\begin{align}
    &\langle \bm{\cdot} \rangle_Y (\mathbf{x})  = \ddfrac{1}{\abs{Y}}\int_{\mathcal{B}^{(i)}(\mathbf{x})} \left( \bm{\cdot} \right)^{(i)} \ \dd \mathbf{y}.
\end{align}
By rearranging the equation~\eqref{eq:general_series}, we can derive the second-order accurate approximation as 
\begin{align}
\label{eq:general_series_mod}
\langle \bm{\cdot} \rangle_Y (\mathbf{x}^+) = 2\langle \bm{\cdot} \rangle_Y \left(\mathbf{x}^+ + \ddfrac{\epsilon}{2}\right) - \langle \bm{\cdot} \rangle_Y (\mathbf{x}^++\epsilon).
\end{align}
To derive the flux boundary conditions for fine-scale equations, we define the {unresolved} flux $q^{\left(p,n\right)}$ of the packing material as 
\begin{align}
q^{\left(p,n\right)} = \phi^{\left(p\right)}\ensmean{\mathbf{J}_\epsilon^{\left(p\right)}}_{Y}(\mathbf{x}^{+})  \bm{\cdot} \mathbf{n}^{\left(p\right)}_\epsilon.
\end{align}
By expanding the integral at the coupling boundary $\mathbf{x}_{HC}$ with Taylor expansion, we approximate the integral as %\textcolor{red}{[I think it should be written $\mathcal{O}\left(\left(\mathbf{y} - \mathbf{x}_{HC}\right)^2\right)$]}
\begin{align}
\label{eq:series-flux-to-ps}
q^{\left(p,n\right)} &=  \ddfrac{\phi^{\left(p\right)}}{\abs{Y}}\int_{\mathcal{B}_\epsilon^{\left(p\right)}}  \mathbf{J}_\epsilon^{\left(p\right)}(\mathbf{y})  \ \dd \mathbf{y} \bm{\cdot} \mathbf{n}^{(p)}_{\epsilon}  \nonumber \\
&= \ddfrac{\phi^{\left(p\right)}}{\abs{Y}} \int_{\mathcal{B}_\epsilon^{\left(p\right)}}  \left[\mathbf{J}_\epsilon^{\left(p\right)}(\mathbf{x}_{HC}) + \pdv{\mathbf{J}_\epsilon^{\left(p\right)}(\mathbf{x}_{HC})}{\mathbf{x}} \left(\mathbf{y} - \mathbf{x}_{HC}\right) + \mathcal{O}\left((\mathbf{y}-\mathbf{x}_{HC})^2\right) \right]\ \dd \mathbf{y}.
\end{align}
By defining the centroid of the unit cell $Y$ as
\begin{align}
    &\mathbf{x}_{c} = \ddfrac{1}{\phi^{\left(p\right)}\abs{Y} } \int_{\mathcal{B}_\epsilon^{\left(p\right)}} \mathbf{y} \ \dd\mathbf{y},
\end{align}
and $\abs{\mathbf{x}_{c} - \mathbf{x}_{HC}} = 0$ due to the collocation of the coupling boundary {$\mathbf{x}_{HC}$} and centroid of the unit cell {$\mathbf{x}_{c}$}, equation~\eqref{eq:series-flux-to-ps} can be simplified to %\textcolor{red}{[I think it should be written $\mathcal{O}\left(\left(\mathbf{y} - \mathbf{x}_{HC}\right)^2\right)$]}
\begin{align}
q^{\left(p,n\right)} &= \left(\phi^{\left(p\right)}\right)^2 \mathbf{J}_\epsilon^{\left(p\right)}(\mathbf{x}_{HC})  + \int_{\mathcal{B}_\epsilon^{\left(p\right)}}   \mathcal{O}\left((\mathbf{y}-\mathbf{x}_{HC})^2\right) \ \dd \mathbf{y}.
\end{align}
By rearranging the equation, we can express the fine-scale packing material flux as 
\begin{align}
& \mathbf{J}_\epsilon^{\left(p\right)}(\mathbf{x}_{HC}) \approx  {\left(\phi^{\left(p\right)}\right)}^{-2} q^{\left(p,n\right)}.
\end{align}
% To ensure continuity in average cell temperature, we add a flux boundary condition for the battery cell temperature as 
% \begin{align}
% & \mathbf{J}_\epsilon^{c}(\mathbf{x}_{HC})= \ddfrac{1}{\phi^c} q^{c,n},
% \end{align}
% where $q^{c,n}$ is the guess flux for the battery cell temperature.

In summary, the second-order accurate coupling conditions with Series expansion approach can be expressed as
\begin{subequations}
\label{eq:series-hc}
\begin{align}
&\ensmean{T_\epsilon^{\left(p\right)}}_{Y_{HC}}(\mathbf{x}^{+}) \approx 2\ensmean{T_\epsilon^{\left(p\right)}}_{Y}(\mathbf{x}^{+}+\ddfrac{\epsilon}{2}) - \ensmean{T_\epsilon^{\left(p\right)}}_{Y}(\mathbf{x}^{+}+\epsilon), \label{eq:series-hc-temp}\\
&\phi^{\left(p\right)}\ensmean{\mathbf{J}_\epsilon^{\left(p\right)}}_{Y}(\mathbf{x}^{+})  \bm{\cdot} \mathbf{n}^{\left(p\right)}_\epsilon \approx \phi^{\left(p\right)}\left(2\ensmean{\mathbf{J}_\epsilon^{\left(p\right)}}_{Y}(\mathbf{x}^{+}+\ddfrac{\epsilon}{2}) - \ensmean{\mathbf{J}_\epsilon^{\left(p\right)}}_{Y}(\mathbf{x}^{+}+\epsilon)\right) \bm{\cdot} \mathbf{n}^{\left(p\right)}_\epsilon, \label{eq:series-hc-flux}\\
&\mathbf{J}_\epsilon^{\left(p\right)}(\mathbf{x}_{HC}) \bm{\cdot} \mathbf{n}_{\epsilon}^{\left(p\right)}  \approx {\left(\phi^{\left(p\right)}\right)}^{-2} q^{\left(p,n\right)},
\end{align}
\end{subequations}

For numerical implementation, we utilized open-source software \fenics ~\cite{Alnaes2009-ty,Alnaes2014-kl,Alnaes2015-hx,Kirby2004-xd,Kirby2006-el,Olgaard2010-pd,Logg2010-mw,Logg2012-ul}. To avoid the time-step size constraint, we discretized the equations with the first-order implicit Euler method. Details on the numerical implementation of fine-scale and upscaled equations in \fenics can be found in Pietrzyk et al.~\cite{Pietrzyk2023-ou}.
\subsection{Summary of the hybrid simulation algorithm}

\begin{figure}
\centerline{
 {\includegraphics[width=0.8\textwidth]{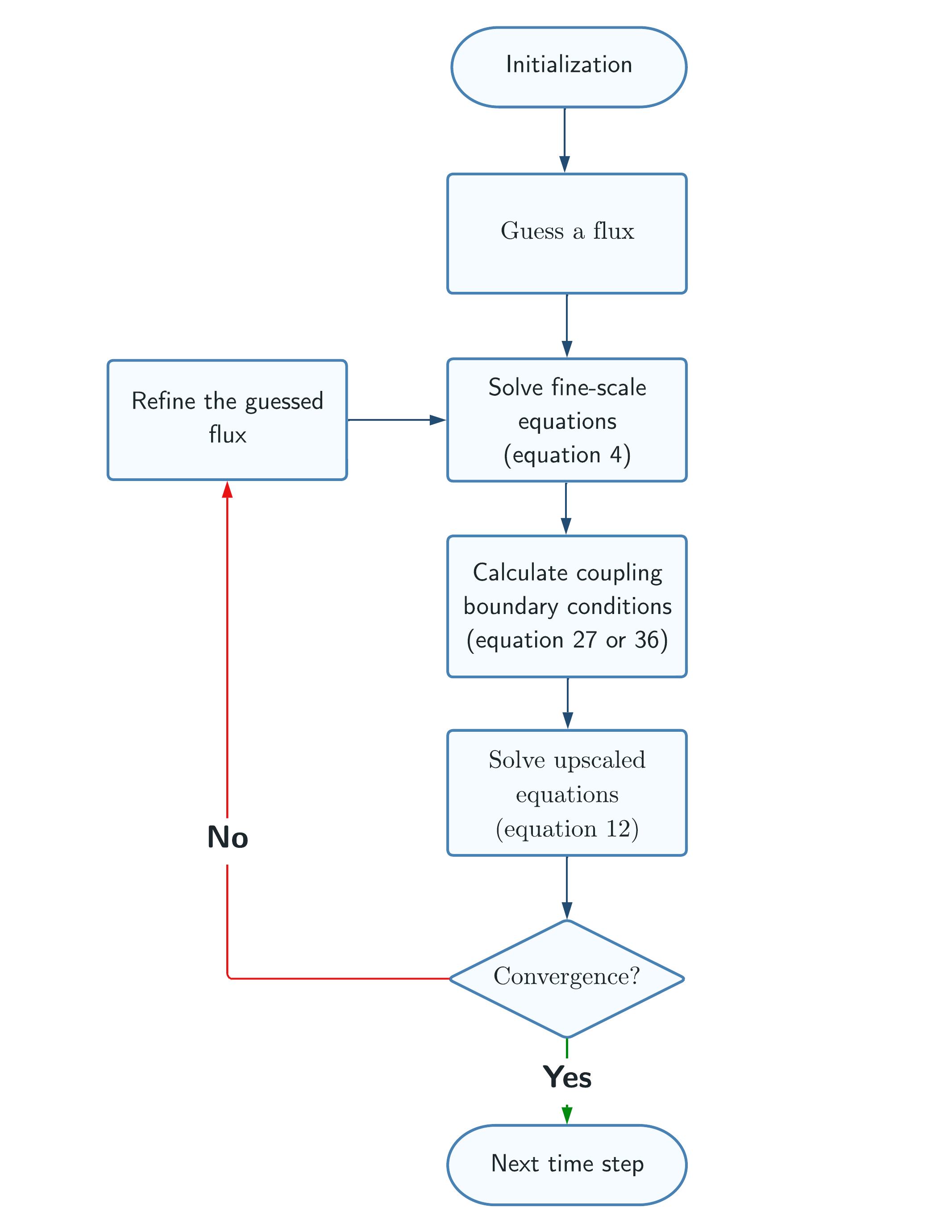}}}
\caption{Flow chart of the hybrid algorithm. }
\label{fig:hybrid_alg}
\end{figure}

Overall, the hybrid simulation is given by the following steps (Figure~\ref{fig:hybrid_alg}):

\begin{enumerate}
    \item Solve fine-scale equations with {a guess for} flux $q^{\left(p,n\right)}$ from the previous time step 
    \item Calculate average packing materials flux with either Taylor (equations~\eqref{eq:taylor-hc-flux}) or Series (equations~\ref{eq:series-hc-flux}) expansion
    \item Solve upscaled equations with calculated average packing material's flux
    \item Compute average packing temperature with either Taylor (equation~\eqref{eq:taylor-hc-temp}) or Series (equation~\eqref{eq:series-hc-temp}) expansion
    \item Compute the error of average temperature continuity $\mathcal{F}$ at the coupling boundary 
    \begin{align}
    & \mathcal{F} = \ensmean{T^{(p)}} - \ensmean{T_\epsilon^{\left(p\right)}}, 
    \end{align}
    
    \item If $\max(\norm{{\mathcal{F}}}_\infty, \norm{{\mathcal{F}}}_2)  > \epsilon_{tol}$ where $\epsilon_{tol}$ is a defined tolerance, refine the unresolved flux $q^{\left(p,n\right)}$ with a zero-finding algorithm (i.e, Broyden's method) and repeat Step 1-5.
\end{enumerate}

%% file: thermal.tex
\section{Thermal runaway problems}
\label{sec:therm-res}
\subsection{Test case description}
\begin{figure}[hbpt!]
\centerline{
 {\includegraphics[width=\textwidth]{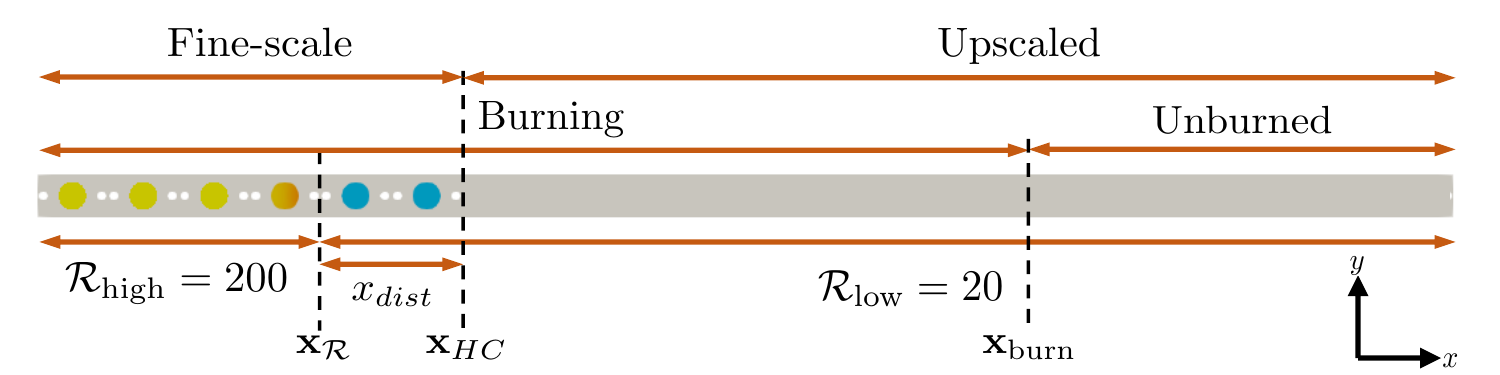}}}

\caption{Schematic diagram of the setup used in the simulations.}
\label{fig:setup-diag}
\end{figure}
To validate the accuracy of the hybrid coupling method developed, we consider a two-dimensional domain with $N_{x}^{\left(c\right)}=20$ and $N_{y}^{\left(c\right)}=1$ (Figure~\ref{fig:setup-diag}). Table~\ref{tab:unit-cell-geom-params} summarizes the parameters used to define the geometry of the unit cell. Periodic boundary conditions are applied at the top and bottom, whereas zero gradient boundary conditions are applied to the left and right boundaries of the domain. To model a battery thermal runaway problem, we modify the power flux source term in the fine-scale equations as 
\begin{align}
\Pi\left(T_{\epsilon}^{\left(c\right)}, \mathbf{x}\right) &=
\begin{cases}
\Pi_{\text{NB}}\left(T_{\epsilon}^{\left(c\right)}\right), \text{ for } \mathbf{x} > \mathbf{x}_{\text{burn}}, \\
\Pi_{\text{FB}}\left(T_{\epsilon}^{\left(c\right)}\right), \text{ for } \mathbf{x} \le \mathbf{x}_{\text{burn}},
\end{cases}
\end{align}
where $\mathbf{x}_{\text{burn}}$ is the location that separates the burning and unburned regions, %\textcolor{red}{[Check which Pi is burning and unburning (i.e., NB vs. FB). It is inconsistent in the following (not sure if it is throughout the manuscript as well).]} 
$\Pi_{\text{NB}}\left(T_{\epsilon}^{\left(c\right)}\right)$ and $\Pi_{\text{FB}}\left(T_{\epsilon}^{\left(c\right)}\right)$ are the dimensionless fine-scale power flux source term for unburned and burning battery cell such that
\begin{subequations}
\begin{align}
\Pi_{\text{NB}}\left(T_{\epsilon}^{\left(c\right)}\right) &= \Pi_{\text{base}} + \frac{1}{2}\left\{\text{Erf}\left[A_1 T_{\epsilon}^{\left(c\right)} + B_1\right] + 1\right\}\left(1 - \Pi_{\text{base}}\right) \nonumber \\ 
& - \frac{1}{2}\left\{\text{Erf}\left[A_2 T_{\epsilon}^{\left(c\right)} + B_2\right] + 1\right\}\\
\Pi_{\text{FB}}\left(T_{\epsilon}^{\left(c\right)}\right) &= \Pi_{\text{burn}} - \frac{1}{2}\left\{\text{Erf}\left[A_2 T_{\epsilon}^{\left(c\right)} + B_2\right] + 1\right\}.\label{eq:new_pi_term_dimless_burn}
\end{align}
\end{subequations}
For the upscaled equations, we use a sigmoid function to calculate the power flux source term to avoid discontinuity such that
\begin{align}
\overline{\Pi} \left(\langle T^{\left(c\right)} \rangle_{Y}, \mathbf{x}\right) &= \ddfrac{ \overline{\Pi}_{\text{FB}} \left(\langle T^{\left(c\right)} \rangle_{Y}, \mathbf{x}\right) -  \overline{\Pi}_{\text{NB}} \left(\langle T^{\left(c\right)} \rangle_{Y}, \mathbf{x}\right)}{1+\exp\left(\gamma(\mathbf{x} - 
\mathbf{x}_{\text{burn}})\right)} + \overline{\Pi}_{\text{NB}} \left(\langle T^{\left(c\right)} \rangle_{Y}, \mathbf{x}\right), \label{eq:upscale_force_pi_term}
\end{align}
{where $\gamma=180$ is a constant to approximate a stepwise function}, $\overline{\Pi}_{\text{NB}} \left(\langle T^{\left(c\right)} \rangle_{Y}, \mathbf{x}\right)$ and $\overline{\Pi}_{\text{FB}} \left(\langle T^{\left(c\right)} \rangle_{Y}, \mathbf{x}\right)$ are the dimensionless upscaled power flux source term for unburned and burning battery cells that are defined as
\begin{subequations}
\begin{align}
\overline{\Pi}_{\text{NB}} \left(\langle T^{\left(c\right)} \rangle_{Y}, \mathbf{x}\right) &= \Pi_{\text{base}}\left(\mathbf{x}\right) + \frac{1}{2}\left\{\text{Erf}\left[\frac{A_1}{\phi^{\left(c\right)}} \langle T^{\left(c\right)} \rangle_{Y} + B_1\right] + 1\right\}\left(1 - \Pi_{\text{base}}\left(\mathbf{x}\right)\right)  \nonumber \\
&- \frac{1}{2}\left\{\text{Erf}\left[\frac{A_2}{\phi^{\left(c\right)}} \langle T^{\left(c\right)} \rangle_{Y} + B_2\right] + 1\right\}, \\
\overline{\Pi}_{\text{FB}} \left(\langle T^{\left(c\right)} \rangle_{Y}, \mathbf{x}\right) &= \Pi_{\text{burn}}\left(\mathbf{x}\right) - \frac{1}{2}\left\{\text{Erf}\left[\frac{A_2}{\phi^{\left(c\right)}} \langle T^{\left(c\right)} \rangle_{Y} + B_2\right] + 1\right\}. 
\end{align}
\end{subequations}
Additionally, in reality, the rate of heat generation for burning battery cells can vary from cell to cell because of potential manufacturing defects and aging. Therefore, we consider that a region of burning battery cells has different rates of heat generation
where 
\begin{align}
    \hat{\Pi}_{\text{burn}} = 
    \begin{cases}
    \hat{\Pi}_{\text{burn},\text{low}}=\hat{\Pi}_{\text{burn},0}, \text{ for } \mathbf{x} > \mathbf{x}_{\mathcal{R}}, \\
    \hat{\Pi}_{\text{burn},\text{high}}=10\hat{\Pi}_{\text{burn},0}, \text{ for } \mathbf{x} \le \mathbf{x}_{\mathcal{R}},
    \end{cases}
\end{align}
such that
\begin{align}
    &\hat{\Pi}_{\text{burn},0} = \ddfrac{\hat{T}_{max} \hat{k}^{\left(p\right)}}{\max{(\hat{\mathcal{L}}_{x}, \hat{\mathcal{L}}_{y})}\hat{\ell}}.
\end{align}
Incorporating the effect of $\hat{\Pi}_{\text{burn},\text{low}}$ and $\hat{\Pi}_{\text{burn},\text{high}}$, we define the dimensionless number $\mathcal{R}$ as 
\begin{align}
&\mathcal{R}(\mathbf{x}) = \left(\ddfrac{\mathcal{R}_{\text{high}} + \mathcal{R}_{\text{low}} }{2}\right) - \left(\ddfrac{\mathcal{R}_{\text{high}} - \mathcal{R}_{\text{low}} }{2}\right)\tanh\left(100(x-\mathbf{x}_{\mathcal{R}}) \right),
\label{eq:Rx-def}
\end{align}
where $\mathcal{R}_{\text{high}} = 10\mathcal{R}_{0}$ and $\mathcal{R}_{\text{low}} = \mathcal{R}_{0}$ represent regions with high and low heat generation rates, respectively, and $\mathcal{R}_{0}$ is defined with $\hat{\Pi}_{\text{burn},0}$ in equation~\eqref{eq:dimless_groups}. By setting $\mathbf{x}_{\text{burn}}$ and $\mathbf{x}_{\mathcal{R}}$, two temperature gradients will be observed due to (1) the difference in the rate of heat generation by burned battery cells at $\mathbf{x}_{\mathcal{R}}$ and (2) the difference in the rate of heat generation by burned and unburned battery cells at $\mathbf{x}_{\text{burn}}$. 
\begin{table}[ht!]
\centering
\caption{
  \label{tab:unit-cell-geom-params} Values of parameters used to define unit-cell geometry in the simulations.}
\begin{tabular}{l|c}
\hline \hline
Cell radius $\hat{r}_\epsilon^{\left(c\right)}$ [\si{L}]  &  0.009  \\
Pipe radius $\hat{r}_\epsilon^{\left(w\right)}$ [\si{L}]  &  0.003  \\
Distance between two battery cells $d^{\left(cc\right)}_\epsilon$ [\si{L}] & 0.009 \\
Distance between battery cell and cooling pipe $d^{\left(1\right)}_\epsilon$ [\si{L}] & 0.001\\
Distance between battery cell and cooling pipe $d^{\left(2\right)}_\epsilon$ [\si{L}] & 0.002 \\
\hline \hline
\end{tabular}%
\end{table}
% The main reason is related to the switching from diffusive to advective behavior, violating the homogeneous assumption and invalidating the upscaled governing equations. Pietrzyk et al.~\cite{Pietrzyk2021-lu} observed the heat front propagate faster in the upscaled than the fine-scale simulations (Figure~\ref{fig:cell_temp_1271_acc}(c) vs (d)). 
In the following simulations, we define an arbitrary simulation end time $t_{final} = 0.2 \hat{t}$ that is sufficient for the temperature to achieve equilibrium, where
\begin{align}
    \hat{t} &= \ddfrac{\hat{\rho}^{\left(p\right)}\hat{C}^{\left(p\right)}\left(\max{(\hat{\mathcal{L}}_{x}, \hat{\mathcal{L}}_{y})}\right)^2}{\hat{k}^{\left(p\right)}},
\end{align}
is the reference timescale for heat transfer. Table~\ref{tab:ref-val-scaling} summarizes the values of the reference parameters used to define the dimensionless numbers in equations~\eqref{eq:dimless_groups} for all subsequent simulations. We further define the following parameters as
\begin{align}
    &\Pi_{\text{base},0} = 0.01\Pi_{\text{burn},0}, \\
    &\hat{U}^{\left(pc\right)} = \ddfrac{\hat{k}^{\left(p\right)}}{\max{(\hat{\mathcal{L}}_{x}, \hat{\mathcal{L}}_{y})}}, \\
    &\hat{Q}^{\left(pw\right)} = 0.00001\ddfrac{\hat{T}_{max} \hat{k}^{\left(p\right)}}{\max{(\hat{\mathcal{L}}_{x}, \hat{\mathcal{L}}_{y})}},
\end{align}
to ensure the applicability of the upscaled equations. The simulations are initialized with zeros for both the packing- and cell-temperature fields. 

% \noindent where $i=p$ or $c$ refers to packing material or the battery cell, repsectively,  $\hat{\rho}^{\left(i\right)}$ [\si{ML\tothe{-3}}] is the density, $\hat{C}^{\left(i\right)}$ [\si{L\tothe{2}T\tothe{-2}\Theta\tothe{-1}}]is the heat capacity, $\hat{T}_{\epsilon}^{\left(i\right)} [\si{\Theta}] \equiv \hat{T}_{\epsilon}^{\left(i\right)}\left(\hat{t}, \hat{\mathbf{x}}\right)$  is the temperature at time $\hat{t} > 0$ and location $\hat{\mathbf{x}} \in \hat{\mathcal{B}}_{\epsilon}^{\left(i\right)}$, $\hat{k}^{\left(i\right)}$ [\si{MLT\tothe{-3}\Theta\tothe{-1}}] is the thermal conductivity, $\n_{\epsilon}^{\left(i\right)} \equiv \n_{\epsilon}^{\left(i\right)}\left(\hat{\mathbf{x}}\right)$ is the normal vector to the interfaces pointed away from the domain, $\hat{U}^{\left(pc\right)}$ [\si{MT\tothe{-3}\Theta\tothe{-1}}] is the total heat transfer coefficient between the packing material and battery cells,  $\hat{q}_{\epsilon}^{\left(pw\right)}\left(\hat{t}, \hat{\mathbf{x}}\right)$ [\si{MT\tothe{-3}}]  is a power flux between the packing material and the cooling water pipes, and $\hat{\Pi}(\hat{t}, \hat{\mathbf{x}})$ [\si{ML\tothe{-1}T\tothe{-3}}] is a power flux source term.

\begin{table}[ht!]
\centering
\caption{
  \label{tab:ref-val-scaling} Values of reference values to nondimensionalize the fine-scale governing equations.}
\begin{tabular}{l|c}
\hline \hline
Density of battery cell $\hat{\rho}^{\left(c\right)}$ [\si{ML\tothe{-3}}]   &  2500  \\
Density of packing materials $\hat{\rho}^{\left(p\right)}$ [\si{ML\tothe{-3}}]   &  1500\\
Heat capacity of battery cell $\hat{C}^{\left(c\right)}$ [\si{L\tothe{2}T\tothe{-2}\Theta\tothe{-1}}]  &  900 \\
Heat capacity of packing materials $\hat{C}^{\left(p\right)}$ [\si{L\tothe{2}T\tothe{-2}\Theta\tothe{-1}}]   & 1500\\
Thermal conductivity of battery cell $\hat{k}^{\left(c\right)}$ [\si{MLT\tothe{-3}\Theta\tothe{-1}}] & 3 \\
Thermal conductivity of packing materials $\hat{k}^{\left(p\right)}$ [\si{MLT\tothe{-3}\Theta\tothe{-1}}] & 3\\
Defined temperature at the battery pack edges $\hat{T}_\infty$ [\si{\Theta}] & 293\\
Temperature ranges over which $\hat{\Pi}(\hat{T}_{\epsilon}^{\left(c\right)},\hat{\mathbf{x}}) = \hat{\Pi}_{\text{base}}(\hat{\mathbf{x}})$, $\hat{T}_a$ [\si{\Theta}] & 0\\
Temperature ranges over which $\hat{\Pi}(\hat{T}_{\epsilon}^{\left(c\right)}, \hat{\mathbf{x}}) = \hat{\Pi}_{\text{burn}}$,  $\hat{T}_b$ [\si{\Theta}] & 0\\
Temperature ranges between $\hat{\Pi}_{\text{base}}(\hat{\mathbf{x}})$ and $\hat{\Pi}_{\text{burn}}$, $\hat{T}_{s1}$ [\si{\Theta}] & 120\\
Temperature ranges between $\hat{\Pi}_{\text{burn}}$ and $0$, $\hat{T}_{s2}$ [\si{\Theta}] & 120\\
Smoothness of the error functions $\epsilon_{s1}$ [-] & 0.0005\\
Smoothness of the error functions $\epsilon_{s2}$ [-] & 0.0005\\
\hline \hline
\end{tabular}%
\end{table}

% {over which} $\hat{\Pi}(\hat{T}_{\epsilon}^{\left(c\right)},\hat{\mathbf{x}}) = \hat{\Pi}_{\text{base}}(\hat{\mathbf{x}})$ and $\hat{\Pi}(\hat{T}_{\epsilon}^{\left(c\right)}, \hat{\mathbf{x}}) = \hat{\Pi}_{\text{burn}}$

% $\hat{\Pi}_{\text{base}}(\hat{\mathbf{x}})$ and $\hat{\Pi}_{\text{burn}}$ are the \textit{base} and \textit{burn} power flux values, respectively, $\hat{T}_{ref}$~[\si{\Theta}] is the reference temperature, $\hat{T}_a$~[\si{\Theta}] and $\hat{T}_b$~[\si{\Theta}] are temperature ranges over with $\hat{\Pi}(\hat{T}_{\epsilon}^{\left(c\right)},\hat{\mathbf{x}}) = \hat{\Pi}_{\text{base}}(\hat{\mathbf{x}})$ and $\hat{\Pi}(\hat{T}_{\epsilon}^{\left(c\right)}, \hat{\mathbf{x}}) = \hat{\Pi}_{\text{burn}}$, respectively, $\hat{T}_{s1}$ and $\hat{T}_{s2}$ are the temperature ranges which $\hat{\Pi}(\hat{T}_{\epsilon}^{\left(c\right)}, \hat{\mathbf{x}})$ transitions from $\hat{\Pi}_{\text{base}}(\hat{\mathbf{x}})$ to $\hat{\Pi}_{\text{burn}}$ and from $\hat{\Pi}_{\text{burn}}$ to $0$, respectively, and $\epsilon_{s1}$~[-]$= 0.0005$ and $\epsilon_{s2}$~[-]$= 0.0005$ are paramters associated with smoothness of the error functions. The detailed formulation and validation of the source term can be found in Pietrzyk et al.~\cite{Pietrzyk2022-ou}.

\subsection{Accuracy of hybrid coupling}
\label{sec:acc-hc}
Following  Pietrzyk et al. (2022), the error between fine-scale and upscaled simulations should be bounded by the theoretical upscaling error $\mathcal{O}\left(\epsilon\right)$ where $\epsilon = 1/\max{(N_{x}^{\left(c\right)}, N_{y}^{\left(c\right)})} = 0.05 \ll 1$ and the magnitudes of dimensionless numbers are within the applicability regime. Since the coupling boundary conditions are derived with an error in $\mathcal{O}(\epsilon)$, hybrid simulation errors are always expected to be within the upscaling error threshold. To validate the accuracy of the developed coupling method, we compare the results of hybrid simulations with fine-scale and upscaled simulations. Here, we define the coupling boundary at $x_{HC}=-0.0875$, which lies between $\mathbf{x}_{\mathcal{R}}=-0.3125$ ($4.5$ unit cells away) and $\mathbf{x}_{\text{burn}}=0.2125$ ($6$ unit cells away). {To automate the detection of $\mathbf{x}_{\mathcal{R}}$ with a distribution of $\mathcal{R}$ (Equation~\eqref{eq:Rx-def}, we compute $x_{\mathcal{R}}$ as}
\begin{align}
    {\mathbf{x}_{\mathcal{R}} =  \argmax_{\mathbf{x}} \left( \left\vert \ddfrac{\mathcal{R}}{\mathcal{R}_{\text{low}}} - 1 - \alpha_1  \right\vert \right),}
\end{align}
{where $\alpha_{1} = 0.01$ is the tolerance factor.} This ensures the validity of the homogeneity assumption at the coupling boundary until the end of the simulations. Table~\ref{tab:val-acc-sim-params} summarizes the parameters used in the simulations.

\begin{figure}
\centerline{
 {\includegraphics{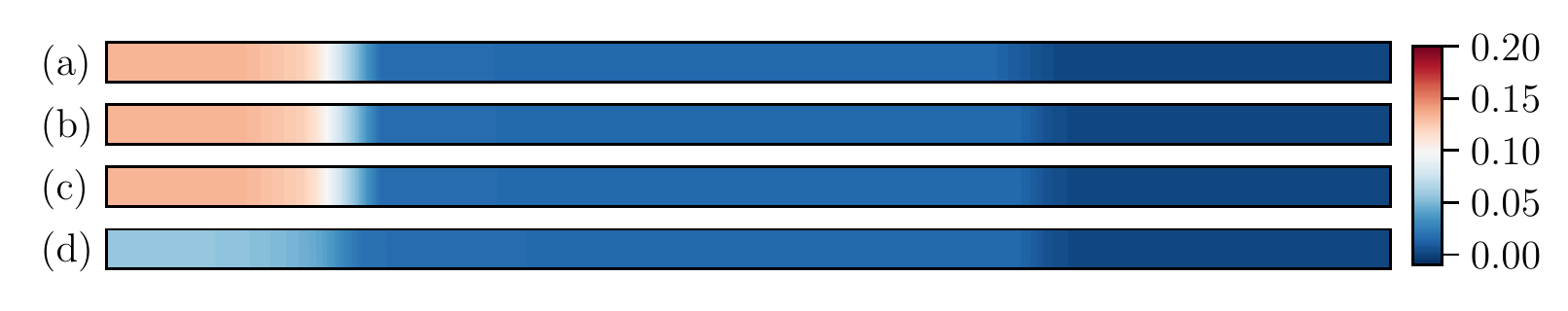}}}

\caption{Average cell temperature computed at $t=0.02$ with a time step size $\Delta{t}=3.15 \times 10^{-5}$ for (a) fine-scale (b) hybrid with Taylor expansion, (c) hybrid with Series expansion and (d) upscaled simulations.}
\label{fig:cell_temp_032_acc}
\end{figure}

\begin{figure}
\centerline{
 {\includegraphics{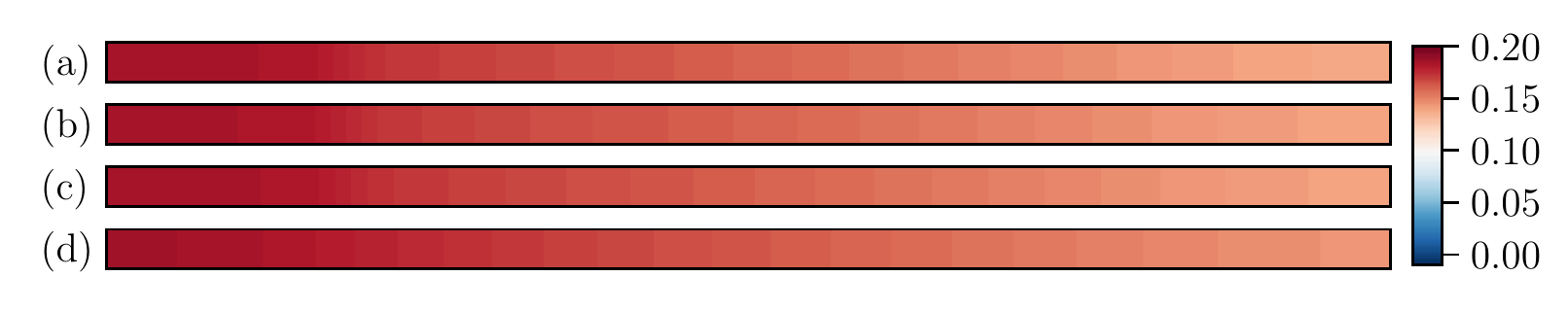}}}

\caption{Average cell temperature computed at $t=0.2$ with a time step size $\Delta{t}=3.15 \times 10^{-5}$ for (a) fine-scale (b) hybrid with Taylor expansion (c) hybrid with Series expansion and (d) upscaled simulations.}
\label{fig:cell_temp_1271_acc}
\end{figure}

\begin{figure}
\centerline{
 {\includegraphics{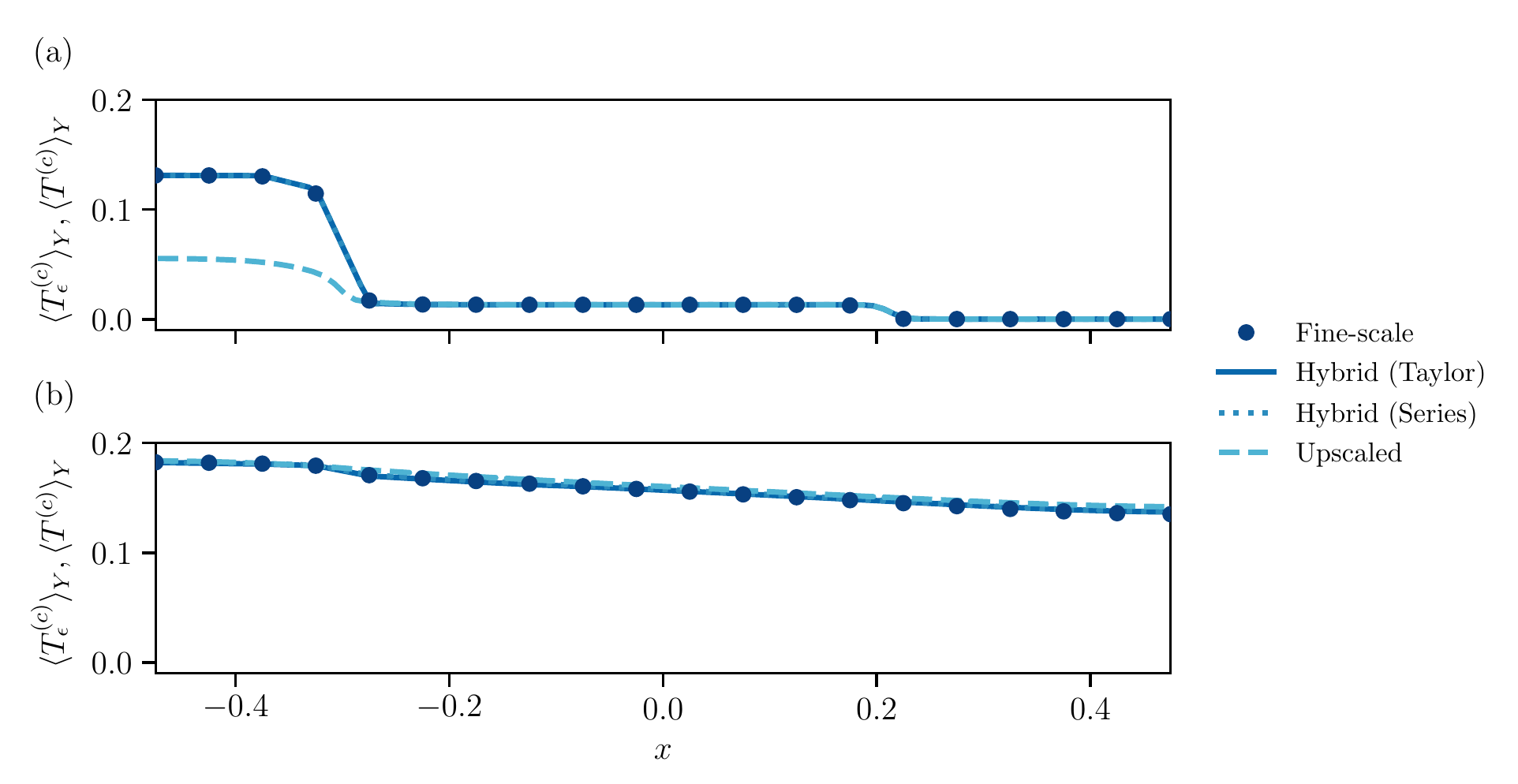}}}

\caption{The centerline ($y=0$) average cell temperature computed at (a) $t=0.02$ and (b) $t=0.20$ with a time step size $\Delta{t}=3.15 \times 10^{-5}$ for the fine-scale, hybrid with Taylor and Series expasnion, respectively, and upscaled simulations.}
\label{fig:cell_temp_line_acc}
\end{figure}

\begin{table}[ht!]
\centering
\caption{
  \label{tab:val-acc-sim-params} Summary of simulation parameters used in evaluating the accuracy of the proposed hybrid algorithm.}
\begin{tabular}{l|c}
\hline \hline
Time step size $\Delta{t}/\hat{t}$ &  $3.15 \times 10^{-5}$  \\
Upscale subdomain minimum grid resolution $h_{\text{up},min}$  &  $1.00 \times 10^{-2}$\\
fine-scale subdomain minimum grid resolution $h_{\text{fine},min}$  &  $7.63\times10^{-4}$\\
Final simulation time $t_{final}$   & 0.2\\
Coupling boundary location $x_{HC}$ & -0.0875 \\
$\mathcal{R}$ location $\mathbf{x}_{\mathcal{R}}$  & -0.3125 \\
burned location $\mathbf{x}_{\text{burn}}$  & 0.2125 \\
Polynomial order $k$ & 1 \\
\hline \hline
\end{tabular}%
\end{table}
To compare fine-scale and hybrid simulations with upscaled simulations, we apply an averaging operator (equation~\eqref{eq:ave-op}) to the fine-scale and hybrid simulation results. Figures~\ref{fig:cell_temp_032_acc} and~\ref{fig:cell_temp_1271_acc} show the contour plots of the average cell temperatures for fine-scale, hybrid (Taylor and Series expansion methods), and upscaled simulations at $t=0.02$ and $0.20$ respectively. At $t=0.02$, we clearly observe that upscaled simulations (Figure~\ref{fig:cell_temp_032_acc}(d)) underpredict the average cell temperature for $\mathbf{x} \le \mathbf{x}_{\mathcal{R}}$, as indicated by the color difference. No significant differences in the average cell temperature are observed near the burned location $\mathbf{x}_{\text{burn}}$. This shows that the upscaled equations can accurately capture the behavior of heat wave propagation due to burned and unburned cells. However, since the value of $\mathcal{R}$ is greater than the applicability regime for $\mathbf{x} \le \mathbf{x}_{\mathcal{R}}$ ($FT \sim \mathcal{O}(\epsilon)$), the upscaled simulations are unable to capture the rate of increase in cell temperature. In contrast to the upscaled simulations, hybrid simulations (Figure~\ref{fig:cell_temp_032_acc}(b) and (c)) are capable of simulating accurate heat propagation behaviors at both $\mathbf{x}_{\mathcal{R}}$ and $\mathbf{x}_{\text{burn}}$, therefore obtaining cell temperatures that are comparable to simulations at the fine scale. Similar trends can be observed for the packing temperature, where the upscaled simulations overpredict the temperatures (\ref{app:pack-temp-appen}). By comparing the centerline average cell temperature ($y=0$) in Figure~\ref{fig:cell_temp_line_acc}, the hybrid simulations have clearly demonstrated their ability to capture the accurate behavior of the propagation of heat waves compared to upscaled simulations. No obvious difference between Taylor and Series expansion approaches has been observed, as expected, because the accuracy of the upscaled equations is only first order. Therefore, the first-order accurate Taylor approach and the second-order accurate Series approach should not result in significant differences. At $t = 0.2$, the cell temperature approach the equilibrium temperature; therefore, the effect of large $\mathcal{R}$ values is no longer significant. Both the hybrid and upscaled simulations are accurate in predicting the average cell temperature (Figure~\ref{fig:cell_temp_1271_acc} and~\ref{fig:cell_temp_line_acc}(b)). Similar observations on the average temperature of the packing materials can be found in Figure~\ref{fig:cell_temp_1271_acc} of~\ref{app:pack-temp-appen}. 

\begin{figure}
\centerline{
 {\includegraphics{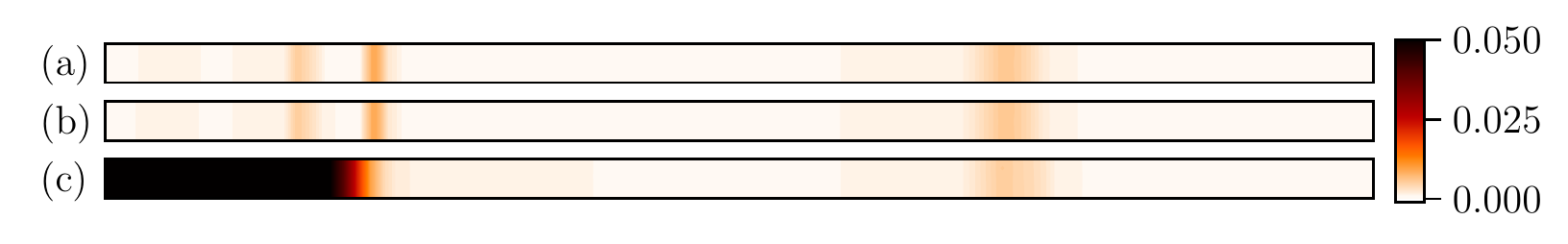}}}

\caption{Absolute error of average cell temperature computed at $t=0.02$ with a time step size $\Delta{t}=3.15 \times 10^{-5}$ for (a) fine-scale (b) hybrid with Taylor expansion, (c) hybrid with Series expansion and (d) upscaled simulations.}
\label{fig:cell_temp_err_032_acc}
\end{figure}

\begin{figure}
\centerline{
 {\includegraphics{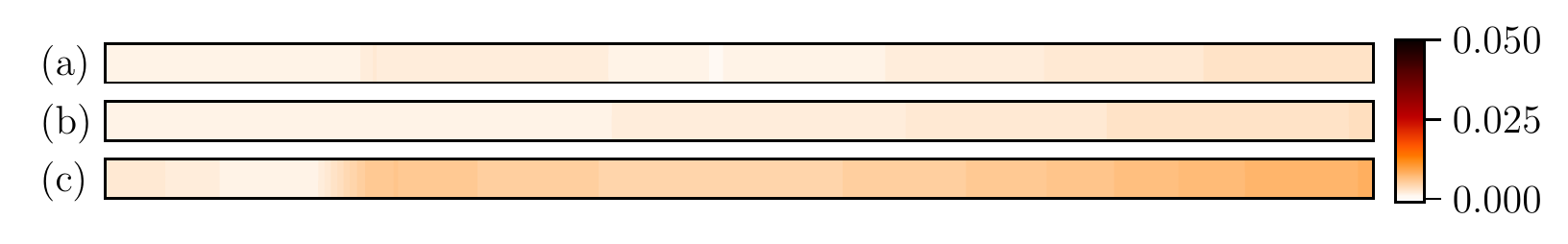}}}

\caption{Absolute error of the average cell temperature computed at $t=0.20$ with a time step size $\Delta{t}=3.15 \times 10^{-5}$ for (a) fine-scale (b) hybrid with Taylor expansion, (c) hybrid with Series expansion and (d) upscaled simulations.}
\label{fig:cell_temp_err_1271_acc}
\end{figure}

\begin{figure}
\centerline{
 {\includegraphics{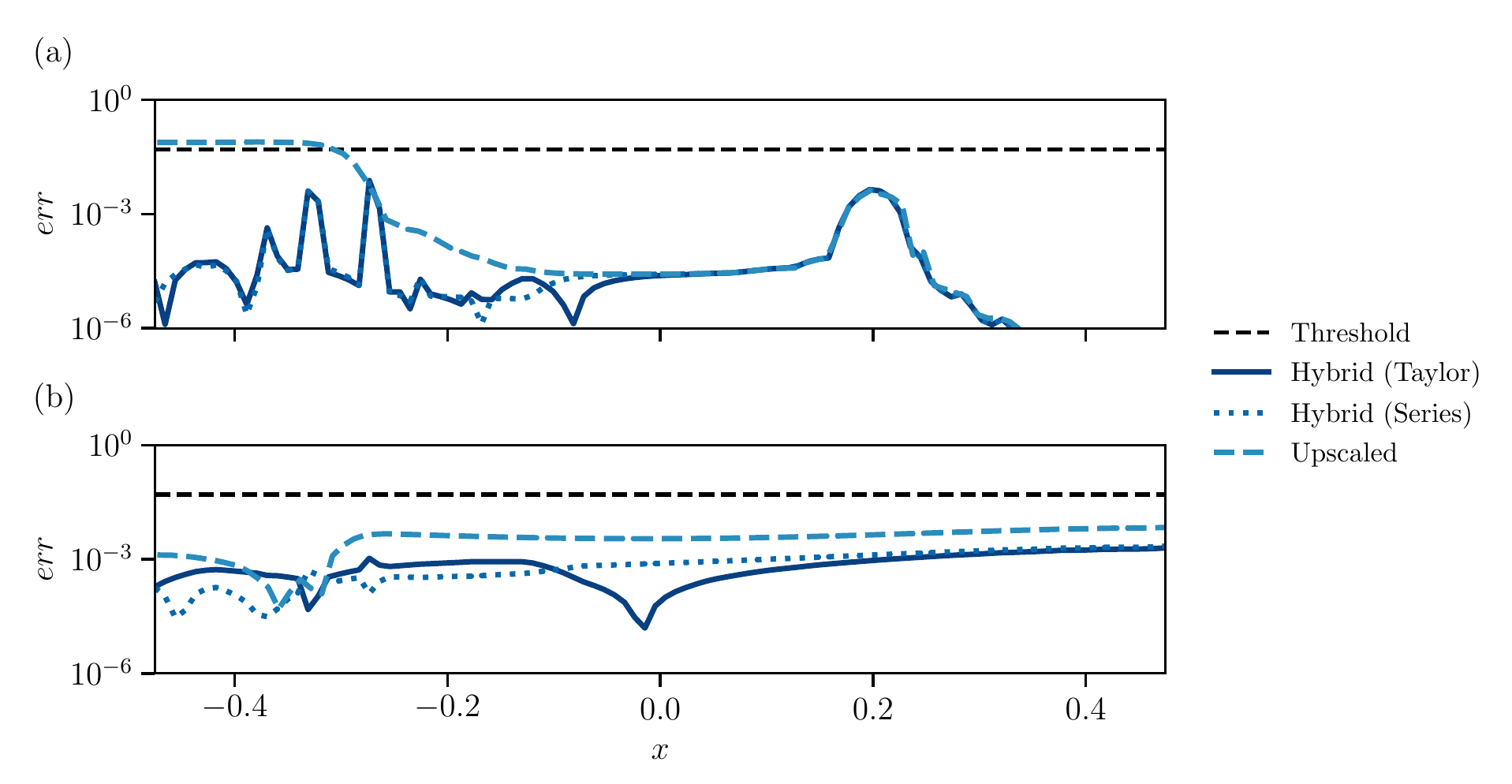}}}
\caption{The centerline ($y=0$) absolute error of the average cell temperature computed at (a) $t=0.02$ and (b) $t=0.20$ with a time step size $\Delta{t}=3.15 \times 10^{-5}$ for hybrid with Taylor and Series expasnion, respectively, and upscaled simulations.}
\label{fig:cell_temp_err_line_acc}
\end{figure}

We now compute the error as 
\begin{align}
\label{eq:err-cal}
    &err = \abs{\ensmean{T^{\left(i\right)}}_Y - \ensmean{T^{\left(i\right)}_\epsilon}_Y},
\end{align}
where $i=c$ or $p$ represent the cell or packing temperature, respectively. Figures~\ref{fig:cell_temp_err_032_acc} and~\ref{fig:cell_temp_err_1271_acc} show the contour plot of the errors at two different time instances. At $t=0.02$ where the battery cells begin to release heat, the errors in the upscaled simulations quickly approach the tolerance, as indicated by the black region for $\mathbf{x} \le \mathbf{x}_{\mathcal{R}}$. On the contrary, the errors in hybrid simulations are still below the tolerance. At $t=0.2$, the errors for both hybrid and upscaled simulations are lower than the threshold $\epsilon=0.05$ because the cell temperature approaches equilibrium. Figure~\ref{fig:cell_temp_err_line_acc} shows the centerline error plot of hybrid and upscaled simulations. The error in the upscaled simulations is more than twice the tolerance. No discernible differences have been observed between the Taylor and Series approaches in the hybrid simulations, as expected. The errors observed at $\mathbf{x}_{\text{burn}}$ are due to the difference in the power flux source term where the function is discontinuous in fine-scale simulations, while continuous in upscaled simulations.

\subsection{Coupling boundary location on the accuracy of hybrid coupling}
\label{sec:xhc-acc}
In the previous section, we evaluated the accuracy of hybrid simulations by setting $x_{HC}$ to be approximately 4.5 unit cells from $\mathbf{x}_{\mathcal{R}}$. This conservative setup ensures the accuracy of hybrid simulations at the expense of computational costs. Here, we define the distance between $x_{HC}$ and ${x}_{\mathcal{R}}$ as 
\begin{align}
    &x_{dist} = \abs{x_{HC} - {x}_{\mathcal{R}}}.
\end{align}
In reality, we would want to minimize $x_{dist}$ to avoid additional computational costs. Theoretically, the minimum distance $x_{dist,\min}$ is $1.0\epsilon$ for the Taylor approach and $1.5\epsilon$ for the Series approach. In this section, we focus on evaluating the effects of $x_{dist}$ on the accuracy of hybrid simulations for both approaches where five different distances ($0\epsilon$, $1.0\epsilon$, $1.5\epsilon$, $3.0\epsilon$, and $4.5\epsilon$) will be used. The simulation parameters are identical to Table~\ref{tab:val-acc-sim-params} except for $x_{HC}$. 

\begin{figure}
\centerline{
 {\includegraphics{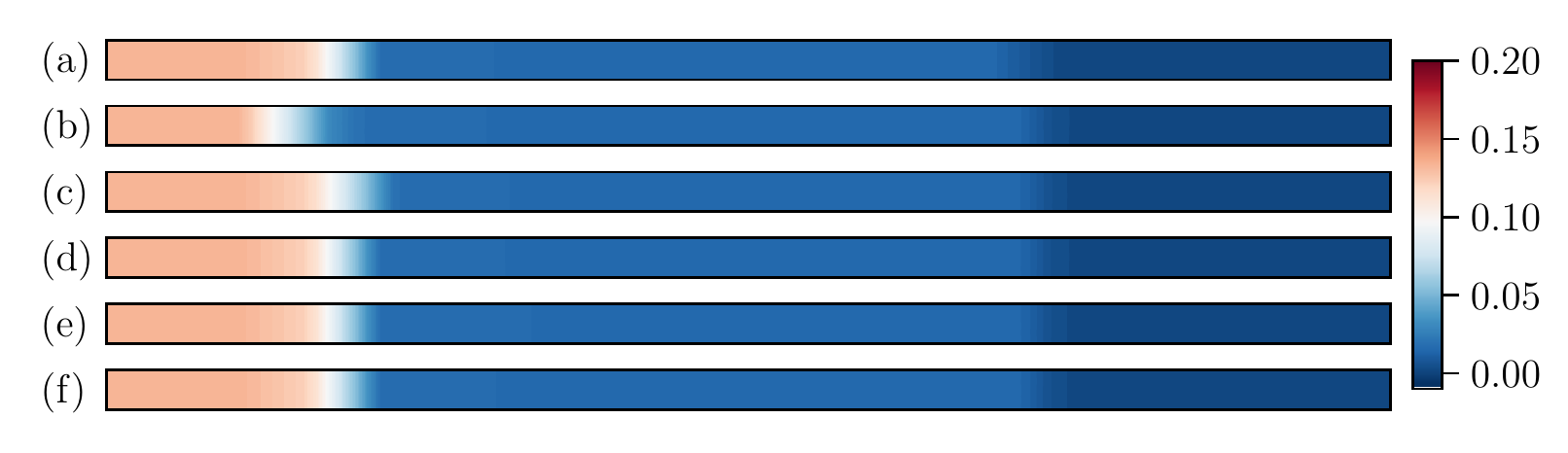}}}

\caption{Average cell temperature computed at $t=0.02$ with a time step size $\Delta{t}=3.15 \times 10^{-5}$ for (a) fine-scale and hybrid simulations with Taylor expansion with $x_{dist}=$ (b) $0.0\epsilon$, (c) $1.0\epsilon$, (d) $1.5\epsilon$, (e) $3.0\epsilon$ and (f) $4.5\epsilon$. }
\label{fig:cell_temp_032_loc_taylor}
\end{figure}

\begin{figure}
\centerline{
 {\includegraphics{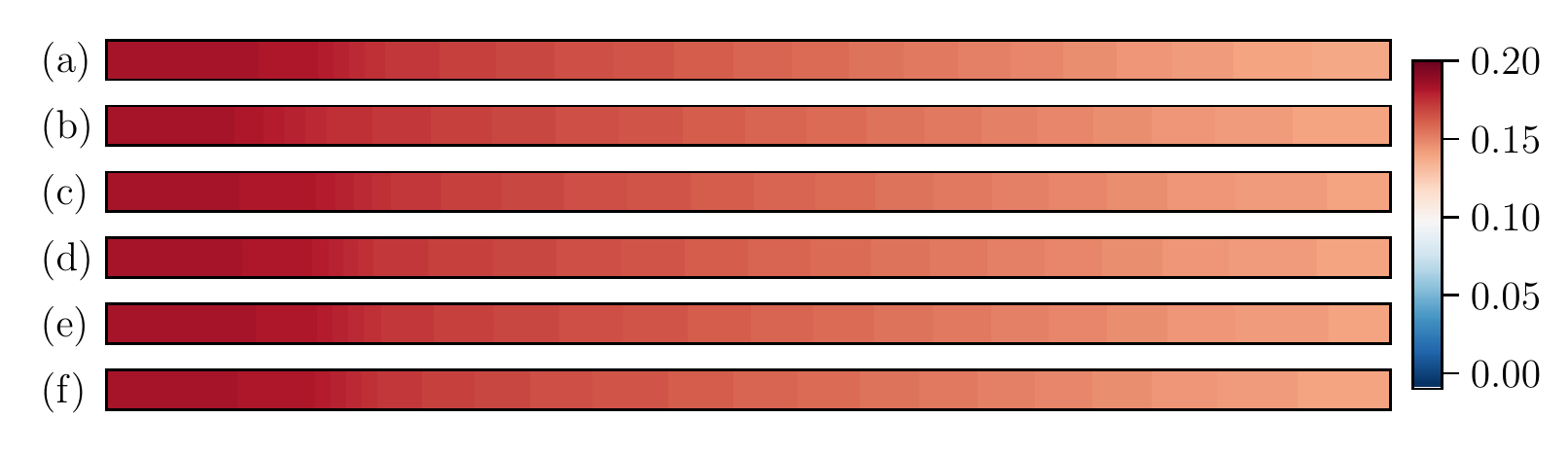}}}

\caption{Average cell temperature computed at $t=0.20$ with a time step size $\Delta{t}=3.15 \times 10^{-5}$ for (a) fine-scale and hybrid simulations with Taylor expansion with $x_{dist}=$ (b) $0.0\epsilon$, (c) $1.0\epsilon$, (d) $1.5\epsilon$, (e) $3.0\epsilon$ and (f) $4.5\epsilon$.}
\label{fig:cell_temp_1271_loc_taylor}
\end{figure}

\begin{figure}
\centerline{
 {\includegraphics{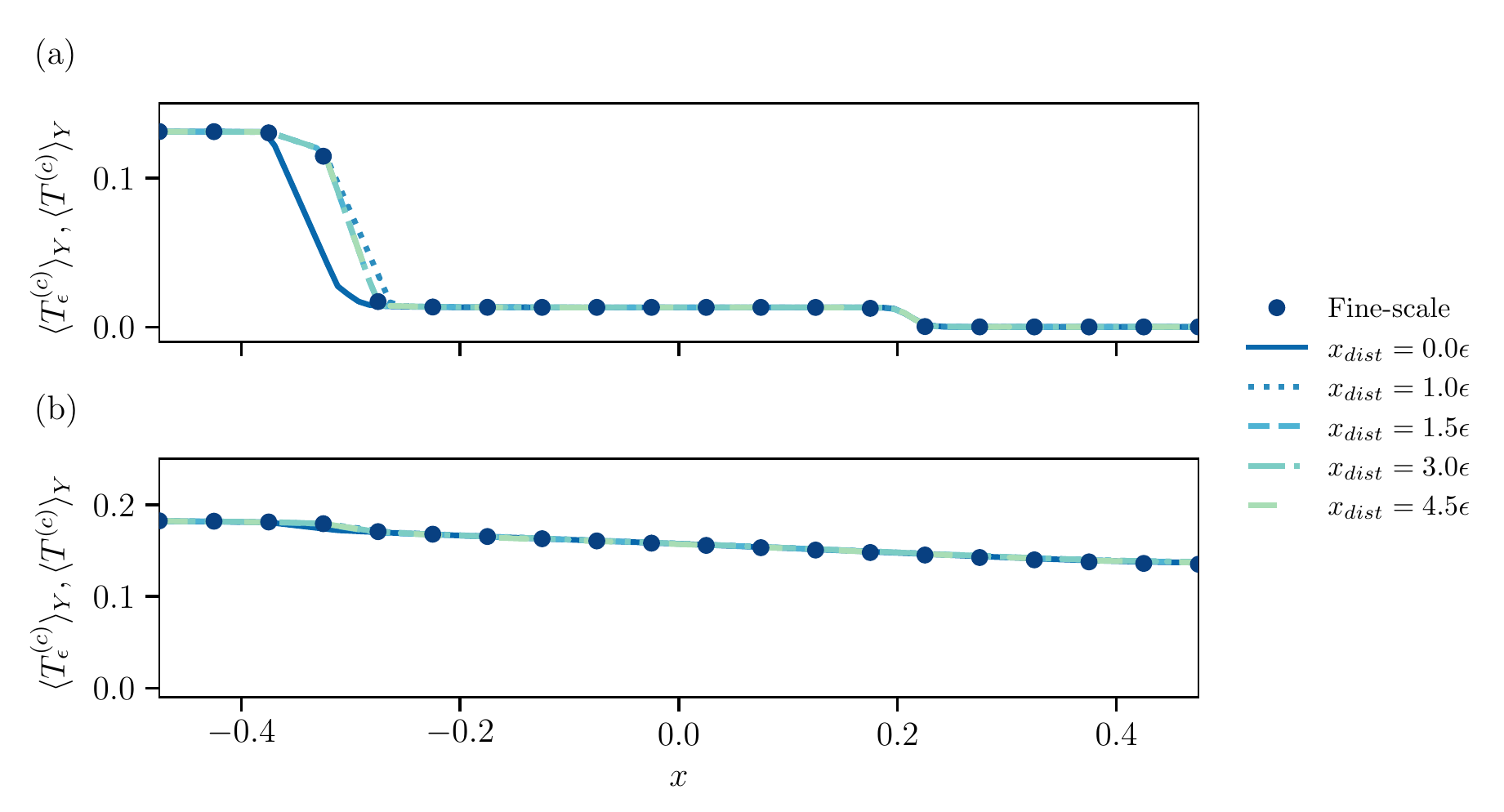}}}

\caption{The centerline ($y=0$) average cell temperature computed at (a) $t=0.02$ and (b) $t=0.20$ with a time step size $\Delta{t}=3.15 \times 10^{-5}$ for fine-scale and hybrid simulations with Taylor expansion with $x_{dist}=0.0\epsilon$, $1.0\epsilon$, $1.5\epsilon$,  $3.0\epsilon$ and $4.5\epsilon$.}
\label{fig:cell_temp_line_loc_taylor}
\end{figure}

\begin{figure}
\centerline{
 {\includegraphics{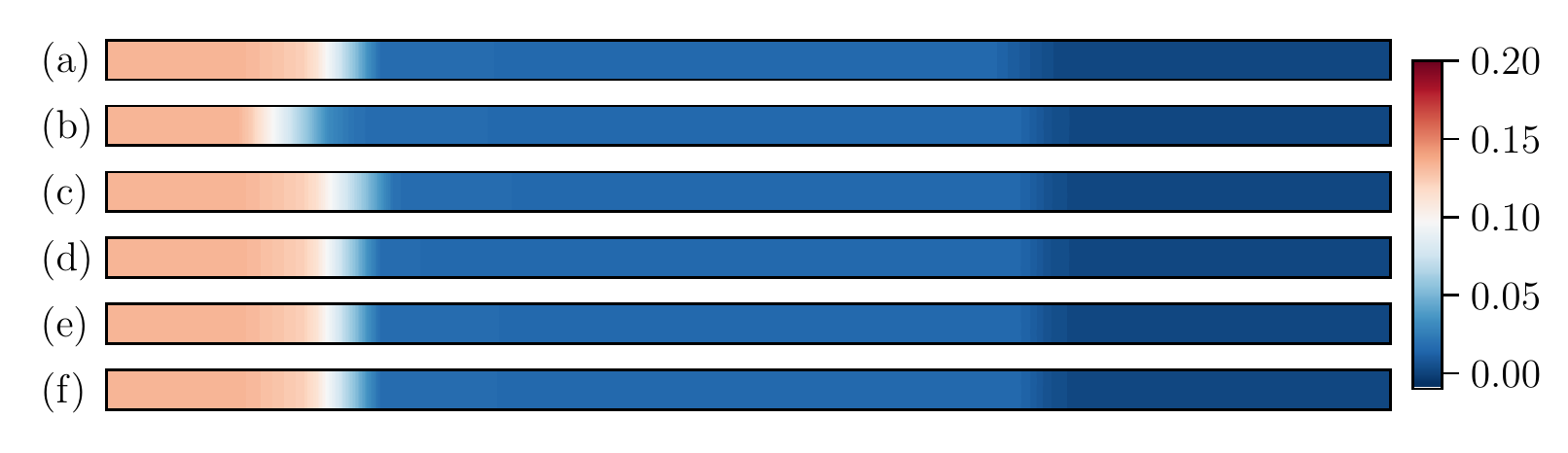}}}
\caption{Average cell temperature computed at $t=0.02$ with a time step size $\Delta{t}=3.15 \times 10^{-5}$ for (a) fine-scale and hybrid simulations with Series expansion with $x_{dist}=$ (b) $0.0\epsilon$, (c) $1.0\epsilon$, (d) $1.5\epsilon$, (e) $3.0\epsilon$ and (f) $4.5\epsilon$.}
\label{fig:cell_temp_032_loc_series}
\end{figure}

\begin{figure}
\centerline{
 {\includegraphics{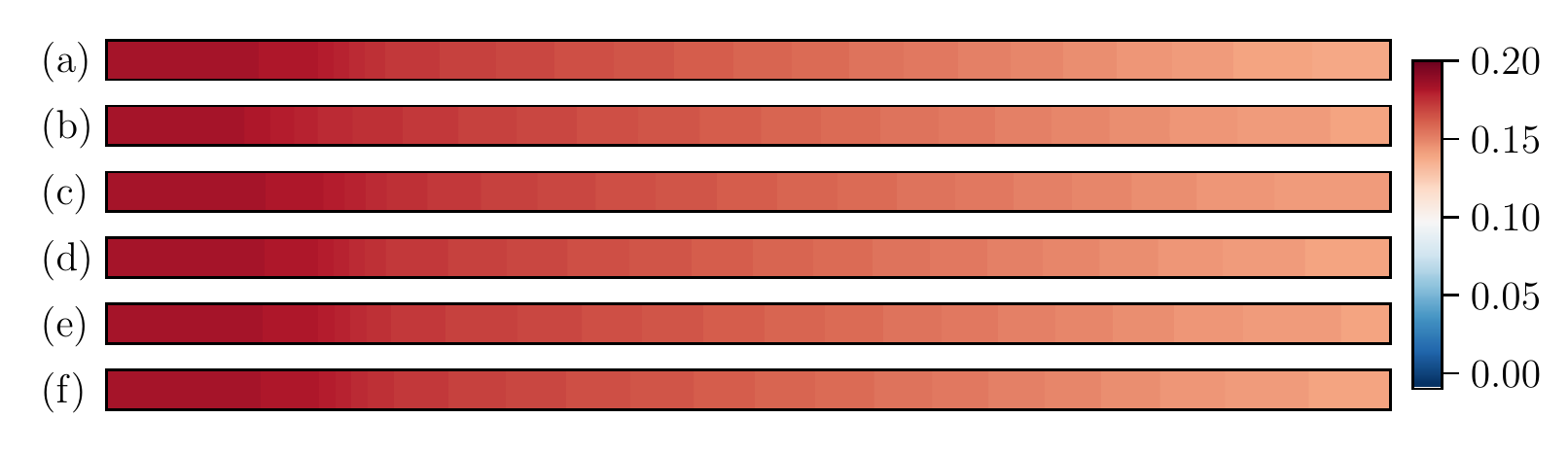}}}
\caption{Average cell temperature computed at $t=0.20$ with a time step size $\Delta{t}=3.15 \times 10^{-5}$ for (a) fine-scale and hybrid simulations with Series expansion with $x_{dist}=$ (b) $0.0\epsilon$, (c) $1.0\epsilon$, (d) $1.5\epsilon$, (e) $3.0\epsilon$ and (f) $4.5\epsilon$.}
\label{fig:cell_temp_1271_loc_series}
\end{figure}

\begin{figure}
\centerline{
 {\includegraphics{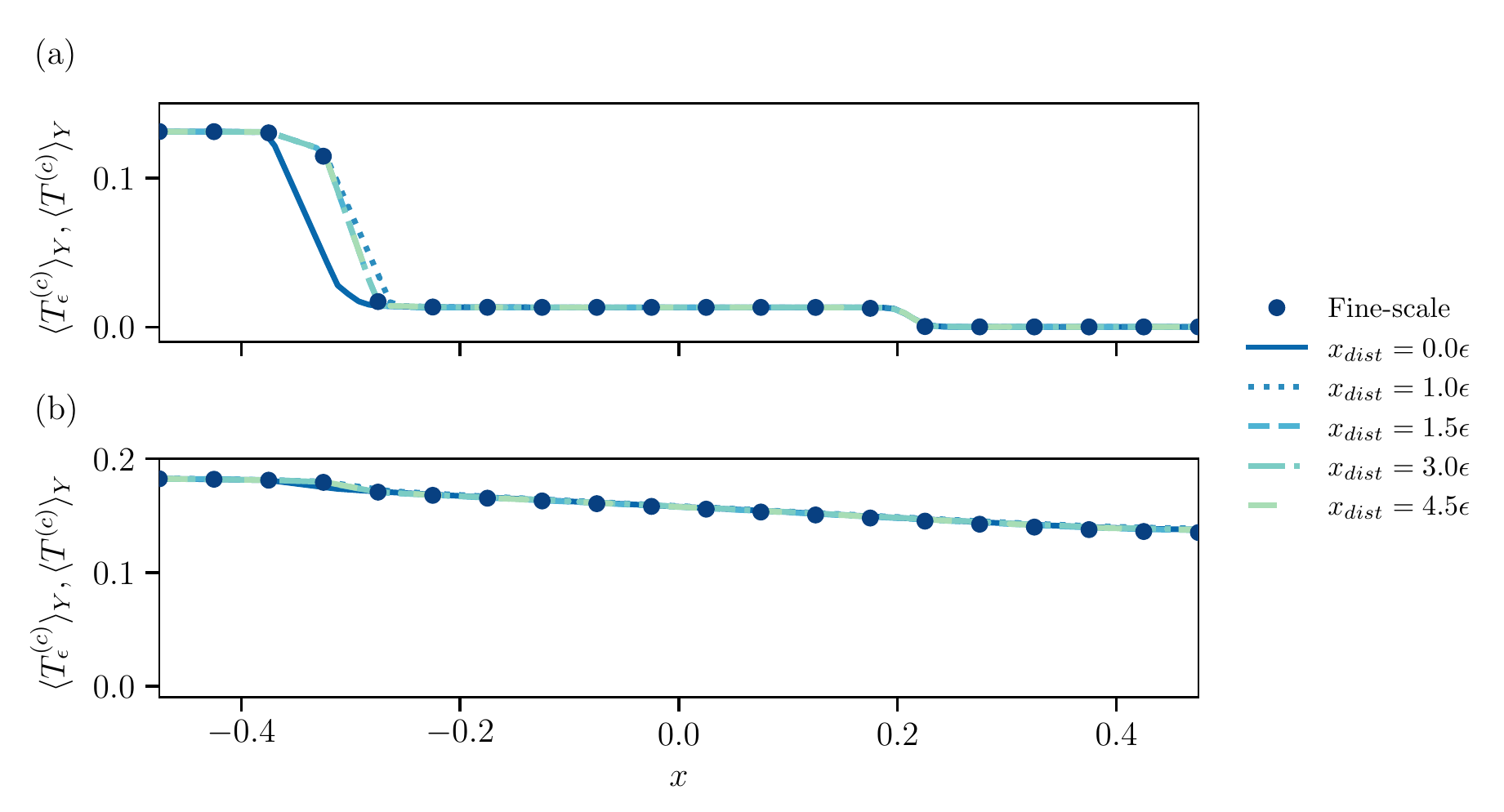}}}
\caption{The centerline ($y=0$) average cell temperature computed at (a) $t=0.02$ and (b) $t=0.20$ with a time step size $\Delta{t}=3.15 \times 10^{-5}$ for fine-scale and hybrid simulations with Series expansion with $x_{dist}=0.0\epsilon$, $1.0\epsilon$, $1.5\epsilon$,  $3.0\epsilon$ and $4.5\epsilon$.}
\label{fig:cell_temp_line_loc_series}
\end{figure}

Figures~\ref{fig:cell_temp_032_loc_taylor} and~\ref{fig:cell_temp_1271_loc_taylor} show the effect of $x_{dist}$ on average cell temperature using the Taylor approach, while Figures~\ref{fig:cell_temp_032_loc_series} and~\ref{fig:cell_temp_1271_loc_series} show the effect using the Series approach for $t=0.02$ and $0.20$. For both Taylor and Series approach, the heat fronts of cases except for $x_{dist} = 0.0\epsilon$ are approximately the same as the fine-scale simulations. Similar trends can be clearly observed with the centerline plots in Figures~\ref{fig:cell_temp_line_loc_taylor} and~\ref{fig:cell_temp_line_loc_series}. 

\begin{figure}
\centerline{
 {\includegraphics{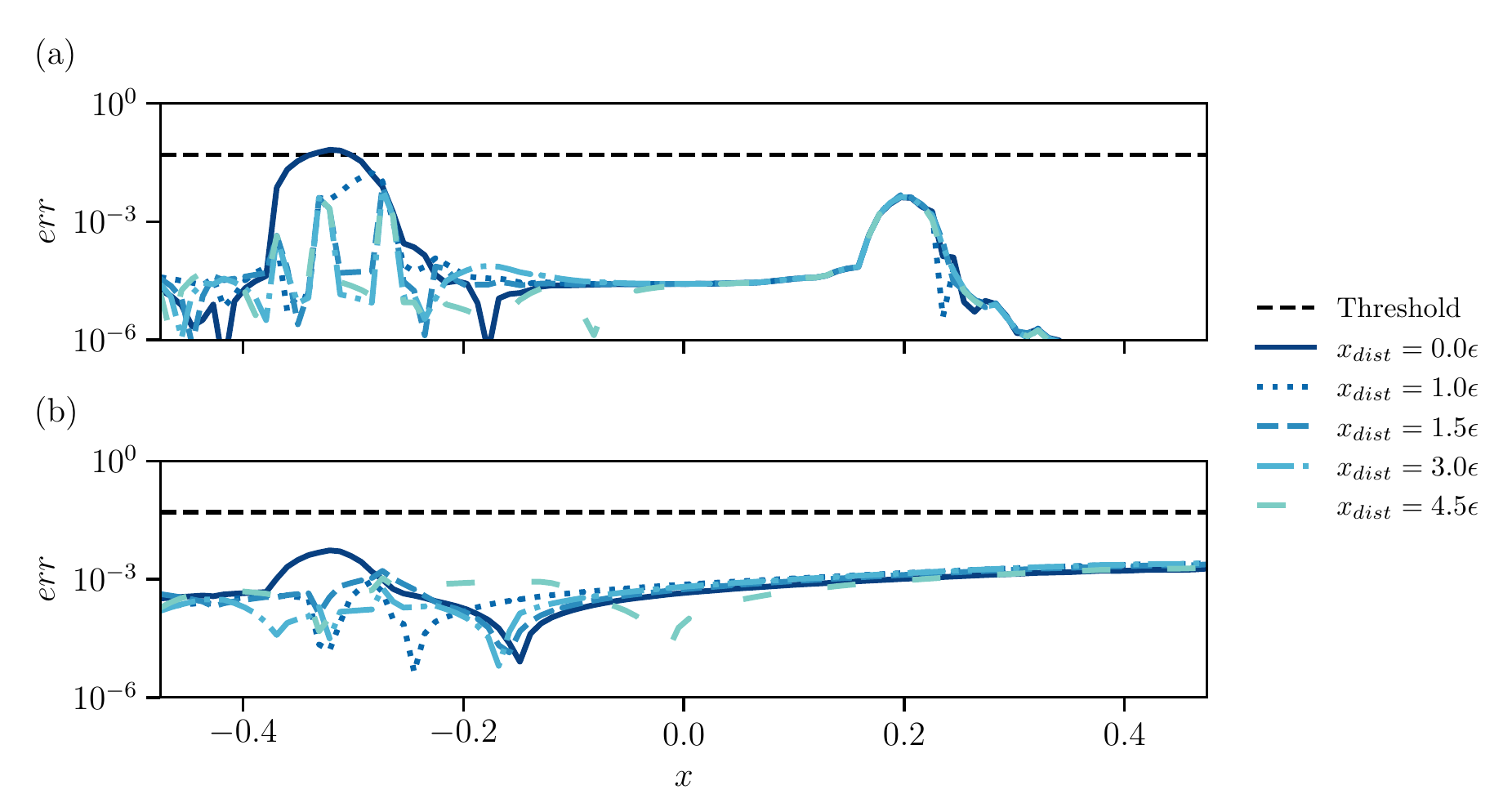}}}

\caption{The centerline ($y=0$) absolute error of the average cell temperature computed at (a) $t=0.02$ and (b) $t=0.20$ with a time step size $\Delta{t}=3.15 \times 10^{-5}$ for hybrid simulations with Taylor expansion with $x_{dist}=0.0\epsilon$, $1.0\epsilon$, $1.5\epsilon$,  $3.0\epsilon$ and $4.5\epsilon$.}
\label{fig:cell_temp_line_loc_taylor_err}
\end{figure}

\begin{figure}
\centerline{
 {\includegraphics{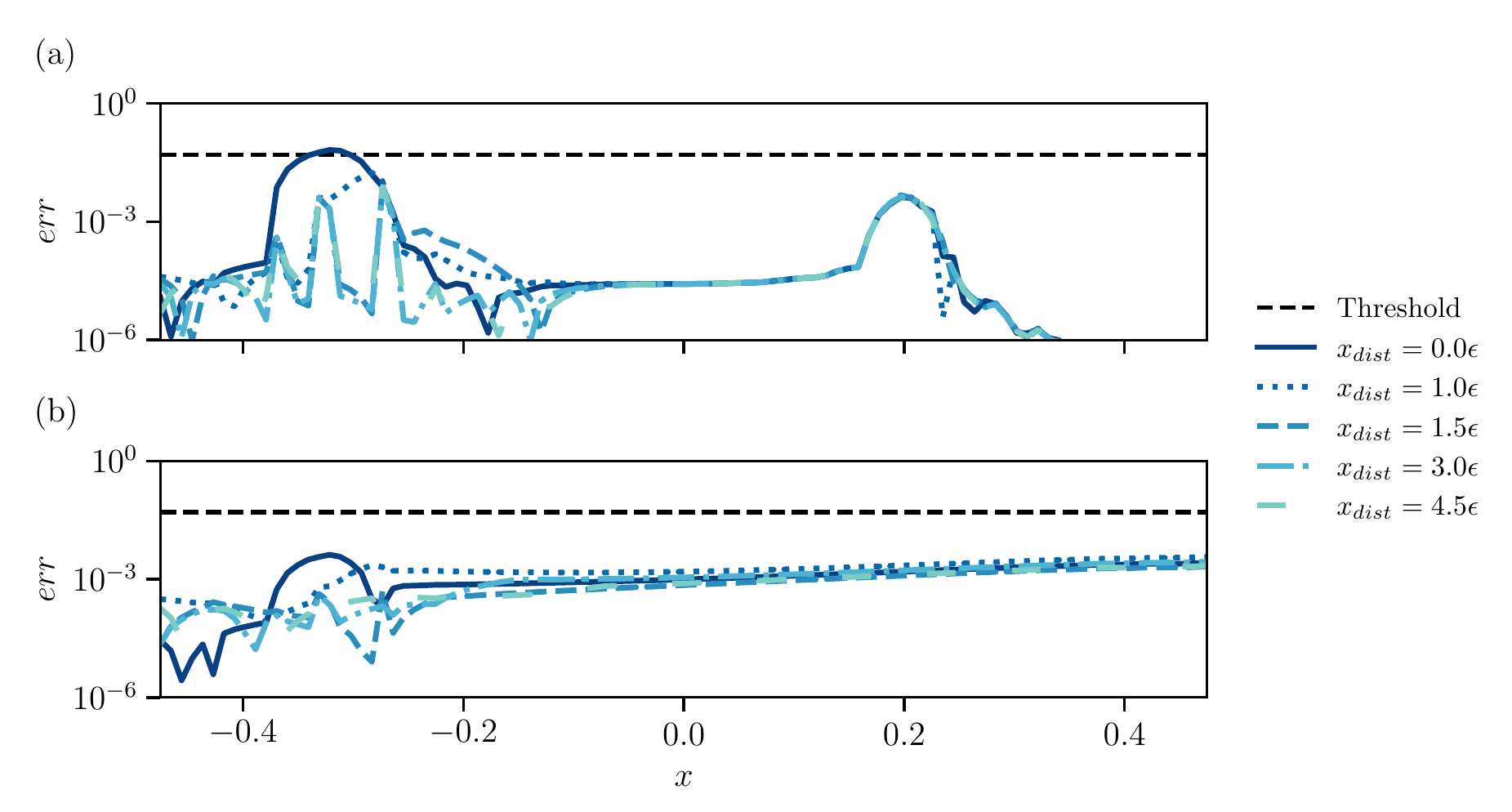}}}

\caption{The centerline ($y=0$) absolute error of the average cell temperature computed at (a) $t=0.02$ and (b) $t=0.20$ with a time step size $\Delta{t}=3.15 \times 10^{-5}$ for hybrid simulations with Series expansion with $x_{dist}=0.0\epsilon$, $1.0\epsilon$, $1.5\epsilon$,  $3.0\epsilon$ and $4.5\epsilon$.}
\label{fig:cell_temp_line_loc_series_err}
\end{figure}

Figures~\ref{fig:cell_temp_line_loc_taylor_err} and~\ref{fig:cell_temp_line_loc_series_err} show the centerline errors computed with equation~\eqref{eq:err-cal} as a function of $x$ for different $x_{dist}$. We can clearly observe that the magnitudes of the errors fall below the upscaling errors when $x_{dist}$ is greater than or equal to the theoretical value ($1.0\epsilon$ for Taylor and $1.5\epsilon$ for Series approaches). Surprisingly, the errors for both the Taylor and Series approaches fall below the threshold when $x_{dist}$ is less than the theoretical value. This shows the robustness of the hybrid coupling algorithm developed. However, the errors are not guaranteed to be bounded within threshold when $x_{dist}$ is less than the theoretical value.   

\subsection{Efficiency of hybrid coupling}
\label{sec:efficiency}
In addition to accuracy, another advantage of hybrid simulation is its reduced computational cost compared to fine-scale simulations~\cite{Battiato2011-zo}. In Section~\ref{sec:acc-hc}, we focus on the accuracy and neglect the computational efficiency by making approximately 80\% of the domain as the fine-scale subdomain. In this section, we instead focus on evaluating the reduction in computational cost by varying the fraction of the fine-scale subdomain. We consider two different domain sizes of $N_{x}^{\left(c\right)}=80$ and $N_{y}^{\left(c\right)}=1$ or $10$, representing small and large domains. The reference values of the parameters are identical to our simulations in Section~\ref{sec:acc-hc} (Table~\ref{tab:ref-val-scaling}). Table~\ref{tab:val-efficiency-sim-params} summarizes the simulation parameters we use. To evaluate reduction in computational cost, the location of the coupling boundary $x_{HC}$ varies from -0.475 to 0.475, resulting in a fine-scale subdomain fraction ranges from 97.5\% to 2.5\%. Since the computational cost highly depends on the number of iterations ($N_{iter}$) in the zero-finding algorithm, we use $N_{iter}=2$ for all simulations to ensure a fair comparison. To evaluate computational cost reduction, the speedup factor is defined as
\begin{align}
   & Speedup = \ddfrac{1}{N_t} \ddfrac{\sum_i^{N_t}t_{\text{fine},i}}{\sum_i^{N_t}t_{HC,i}},
\end{align}
where $N_t$ is the number of time steps and $t_{\text{fine},i}$ and $t_{HC,i}$ are the wall clock time at time step $i$ for fine-scale and hybrid simulations, respectively.
\begin{table}[ht!]
\centering
\caption{
  \label{tab:val-efficiency-sim-params} Summary of simulation parameters used in evaluating the efficiency of the proposed hybrid algorithm.}
\begin{tabular}{l|c}
\hline \hline
Time step size $\Delta{t}/\hat{t}$ &  $3.15 \times 10^{-5}$  \\
Number of unit cells in x-direction ($N_x^{\left(c\right)}$) & $80$ \\
Number of unit cells in y-direction ($N_y^{\left(c\right)}$) & $1$ or $10$ \\
Upscale subdomain minimum grid resolution $h_{\text{up},min}$  &  $1.00 \times 10^{-2}$\\
fine-scale subdomain minimum grid resolution $h_{\text{fine},min}$  &  $2.50\times10^{-4}$\\
Number of timesteps   & 50 \\
Coupling boundary location $x_{HC}$ & -0.475 - 0.475 \\
$\mathcal{R}$ location $\mathbf{x}_{\mathcal{R}}$  & -0.3125 \\
burned location $\mathbf{x}_{\text{burn}}$  & 0.2125 \\
Polynomial order $k$ & 1 \\
\hline \hline
\end{tabular}%
\end{table}

\begin{figure}
\centerline{
 {\includegraphics{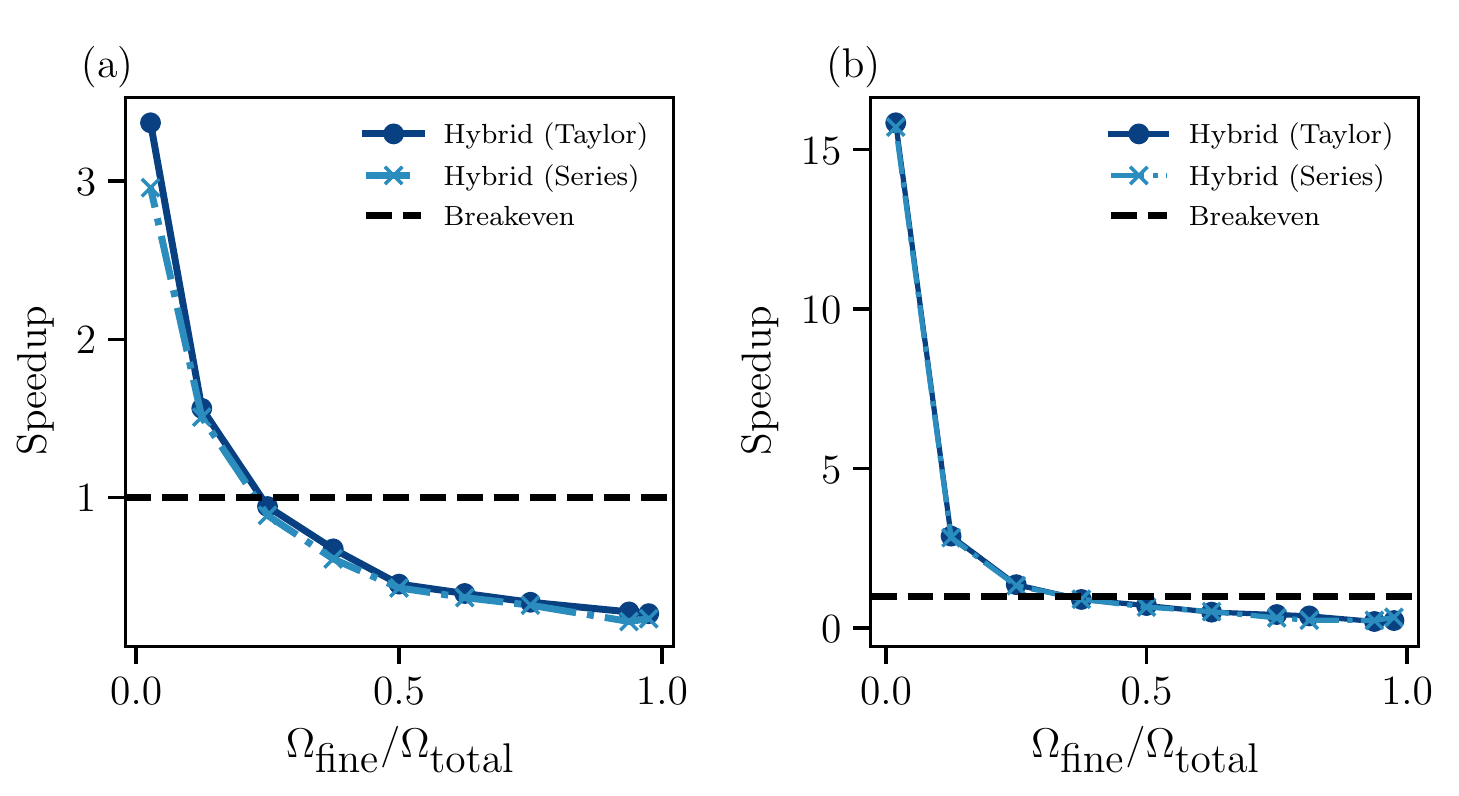}}}

\caption{The speedup factor as a function of fine-scale subdomain volume fraction for domain of (a) $80\times1$ and (b) $80\times10$.}
\label{fig:efficiency}
\end{figure}

Figure~\ref{fig:efficiency} shows the speed-up factor as a function of fine-scale subdomain fraction for two different setups. For the speedup factor $>$ 1, the hybrid approach is favored over the fine-scale approach. Breakeven points (speedup = 1) are approximately 0.2 and 0.4 for $N_{y}^{\left(c\right)}=1$ and $10$, respectively. As the domain size in the $y$ direction increases from $N_{y}^{\left(c\right)}=1$ to $N_{y}^{\left(c\right)}=10$, the breakeven point of the fine-scale subdomain increases from 0.2 to 0.4, indicating that hybrid simulations are more favorable in larger domains. The difference in computational cost between the Taylor and Series expansion approach is negligible except for the small volume fraction. This is expected because the dominant computational cost is related to the linear solvers for high volume fractions. Therefore, computational costs in calculating the coupling conditions are insignificant. 

%% file: conclusion.tex
\section{Conclusion}
We develop a two-way coupled non-intrusive hybrid algorithm to model heat transfer in battery pack. The coupling boundary conditions are derived on the basis of heat and flux conservation. To ensure continuity at the coupling boundary, two methods with different orders of accuracy are proposed by using Taylor and Series expansion approaches, respectively. Weak formulations of the governing equations are derived and implemented in FEniCS for solving the equations numerically.  To validate the accuracy of the proposed two-way coupled hybrid algorithm, we simulated a thermal runway problem with a battery pack. The average temperature profiles of hybrid simulations with both Taylor and Series expansion approaches have been compared with fine-scale and upscaled simulations. The errors are quantified as the absolute difference between hybrid/upscaled and fine-scale simulations. At $t=0.02$, hybrid simulations are able to accurately predict cell and packing temperatures (indicated by negligible errors), while upscaled simulations significantly overpredict the results, as indicated by errors that are significantly higher than upscaling errors. At $t=0.20$, hybrid and upscaled simulations are able to accurately predict the temperature of the battery cell and the packing material by limiting the errors below the upscaling error.  In addition to accuracy, the efficiency of the proposed hybrid algorithms has been investigated. We considered two different configurations with $N_{y}^{\left(c\right)} = 1$ and $10$, representing small and large domains. By varying the fraction of the subdomain on the fine scale from 0.975 to 0.025, we have shown that the breakeven points for $N_{y}^{\left(c\right)} = 1$ and $10$ are 0.2 and 0.4, respectively. Hybrid simulations are favored when the fraction of the fine-scale subdomain is below the breakeven point. In  follow-up studies, an adaptive hybrid algorithm, which dynamically tracks the spatio-temporal region in which applicability conditions are violated, will be developed.

%% file: BCs_pack_temp.tex
\section{Derivation of flux boundary condition for packing temperature}
\label{app:packing-flux-proof}
Apply the averaging operator (equation~\eqref{eq:ave-op}) and multiply by the volume fraction $\phi^{\left(p\right)}$ to the fine-scale packing temperature equations (equation~\eqref{eq:pore_goven_eqs}) to derive the averaged expression as %\textcolor{red}{[$\nabla$ should not have a subscript $\epsilon$ in any of these (see notation in equation (4))]}
\begin{align}
 &\phi^{\left(p\right)}\derive{ \left\langle T_{\epsilon}^{\left(p\right)} \right\rangle_Y }{t} = \phi^{\left(p\right)}\ensmean{k^{\left(p\right)}\grad \bm{\cdot} \grad{T_\epsilon^{\left(p\right)}}}_Y. 
\end{align}
{Here, we would like emphasize that the averaging operator applied to the fine-scale equations depends on $\mathbf{x}$.} 
Integrate over an arbitrary coupling volume $J$ that contains the hybrid coupling boundary $\Gamma^{\left(HC\right)}$, one obtains %\textcolor{red}{[This might be too little too late, but it could be good to somehow clarify (in the notation or otherwise) that $\langle \bm{\cdot} \rangle_{Y}$ is $\mathbf{x}$ dependent. Otherwise the double integral is a little difficult to understand. See previous comment on defining integral over $\mathcal{W}$.]}
\begin{align}
 \int_J \phi^{\left(p\right)}\derive{\left\langle T_{\epsilon}^{\left(p\right)} \right\rangle_Y}{t} \ \dd\mathbf{y} = &\int_{J_{in}} \phi^{\left(p\right)}\ensmean{k^{\left(p\right)}\grad \bm{\cdot} \grad{T_\epsilon^{\left(p\right)}}}_Y \ \dd\mathbf{y} + \int_{J_{out}} \phi^{\left(p\right)}\ensmean{k^{\left(p\right)}\grad \bm{\cdot} \grad{T_\epsilon^{\left(p\right)}}}_Y \ \dd\mathbf{y},   
\end{align}
where $J_{in} = J \cap \Omega_{\text{fine}}$ and $J_{out} = J \cap \Omega_{\text{up}}$ are the partitions of the coupling volume that intersect with the fine-scale and the upscaled domain, respectively. Here, we apply the spatial averaging theorem~{\cite{Yousefzadeh2017-yc,Howes1985-nr}} to the flux term in $J_{in}$ and obtain% \textcolor{red}{[There are two $dy$'s on the LHS. Also, is $\mathcal{V}$ defined anywhere?]}
\begin{align}
 \int_J \phi^{\left(p\right)}\derive{\left\langle T_{\epsilon}^{\left(p\right)} \right\rangle_Y}{t} \ \dd\mathbf{y} &= \int_{J_{in}} \phi^{\left(p\right)}k^{\left(p\right)}\grad \bm{\cdot} \ensmean{\grad{T_\epsilon^{\left(p\right)}}}_Y \ \dd\mathbf{y} \nonumber \\
 &+ \ddfrac{\phi^{\left(p\right)}}{\abs{\mathcal{V}}\abs{J_{in}}} \int_{J_{in}} \int_{\Gamma^{pc}} k^{\left(p\right)} \grad{T_\epsilon^{\left(p\right)}} \bm{\cdot} \mathbf{n}^{\left(p\right)}_\epsilon \ \dd\mathbf{y} \nonumber \\ 
 &+ \int_{J_{out}} \phi^{\left(p\right)}\ensmean{k^{\left(p\right)}\grad \bm{\cdot} \grad{T_\epsilon^{\left(p\right)}}}_Y \ \dd\mathbf{y},
\end{align}
{where $\abs{\mathcal{V}}$ is the volume of battery cells.} Since the upscaled equations are valid in $J_{out}$, we can replace the terms in $J_{out}$ with the upscaled equations such that %\textcolor{red}{[There are two $dy$'s on the LHS.]} 
\begin{align}
 &\int_J \phi^{\left(p\right)}\derive{\left\langle T_{\epsilon}^{\left(p\right)} \right\rangle_Y}{t} \ \dd\mathbf{y} = \int_{J_{in}} \phi^{\left(p\right)}k^{\left(p\right)}\grad \bm{\cdot} \ensmean{\grad{T_\epsilon^{\left(p\right)}}}_Y \ \dd\mathbf{y} \nonumber \\
 &+ \ddfrac{\phi^{\left(p\right)}}{\abs{\mathcal{V}}\abs{J_{in}}} \int_{J_{in}} \int_{\Gamma^{pc}} k^{\left(p\right)} \grad{T_\epsilon^{\left(p\right)}} \bm{\cdot} \mathbf{n}^{\left(p\right)}_\epsilon \ \dd\mathbf{y} \nonumber \\
 &+ \int_{J_{out}} -\left(\U^{\left(p\right)} \bm{\cdot} \nabla_{\mathbf{x}} \langle T^{\left(p\right)} \rangle_{Y} - \V^{\left(p\right)} \bm{\cdot} \nabla_{\mathbf{x}} \langle T^{\left(c\right)} \rangle_{Y} - \nabla_{\mathbf{x}} \bm{\cdot} \left(\K^{\left(p\right)} \bm{\cdot} \nabla_{\mathbf{x}} \langle T^{\left(p\right)} \rangle_{Y}\right) \right)  \ \dd\mathbf{y} \nonumber \\ 
 &+ \int_{J_{out}} \boldsymbol{\mathbb{S}}^{\left(p\right)}\left(t, \mathbf{x}\right) \ \dd\mathbf{y},   
\end{align}
where 
\begin{align}
\boldsymbol{\mathbb{S}}^{\left(p\right)}\left(t, \mathbf{x}\right) = -R_1^{\left(p\right)}\langle T^{\left(p\right)} \rangle_{Y} + R_2^{\left(p\right)}\langle T^{\left(c\right)} \rangle_{Y} - R_3^{\left(p\right)}q^{\left(pw\right)}\left(t, \mathbf{x}\right) + \R_4^{\left(p\right)} \bm{\cdot} \nabla_{\mathbf{x}} q^{\left(pw\right)}\left(t, \mathbf{x}\right).
\end{align}
By applying the divergence theorem and dividing the boundaries into the external boundaries of the coupling volume $\Gamma_{in}$ and $\Gamma_{out}$ and the internal coupling boundary $\Gamma^{\left(HC\right)}$, one obtains %\textcolor{red}{[There are two $dy$'s on the LHS.]}
\begin{align}
\label{eq:packing-flux-rhs}
 \int_J \phi^{\left(p\right)} &\derive{\left\langle T_{\epsilon}^{\left(p\right)} \right\rangle_Y}{t} \ \dd\mathbf{y}  = \int_{\Gamma_{in}} \left(\phi^{\left(p\right)} k^{\left(p\right)} \ensmean{\grad{T_\epsilon^{\left(p\right)}}}_Y\right) \bm{\cdot} \mathbf{n}^{\left(p\right)}_\epsilon \ \dd\mathbf{y} \nonumber \\ 
 &+ \ddfrac{\phi^{\left(p\right)}}{\abs{\mathcal{V}}\abs{J_{in}}} \int_{J_{in}} \int_{\Gamma^{pc}} k^{\left(p\right)} \grad{T_\epsilon^{\left(p\right)}} \bm{\cdot} \mathbf{n}^{\left(p\right)}_\epsilon \ \dd\mathbf{y} \nonumber \\
 &+ \int_{\Gamma_{out}} \left(-\U^{\left(p\right)}\langle T^{\left(p\right)} \rangle_{Y} + \V^{\left(p\right)} \langle T^{\left(c\right)} \rangle_{Y} +  \left(\K^{\left(p\right)} \bm{\cdot} \nabla_{\mathbf{x}} \langle T^{\left(p\right)} \rangle_{Y}\right) \right) \bm{\cdot} \mathbf{n} \ \dd\mathbf{y} \nonumber \\
 &+ \int_{J_{out}} \langle T^{\left(p\right)} \rangle_{Y} \nabla_{\mathbf{x}} \bm{\cdot} \U^{\left(p\right)}  \ \dd\mathbf{y} - \int_{J_{out}} \langle T^{\left(p\right)} \rangle_{Y} \nabla_{\mathbf{x}} \bm{\cdot} \V^{\left(p\right)}  \ \dd\mathbf{y} + \int_{J_{out}} \boldsymbol{\mathbb{S}}^{\left(p\right)}\left(t, \mathbf{x}\right) \ \dd\mathbf{y} \nonumber \\
 &+ \int_{\Gamma^{\left(HC\right)}} \left(\phi^{\left(p\right)} k^{\left(p\right)} \ensmean{\grad{T_\epsilon^{\left(p\right)}}}_Y\right) \bm{\cdot} \mathbf{n}^{\left(p\right)}_\epsilon \ \dd\mathbf{y}  \nonumber \\
 &+ \int_{\Gamma^{\left(HC\right)}} \left(-\U^{\left(p\right)}\langle T^{\left(p\right)} \rangle_{Y} + \V^{\left(p\right)} \langle T^{\left(c\right)} \rangle_{Y} +  \left(\K^{\left(p\right)} \bm{\cdot} \nabla_{\mathbf{x}} \langle T^{\left(p\right)} \rangle_{Y}\right) \right) \bm{\cdot} \mathbf{n}  \ \dd\mathbf{y},   
\end{align}
where $\mathbf{n}$ is the normal vector pointing outwards of the upscaled domain. By applying the spatial average theorem, and divergence theorem to the time-derivative term in equation~\eqref{eq:packing-flux-rhs}, we obtain %\textcolor{red}{[There are two $dy$'s on the LHS.]}
\begin{align}
\label{eq:packing-flux-lhs}
 \int_J \phi^{\left(p\right)} \derive{\left\langle T_{\epsilon}^{\left(p\right)} \right\rangle_Y}{t} \ \dd\mathbf{y}  &= \int_{\Gamma_{in}} \left( \phi^{\left(p\right)} k^{\left(p\right)} \ensmean{\grad{T_\epsilon^{\left(p\right)}}}_Y\right) \bm{\cdot} \mathbf{n}^{\left(p\right)}_\epsilon \ \dd\mathbf{y} \nonumber \\
 &+ \ddfrac{\phi^{\left(p\right)}}{\abs{\mathcal{V}}\abs{J_{in}}}  \int_{J_{in}} \int_{\Gamma^{pc}} k^{\left(p\right)} \grad{T_\epsilon^{\left(p\right)}} \bm{\cdot} \mathbf{n}^{\left(p\right)}_\epsilon \ \dd\mathbf{y} \nonumber \\
 &+\int_{\Gamma_{out}} \left(\phi^{\left(p\right)}k^{\left(p\right)} \ensmean{\grad{T_\epsilon^{\left(p\right)}}}_Y\right) \bm{\cdot} \mathbf{n}^{\left(p\right)}_\epsilon \ \dd\mathbf{y} \nonumber \\
 &+ \ddfrac{\phi^{\left(p\right)}}{\abs{\mathcal{V}}\abs{J_{out}}} \int_{J_{out}} \int_{\Gamma^{pc}} k^{\left(p\right)} \grad{T_\epsilon^{\left(p\right)}} \bm{\cdot} \mathbf{n}^{\left(p\right)}_\epsilon \ \dd\mathbf{y},   
\end{align}
By equating equation~\eqref{eq:packing-flux-rhs} and equation~\eqref{eq:packing-flux-lhs}, the simplified equation is expressed as
\begin{align}
& \int_{\Gamma_{out}} \left(\phi^{\left(p\right)} 
 k^{\left(p\right)} \ensmean{\grad{T_\epsilon^{\left(p\right)}}}_Y\right) \bm{\cdot} \mathbf{n}^{\left(p\right)}_\epsilon \ \dd\mathbf{y} + \ddfrac{\phi^{\left(p\right)}}{\abs{\mathcal{V}}\abs{J_{out}}}\int_{J_{out}} \int_{\Gamma^{pc}} k^{\left(p\right)} \grad{T_\epsilon^{\left(p\right)}} \bm{\cdot} \mathbf{n}^{\left(p\right)}_\epsilon \ \dd\mathbf{y} \nonumber \\ 
&- \int_{\Gamma_{out}} \left(-\U^{\left(p\right)}\langle T^{\left(p\right)} \rangle_{Y} + \V^{\left(p\right)} \langle T^{\left(c\right)} \rangle_{Y} +  \left(\K^{\left(p\right)} \bm{\cdot} \nabla_{\mathbf{x}} \langle T^{\left(p\right)} \rangle_{Y}\right) \right) \bm{\cdot} \mathbf{n} \ \dd\mathbf{y} \nonumber \\
 &- \int_{J_{out}} \langle T^{\left(p\right)} \rangle_{Y} \nabla_{\mathbf{x}} \bm{\cdot} \U^{\left(p\right)}  \ \dd\mathbf{y} + \int_{J_{out}} \langle T^{\left(p\right)} \rangle_{Y} \nabla_{\mathbf{x}} \bm{\cdot} \V^{\left(p\right)}  \ \dd\mathbf{y} - \int_{J_{out}} \boldsymbol{\mathbb{S}}^{\left(p\right)}\left(t, \mathbf{x}\right) \ \dd\mathbf{y} =\nonumber \\
 &+\int_{\Gamma^{\left(HC\right)}} \left(\phi^{\left(p\right)} k^{\left(p\right)} \ensmean{\grad{T_\epsilon^{\left(p\right)}}}_Y\right) \bm{\cdot} \mathbf{n}^{\left(p\right)}_\epsilon \ \dd\mathbf{y} \nonumber \\
 &+ \int_{\Gamma^{\left(HC\right)}} \left(-\U^{\left(p\right)}\langle T^{\left(p\right)} \rangle_{Y} + \V^{\left(p\right)} \langle T^{\left(c\right)} \rangle_{Y} +  \left(\K^{\left(p\right)} \bm{\cdot} \nabla_{\mathbf{x}} \langle T^{\left(p\right)} \rangle_{Y}\right) \right) \bm{\cdot} \mathbf{n}  \ \dd\mathbf{y} , 
\end{align}
The integral over $J_{out}$ and $\Gamma_{out}$ is identically 0 because of the equvalence between fine-scale and upscaled equations, therefore, the packing flux condtion is obtained as
\begin{align}
& \int_{\Gamma^{\left(HC\right)}} \left(\phi^{\left(p\right)} k^{\left(p\right)} \ensmean{\grad{T_\epsilon^{\left(p\right)}}}_Y\right) \bm{\cdot} \mathbf{n}^{\left(p\right)}_\epsilon + \left(-\U^{\left(p\right)}\langle T^{\left(p\right)} \rangle_{Y}   + \V^{\left(p\right)} \langle T^{\left(c\right)} \rangle_{Y} +  \left(\K^{\left(p\right)} \bm{\cdot} \nabla_{\mathbf{x}} \langle T^{\left(p\right)} \rangle_{Y}\right)  \right) \bm{\cdot}  \mathbf{n}   \ \dd\mathbf{y} =0,
 \end{align}
which can be simplified as 
\begin{subequations}
\begin{align}
& \phi^{\left(p\right)} \ensmean{\mathbf{J}_\epsilon^{\left(p\right)}}_{Y} \bm{\cdot} \mathbf{n}^{\left(p\right)}_\epsilon = \ensmean{\mathbf{J}^{\left(p\right)}}_{Y} \bm{\cdot} \mathbf{n}^{\left(p\right)}_\epsilon,
\end{align}
\end{subequations}
where 
\begin{align}
&\ensmean{\mathbf{J}_\epsilon^{\left(p\right)}}_{Y}  =- \left(k^{\left(p\right)}\ensmean{\grad{T_\epsilon^{\left(p\right)}}}_Y\right), \\
&\ensmean{\mathbf{J}^{\left(p\right)}}_{Y} = \U^{\left(p\right)}\langle T^{\left(p\right)} \rangle_{Y}   - \V^{\left(p\right)} \langle T^{\left(c\right)} \rangle_{Y} -  \left(\K^{\left(p\right)} \bm{\cdot} \nabla_{\mathbf{x}} \langle T^{\left(p\right)} \rangle_{Y}\right),
\end{align}
$\mathbf{n}^{\left(p\right)}_\epsilon = -\mathbf{n}$ refers to the relationship between normal vector of the fine-scale and upscaled domains.

%% file: BCs_cell_temp.tex
\section{Derivation of flux boundary condition for cell temperature}
\label{app:cell-flux-proof}
We follow a similar approach in the derivation of the flux boundary conditions for the cell temperature. The main difference between packing and cell temperature is that the cell temperature domain is discontinuous in contrast to the continuous packing domain. We first apply the averaging operator (equation~\eqref{eq:ave-op}) to the fine-scale cell temperature equations (equation~\eqref{eq:pore_goven_eqs}) to derive the averaged expression as 
\begin{align}
 &\phi^{\left(c\right)}\derive{ \left\langle T_{\epsilon}^{\left(c\right)} \right\rangle_Y}{t} = \left(\phi^{\left(c\right)}\varrho \bm{\cdot} \varsigma\right) \left\langle \nabla \bm{\cdot} \left(k^{\left(c\right)} \nabla T_{\epsilon}^{\left(c\right)}\right) \right\rangle_Y + \left(\phi^{\left(c\right)}\varrho \bm{\cdot} \mathcal{R}\right)\left\langle \Pi\left(T_{\epsilon}^{\left(c\right)}, \mathbf{x}\right) \right\rangle_Y,
\end{align}
Integrating over an arbitrary coupling volume $J$ that contains the hybrid coupling boundary $\Gamma^{\left(HC\right)}$ and replacing the fine-scale equations in $J_{out}$ with upscaled equations, one obtains

\begin{align}
 \int_J \phi^{\left(c\right)}\derive{ \left\langle T_{\epsilon}^{\left(c\right)} \right\rangle_Y}{t}  \ \dd\mathbf{y} &= \int_{J_{in}}\left(\phi^{\left(c\right)}\varrho \bm{\cdot} \varsigma\right) \left\langle \nabla \bm{\cdot} \left(k^{\left(c\right)} \nabla T_{\epsilon}^{\left(c\right)}\right) \right\rangle_Y \ \dd\mathbf{y} \\
 &+ \int_{J_{in}} \left(\phi^{\left(c\right)}\varrho \bm{\cdot} \mathcal{R}\right)\left\langle\Pi\left(T_{\epsilon}^{\left(c\right)}, \mathbf{x}\right) \right\rangle_Y \ \dd\mathbf{y} \nonumber   \\
&+ \int_{J_{out}} -\left(\U^{\left(c\right)} \bm{\cdot} \nabla_{\mathbf{x}} \langle T^{\left(c\right)} \rangle_{Y} - \V^{\left(c\right)} \bm{\cdot} \nabla_{\mathbf{x}} \langle T^{\left(p\right)} \rangle_{Y} - \nabla_{\mathbf{x}} \bm{\cdot} \left(\K^{\left(c\right)} \bm{\cdot} \nabla_{\mathbf{x}} \langle T^{\left(c\right)} \rangle_{Y}\right) \right)  \ \dd\mathbf{y} \nonumber \\ 
 &+ \int_{J_{out}} \boldsymbol{\mathbb{S}}^{\left(c\right)}\left(t, \mathbf{x}\right) \ \dd\mathbf{y},
\end{align}
where 
\begin{align}
    \boldsymbol{\mathbb{S}}^{\left(c\right)}\left(t, \mathbf{x}\right) = R_1^{\left(c\right)} \langle T^{\left(p\right)} \rangle_{Y} - R_2^{\left(c\right)} \langle T^{\left(c\right)} \rangle_{Y} + R_3^{\left(c\right)} q^{\left(pw\right)}\left(t, \mathbf{x}\right) + R_4^{\left(c\right)} \overline{\Pi}\left(\langle T^{\left(c\right)} \rangle_{Y}, \mathbf{x}\right).
\end{align}
Applying the spatial averaging theorem, the equations can be expressed as
\begin{align}
 \int_J \phi^{\left(c\right)}\derive{ \left\langle T_{\epsilon}^{\left(c\right)} \right\rangle_Y}{t}  \ \dd\mathbf{y} &= \int_{J_{in}}\left(\phi^{\left(c\right)}\varrho \bm{\cdot} \varsigma\right) \nabla \bm{\cdot} \left\langle \left(k^{\left(c\right)} \nabla T_{\epsilon}^{\left(c\right)}\right) \right\rangle_Y \ \dd\mathbf{y} \nonumber \\
 &+ \ddfrac{\phi^{\left(c\right)}}{\abs{\mathcal{V}}\abs{J_{in}}} \int_{J_{in}} \int_{\Gamma^{\left(pc\right)}} \left(\varrho \bm{\cdot} \varsigma\right) \left(k^{\left(c\right)} \nabla T_{\epsilon}^{\left(c\right)}\right) \bm{\cdot} \mathbf{n}^{\left(c\right)}_\epsilon \ \dd\mathbf{y} \nonumber \\
&+ \int_{J_{in}} \left(\phi^{\left(c\right)}\varrho \bm{\cdot} \mathcal{R}\right)\left\langle\Pi\left(T_{\epsilon}^{\left(c\right)}, \mathbf{x}\right) \right\rangle_Y \ \dd\mathbf{y} \nonumber   \\
&+ \int_{J_{out}} -\left(\U^{\left(c\right)} \bm{\cdot} \nabla_{\mathbf{x}} \langle T^{\left(c\right)} \rangle_{Y} - \V^{\left(c\right)} \bm{\cdot} \nabla_{\mathbf{x}} \langle T^{\left(p\right)} \rangle_{Y} - \nabla_{\mathbf{x}} \bm{\cdot} \left(\K^{\left(c\right)} \bm{\cdot} \nabla_{\mathbf{x}} \langle T^{\left(c\right)} \rangle_{Y}\right) \right)  \ \dd\mathbf{y} \nonumber \\ 
 &+ \int_{J_{out}} \boldsymbol{\mathbb{S}}^{\left(c\right)}\left(t, \mathbf{x}\right) \ \dd\mathbf{y}.
\end{align}
Since the cell domain is discontinous such that $T^c_\epsilon$ exists only for $x \in \mathcal{B}^c$, the divergence theorem cannot be applied. To overcome this, we consider another domain $\mathcal{B}^{c,*}$ such that 
\begin{align}
    &T^{\left(c,*\right)}_\epsilon =
    \begin{cases}
    T^{\left(c\right)}_\epsilon \text{ for } \mathbf{x} \in \mathcal{B}^{\left(c\right)}, \\
    0 \text{ for } \mathbf{x} \notin \mathcal{B}^{\left(c\right)},
    \end{cases}
\end{align}
where $T^{\left(c,*\right)}_\epsilon$ is the cell temperature in the modified domain $\mathcal{B}^{\left(c,*\right)}$. By replacing $\mathcal{B}^{\left(c\right)}$ with $\mathcal{B}^{\left(c,*\right)}$, the cell temperature changes from a discontinuous function to a function with jumps at the interfaces. Here, we introduce a divergence theorem for functions with jumps~\cite{Arnold1982-pq} as %\textcolor{red}{[There's an unlabeled $\mathbf{n}$ here, not sure if you want to specify what it is. It continues into the next couple of equations.]}
\begin{align}
\label{eq:jump-div-theom}
   & \int_V \div{\mathbb{F}} \ \dd\mathbf{x} = \int_{\partial{V}} \mathbb{F}\cdot \mathbf{n} \ \dd\mathbf{x} + \sum_{i,j} \int_{\Gamma^{ij}} (\mathbb{F}^{ij} - \mathbb{F}^{ji}) \bm{\cdot} \mathbf{n}^{ij} \ \dd\mathbf{x},
\end{align}
where $\mathbb{F}$ is an arbitrary flux, $\Gamma^{ij}$ represents the interface with jump and $\mathbf{n}^{ij}$ is the normal vector of interface $i$ pointing toward interface $j$. The modified divergence theorem (equation~\eqref{eq:jump-div-theom}) is reduced to the original divergence theorem if $\Gamma^{ij}$ does not exist. By applying the modified divergence theorem, we obtain
\begin{align}
\label{eq:cell-flux-rhs}
 \int_J \phi^{\left(c\right)}\derive{ \left\langle T_{\epsilon}^{\left(c\right)} \right\rangle_Y}{t}  \ \dd\mathbf{y} &= \int_{\Gamma_{in}}\left(\phi^{\left(c\right)}\varrho \bm{\cdot} \varsigma\right) \left\langle \left(k^{\left(c\right)} \nabla T_{\epsilon}^{\left(c\right)}\right) \right\rangle_Y  \bm{\cdot} \mathbf{n}_\epsilon^{\left(c\right)} \ \dd\mathbf{y} \nonumber \\
 &+ \left(\phi^{\left(c\right)}\varrho \bm{\cdot} \varsigma\right) \sum_{i,j} \int_{\Gamma^{ij}_{in}} k^{\left(c\right)}\left(\ensmean{\nabla T_{\epsilon}^{\left(c\right)}}^{ij}_Y - \ensmean{\nabla T_{\epsilon}^{\left(c\right)}}^{ji}_Y \right) \bm{\cdot}  \mathbf{n}^{ij} \ \dd\mathbf{y} \nonumber \\
 &+ \ddfrac{\phi^{\left(c\right)}}{\abs{\mathcal{V}}\abs{J_{in}}} \int_{J_{in}} \int_{\Gamma^{\left(pc\right)}} \left(\varrho \bm{\cdot} \varsigma\right) \left(k^{\left(c\right)} \nabla T_{\epsilon}^{\left(c\right)}\right) \bm{\cdot} \mathbf{n}^{\left(c\right)}_\epsilon \ \dd\mathbf{y} \nonumber \\
&+ \int_{J_{in}} \left(\phi^{\left(c\right)}\varrho \bm{\cdot} \mathcal{R}\right)\left\langle\Pi\left(T_{\epsilon}^{\left(c\right)}, \mathbf{x}\right) \right\rangle_Y \ \dd\mathbf{y} \nonumber   \\
 &+ \int_{\Gamma_{out}} \left(-\U^{\left(c\right)}\langle T^{\left(c\right)} \rangle_{Y} + \V^{\left(c\right)} \langle T^{\left(p\right)} \rangle_{Y} +  \left(\K^{\left(c\right)} \bm{\cdot} \nabla_{\mathbf{x}} \langle T^{\left(c\right)} \rangle_{Y}\right) \right) \bm{\cdot} \mathbf{n} \ \dd\mathbf{y} \nonumber \\
 &+ \int_{J_{out}} \langle T^{\left(c\right)} \rangle_{Y} \nabla_{\mathbf{x}} \bm{\cdot} \U^{\left(c\right)}  \ \dd\mathbf{y} - \int_{J_{out}} \langle T^{\left(p\right)} \rangle_{Y} \nabla_{\mathbf{x}} \bm{\cdot} \V^{\left(c\right)}  \ \dd\mathbf{y} + \int_{J_{out}} \boldsymbol{\mathbb{S}}^{c}\left(t, \mathbf{x}\right) \ \dd\mathbf{y} \nonumber \\
&+\int_{\Gamma^{\left(HC\right)}}\left(\phi^{\left(c\right)}\varrho \bm{\cdot} \varsigma\right) \left\langle \left(k^{\left(c\right)} \nabla T_{\epsilon}^{\left(c\right)}\right) \right\rangle_Y  \bm{\cdot} \mathbf{n}_\epsilon^{\left(c\right)} \ \dd\mathbf{y} \nonumber \\
 &+ \int_{\Gamma^{\left(HC\right)}}\left(-\U^{\left(c\right)}\langle T^{\left(c\right)} \rangle_{Y} + \V^{\left(c\right)} \langle T^{\left(p\right)} \rangle_{Y} +  \left(\K^{\left(c\right)} \bm{\cdot} \nabla_{\mathbf{x}} \langle T^{\left(c\right)} \rangle_{Y}\right) \right) \bm{\cdot} \mathbf{n} \ \dd\mathbf{y}.
\end{align}
% \begin{align}
%  \int_J \pdv{\ensmean{T_\epsilon^{c}}_Y}{t_\epsilon} \ \dd\mathbf{y} = &\int_{\Gamma_{in}} \left(k^{c} \left(\varrho \bm{\cdot} \varsigma\right) \ensmean{\grad{T_\epsilon^{c}}}_Y \right) \bm{\cdot}  \mathbf{n}^c_\epsilon \ \dd\mathbf{y} \nonumber\\  
%  &+ k^{c} \left(\varrho \bm{\cdot} \varsigma\right) \sum_{i,j} \int_{\Gamma^{ij}_{in}} \left(\ensmean{\grad{T_\epsilon^{c}}}^{ij}_Y - \ensmean{\grad{T_\epsilon^{c}}}^{ji}_Y \right) \bm{\cdot}  \mathbf{n}^{ij} \ \dd\mathbf{y} \nonumber \\
%  &+ \ddfrac{1}{\abs{\mathcal{V}}\abs{J_{in}}} \int_{J_{in}} \int_{\Gamma^{pc}} k^{c}\left(\varrho \bm{\cdot} \varsigma\right) \grad{T_\epsilon^{c}} \bm{\cdot} \mathbf{n}^c_\epsilon \ \dd\mathbf{y} + \int_{J_{in}} \ensmean{\left(\varrho \bm{\cdot} \mathcal{R}\right)\Pi\left(t_\epsilon, \mathbf{x}_\epsilon\right)}_Y \ \dd\mathbf{y} \nonumber\\ 
%  &- \int_{\Gamma_{out}}  \left(  \boldsymbol{\theta}^{cp} \ensmean{T^{p}}_Y  \right) \bm{\cdot} \mathbf{n}^c_\epsilon \ \dd\mathbf{y} + \int_{J_{out}} \boldsymbol{\mathbb{S}}^{c}\left(t, \mathbf{x}\right) \ \dd\mathbf{y} \nonumber \\
%  &+ \int_{\Gamma^{\left(HC\right)}} \left(k^{c} \left(\varrho \bm{\cdot} \varsigma\right) \ensmean{\grad{T_\epsilon^{c}}}_Y \right) \bm{\cdot}  \mathbf{n}^c_\epsilon \ \dd\mathbf{y} - \int_{\Gamma^{\left(HC\right)}}  \left(  \boldsymbol{\theta}^{cp} \ensmean{T^{p}}_Y  \right) \bm{\cdot} \mathbf{n}^c_\epsilon \ \dd\mathbf{y},   
% \end{align}
By applying the spatial average theorem, and modified divergence theorem to the time-derivative term in equation~\eqref{eq:cell-flux-rhs}, we obtain
\begin{align}
\label{eq:cell-flux-lhs}
 &\int_J \phi^{\left(c\right)}\derive{ \left\langle T_{\epsilon}^{\left(c\right)} \right\rangle_Y}{t}  \ \dd\mathbf{y} = \int_{\Gamma_{in}}\left(\phi^{\left(c\right)}\varrho \bm{\cdot} \varsigma\right) \left\langle \left(k^{\left(c\right)} \nabla T_{\epsilon}^{\left(c\right)}\right) \right\rangle_Y  \bm{\cdot} \mathbf{n}_\epsilon^{\left(c\right)} \ \dd\mathbf{y} \nonumber \\
 &+ \left(\phi^{\left(c\right)}\varrho \bm{\cdot} \varsigma\right) \sum_{i,j} \int_{\Gamma^{ij}_{in}} k^{\left(c\right)}\left(\ensmean{\nabla T_{\epsilon}^{\left(c\right)}}^{ij}_Y - \ensmean{\nabla T_{\epsilon}^{\left(c\right)}}^{ji}_Y \right) \bm{\cdot}  \mathbf{n}^{ij} \ \dd\mathbf{y} \nonumber \\
 &+ \ddfrac{\phi^{\left(c\right)}}{\abs{\mathcal{V}}\abs{J_{in}}} \int_{J_{in}} \int_{\Gamma^{\left(pc\right)}} \left(\varrho \bm{\cdot} \varsigma\right) \left(k^{\left(c\right)} \nabla T_{\epsilon}^{\left(c\right)}\right) \bm{\cdot} \mathbf{n}^{\left(c\right)}_\epsilon \ \dd\mathbf{y} \nonumber \\
&+ \int_{J_{in}} \left(\phi^{\left(c\right)}\varrho \bm{\cdot} \mathcal{R}\right)\left\langle\Pi\left(T_{\epsilon}^{\left(c\right)}, \mathbf{x}\right) \right\rangle_Y \ \dd\mathbf{y} \nonumber   \\
&+\int_{\Gamma_{out}}\left(\phi^{\left(c\right)}\varrho \bm{\cdot} \varsigma\right) \left\langle \left(k^{\left(c\right)} \nabla T_{\epsilon}^{\left(c\right)}\right) \right\rangle_Y  \bm{\cdot} \mathbf{n}_\epsilon^{\left(c\right)} \ \dd\mathbf{y} \nonumber \\
 &+ \left(\phi^{\left(c\right)}\varrho \bm{\cdot} \varsigma\right) \sum_{i,j} \int_{\Gamma^{ij}_{out}} k^{\left(c\right)}\left(\ensmean{\nabla T_{\epsilon}^{\left(c\right)}}^{ij}_Y - \ensmean{\nabla T_{\epsilon}^{\left(c\right)}}^{ji}_Y \right) \bm{\cdot}  \mathbf{n}^{ij} \ \dd\mathbf{y} \nonumber \\
 &+ \ddfrac{\phi^{\left(c\right)}}{\abs{\mathcal{V}}\abs{J_{out}}} \int_{J_{out}} \int_{\Gamma^{\left(pc\right)}} \left(\varrho \bm{\cdot} \varsigma\right) \left(k^{\left(c\right)} \nabla T_{\epsilon}^{\left(c\right)}\right) \bm{\cdot} \mathbf{n}^{\left(c\right)}_\epsilon \ \dd\mathbf{y} \nonumber \\
&+ \int_{J_{out}} \left(\phi^{\left(c\right)}\varrho \bm{\cdot} \mathcal{R}\right)\left\langle\Pi\left(T_{\epsilon}^{\left(c\right)}, \mathbf{x}\right) \right\rangle_Y \ \dd\mathbf{y},
\end{align}
By equating equation~\eqref{eq:cell-flux-rhs} and equation~\eqref{eq:cell-flux-lhs}, we obtain
\begin{align}
&\int_{\Gamma_{out}}\left(\phi^{\left(c\right)}\varrho \bm{\cdot} \varsigma\right) \left\langle \left(k^{\left(c\right)} \nabla T_{\epsilon}^{\left(c\right)}\right) \right\rangle_Y  \bm{\cdot} \mathbf{n}_\epsilon^{\left(c\right)} \ \dd\mathbf{y} \nonumber \\
 &+ \left(\phi^{\left(c\right)}\varrho \bm{\cdot} \varsigma\right) \sum_{i,j} \int_{\Gamma^{ij}_{out}} k^{\left(c\right)}\left(\ensmean{\nabla T_{\epsilon}^{\left(c\right)}}^{ij}_Y - \ensmean{\nabla T_{\epsilon}^{\left(c\right)}}^{ji}_Y \right) \bm{\cdot}  \mathbf{n}^{ij} \ \dd\mathbf{y} \nonumber \\
 &+ \ddfrac{\phi^{\left(c\right)}}{\abs{\mathcal{V}}\abs{J_{out}}} \int_{J_{out}} \int_{\Gamma^{\left(pc\right)}} \left(\varrho \bm{\cdot} \varsigma\right) \left(k^{\left(c\right)} \nabla T_{\epsilon}^{\left(c\right)}\right) \bm{\cdot} \mathbf{n}^{\left(c\right)}_\epsilon \ \dd\mathbf{y} \nonumber \\
&+ \int_{J_{out}} \left(\phi^{\left(c\right)}\varrho \bm{\cdot} \mathcal{R}\right)\left\langle\Pi\left(T_{\epsilon}^{\left(c\right)}, \mathbf{x}\right) \right\rangle_Y \ \dd\mathbf{y}= \nonumber \\
&+ \int_{\Gamma_{out}} \left(-\U^{\left(c\right)}\langle T^{\left(c\right)} \rangle_{Y} + \V^{\left(c\right)} \langle T^{\left(p\right)} \rangle_{Y} +  \left(\K^{\left(c\right)} \bm{\cdot} \nabla_{\mathbf{x}} \langle T^{\left(c\right)} \rangle_{Y}\right) \right) \bm{\cdot} \mathbf{n} \ \dd\mathbf{y} \nonumber \\
 &+ \int_{J_{out}} \langle T^{\left(c\right)} \rangle_{Y} \nabla_{\mathbf{x}} \bm{\cdot} \U^{\left(c\right)}  \ \dd\mathbf{y} - \int_{J_{out}} \langle T^{\left(p\right)} \rangle_{Y} \nabla_{\mathbf{x}} \bm{\cdot} \V^{\left(c\right)}  \ \dd\mathbf{y} + \int_{J_{out}} \boldsymbol{\mathbb{S}}^{c}\left(t, \mathbf{x}\right) \ \dd\mathbf{y} \nonumber \\
&+\int_{\Gamma^{\left(HC\right)}}\left(\phi^{\left(c\right)}\varrho \bm{\cdot} \varsigma\right) \left\langle \left(k^{\left(c\right)} \nabla T_{\epsilon}^{\left(c\right)}\right) \right\rangle_Y  \bm{\cdot} \mathbf{n}_\epsilon^{\left(c\right)} \ \dd\mathbf{y} \nonumber \\
 &+ \int_{\Gamma^{\left(HC\right)}}\left(-\U^{\left(c\right)}\langle T^{\left(c\right)} \rangle_{Y} + \V^{\left(c\right)} \langle T^{\left(p\right)} \rangle_{Y} +  \left(\K^{\left(c\right)} \bm{\cdot} \nabla_{\mathbf{x}} \langle T^{\left(c\right)} \rangle_{Y}\right) \right) \bm{\cdot} \mathbf{n} \ \dd\mathbf{y}
\end{align}

Since integrating over $J_{out}$ and $\Gamma_{out}$ gives identically 0, the equation can be simplified as
\begin{align}
\label{eq:cell-flux-cond}
&\int_{\Gamma^{\left(HC\right)}}\left(\phi^{\left(c\right)}\varrho \bm{\cdot} \varsigma\right) \left\langle \left(k^{\left(c\right)} \nabla T_{\epsilon}^{\left(c\right)}\right) \right\rangle_Y  \bm{\cdot} \mathbf{n}_\epsilon^{\left(c\right)} \ \dd\mathbf{y} \nonumber \\
 &=- \int_{\Gamma^{\left(HC\right)}}\left(-\U^{\left(c\right)}\langle T^{\left(c\right)} \rangle_{Y} + \V^{\left(c\right)} \langle T^{\left(p\right)} \rangle_{Y} +  \left(\K^{\left(c\right)} \bm{\cdot} \nabla_{\mathbf{x}} \langle T^{\left(c\right)} \rangle_{Y}\right) \right) \bm{\cdot} \mathbf{n} \ \dd\mathbf{y}.
\end{align}
{Since $\Gamma^{\left(HC\right)}$ is the only defined in the packing materials of the fine-scale battery packing material domain}, the integral on the LHS of equation~\eqref{eq:cell-flux-cond} vanishes, then the integral on the RHS must vanish to satisfy the relationship. Therefore, there is no flux coupling for the temperature of battery cells.

%% file: closure.tex
\section{The Closure Problems}
\label{subsection:Appendix_E_Closure_Problems}
In this section, we describe the formulation of the closure problems required to solve the upscaled equations. The formulation of the closure problem is identical to Pietrzyk et al.~\cite{Pietrzyk2023-ou}.
Additionally, we use periodic boundary conditions at the edges of the rectangular unit-cell domain such that $\langle \chi^{\left(p\right)\left[i\right]} \rangle_{\mathcal{B}^{\left(p\right)}} = 0$, $\langle \chiv^{\left(p\right)\left[3\right]} \rangle_{\mathcal{B}^{\left(p\right)}} = \bm{0}$, $\langle \chi^{\left(c\right)\left[1\right]} \rangle_{\mathcal{B}^{\left(c\right)}} = 0$, and $\langle \chiv^{\left(c\right)\left[2\right]} \rangle_{\mathcal{B}^{\left(c\right)}} = \bm{0}$ for $i \in \{1,2\}$.

\subsection{Closure Problem for $\bm{\chi^{\left(p\right)\left[1\right]}}$}
\begin{subequations}
\neweq{}{-\frac{\mathcal{Q}}{|\mathcal{B}^{\left(p\right)}|}|\Gamma^{\left(pw\right)}| - k^{\left(p\right)} \nabla_{\xiv} \bm{\cdot} \nabla_{\xiv} \chi^{\left(p\right)\left[1\right]} = 0 \quad \text{for } \xiv \in \mathcal{B}^{\left(p\right)},}
\neweq{}{-k^{\left(p\right)}\n^{\left(p\right)} \bm{\cdot} \nabla_{\xiv}\chi^{\left(p\right)\left[1\right]} = 0 \quad \text{for } \xiv \in \Gamma^{\left(pc\right)},}
\neweq{}{-k^{\left(p\right)} \n^{\left(p\right)} \bm{\cdot} \nabla_{\xiv}\chi^{\left(p\right)\left[1\right]} = \mathcal{Q} \quad \text{for } \xiv \in \Gamma^{\left(pw\right)}.}
\end{subequations}

\subsection{Closure Problem for $\bm{\chi^{\left(p\right)\left[2\right]}}$}
\begin{subequations}
\neweq{}{-\frac{\text{Bi}^{\left(p\right)}}{|\mathcal{B}^{\left(p\right)}|}|\Gamma^{\left(pc\right)}| - k^{\left(p\right)} \nabla_{\xiv} \bm{\cdot} \nabla_{\xiv} \chi^{\left(p\right)\left[2\right]} = 0 \quad \text{for } \xiv \in \mathcal{B}^{\left(p\right)},}
\neweq{}{-k^{\left(p\right)}\n^{\left(p\right)} \bm{\cdot} \nabla_{\xiv}\chi^{\left(p\right)\left[2\right]} = \text{Bi}^{\left(p\right)} \quad \text{for } \xiv \in \Gamma^{\left(pc\right)},}
\neweq{}{-k^{\left(p\right)}\n^{\left(p\right)} \bm{\cdot} \nabla_{\xiv}\chi^{\left(p\right)\left[2\right]} = 0 \quad \text{for } \xiv \in \Gamma^{\left(pw\right)}.}
\end{subequations}

\subsection{Closure Problem for $\bm{\chiv^{\left(p\right)\left[3\right]}}$}
\begin{subequations}
\neweq{}{-k^{\left(p\right)} \nabla_{\xiv} \bm{\cdot} \left(\I + \nabla_{\xiv} \chiv^{\left(p\right)\left[3\right]}\right) = \mathbf{0} \quad \text{for } \xiv \in \mathcal{B}^{\left(p\right)},}
\neweq{}{-k^{\left(p\right)} \n^{\left(p\right)} \bm{\cdot} \left(\I + \nabla_{\xiv}\chiv^{\left(p\right)\left[3\right]}\right) = \mathbf{0} \quad \text{for } \xiv \in \Gamma^{\left(pc\right)},}
\neweq{}{-k^{\left(p\right)} \n^{\left(p\right)} \bm{\cdot} \left(\I + \nabla_{\xiv}\chiv^{\left(p\right)\left[3\right]}\right) = \mathbf{0} \quad \text{for } \xiv \in \Gamma^{\left(pw\right)}.}
\end{subequations}

\subsection{Closure Problem for $\bm{\chi^{\left(c\right)\left[1\right]}}$}
\begin{subequations}
\neweq{cell_closure0_1}{\frac{\text{Bi}^{\left(c\right)} \varrho \varsigma}{|\mathcal{B}^{\left(c\right)}|}|\Gamma^{\left(pc\right)}| - k^{\left(c\right)} \varrho \varsigma \nabla_{\xiv} \bm{\cdot} \nabla_{\xiv} \chi^{\left(c\right)\left[1\right]} = 0 \quad \text{for } \xiv \in \mathcal{B}^{\left(c\right)},}
\neweq{cell_closure0_2}{-k^{\left(c\right)}\n^{\left(c\right)} \bm{\cdot} \nabla_{\xiv}\chi^{\left(c\right)\left[1\right]} = -\text{Bi}^{\left(c\right)} \quad \text{for } \xiv \in \Gamma^{\left(pc\right)}.}
\end{subequations}

\subsection{Closure Problem for $\bm{\chiv^{\left(c\right)\left[2\right]}}$}
\begin{subequations}
\neweq{cell_closure1_1}{-k^{\left(c\right)} \varrho  \varsigma \nabla_{\xiv} \bm{\cdot} \left(\I + \nabla_{\xiv} \chiv^{\left(c\right)\left[2\right]}\right) = \mathbf{0} \quad \text{for } \xiv \in \mathcal{B}^{\left(c\right)},}
\neweq{cell_closure1_2}{-k^{\left(c\right)} \n^{\left(c\right)} \bm{\cdot} \left(\I + \nabla_{\xiv}\chiv^{\left(c\right)\left[2\right]}\right) = \mathbf{0} \quad \text{for } \xiv \in \Gamma^{\left(pc\right)}.}
\end{subequations}

%% file: accuracy.tex
\section{Accuracy of hybrid coupling for Packing temperature}
\label{app:pack-temp-appen}

\begin{figure}[H]
\centerline{
 {\includegraphics{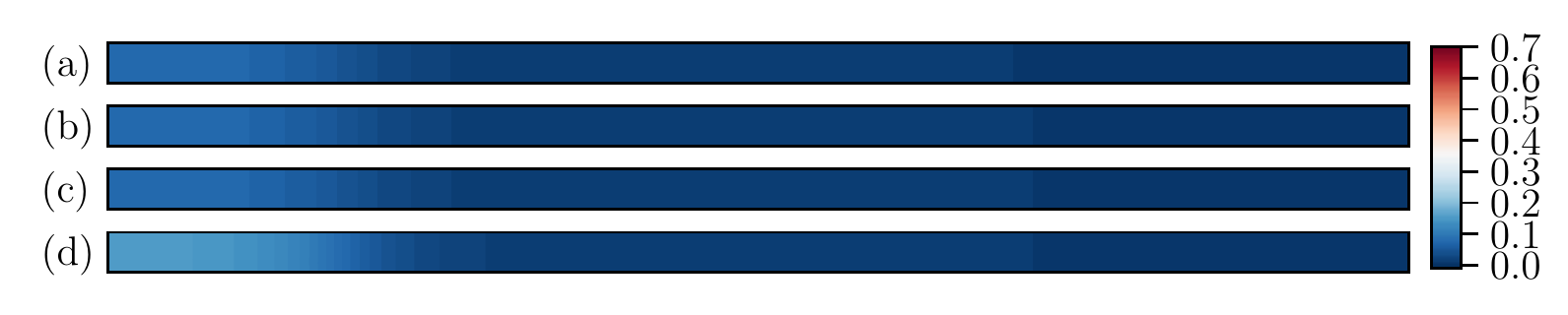}}}

\caption{Average packing temperature computed at $t=0.02$ with a time step size $\Delta{t}=3.15 \times 10^{-5}$ for (a) fine-scale (b) hybrid with Taylor expansion, (c) hybrid with Series expansion and (d) upscaled simulations.}
\label{fig:packing_temp_032_acc}
\end{figure}

\begin{figure}[H]
\centerline{
 {\includegraphics{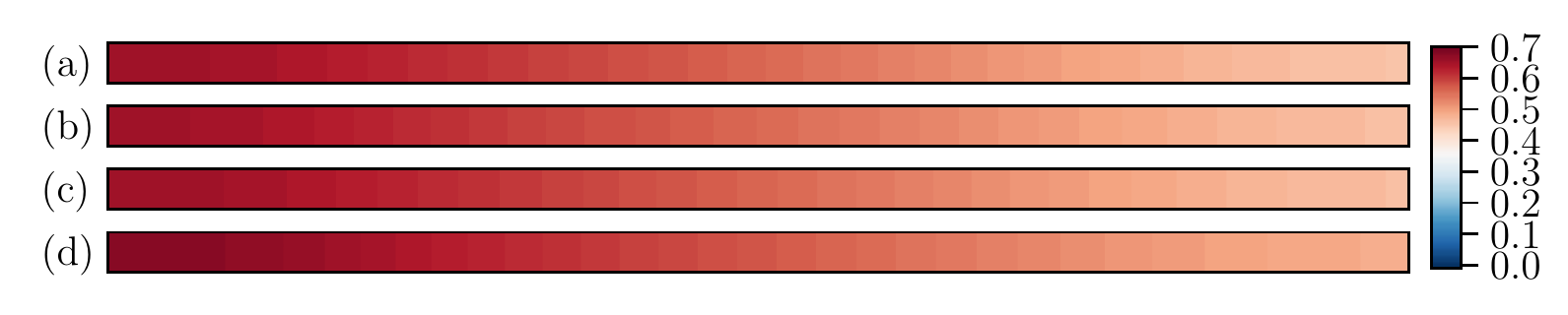}}}

\caption{The average packing temperature computed at $t=0.2$ with a time step size $\Delta{t}=3.15 \times 10^{-5}$ for (a) fine-scale (b) hybrid with Taylor expansion (c) hybrid with Series expansion and (d) upscaled simulations.}
\label{fig:packing_temp_1271_acc}
\end{figure}

\begin{figure}[H]
\centerline{
 {\includegraphics{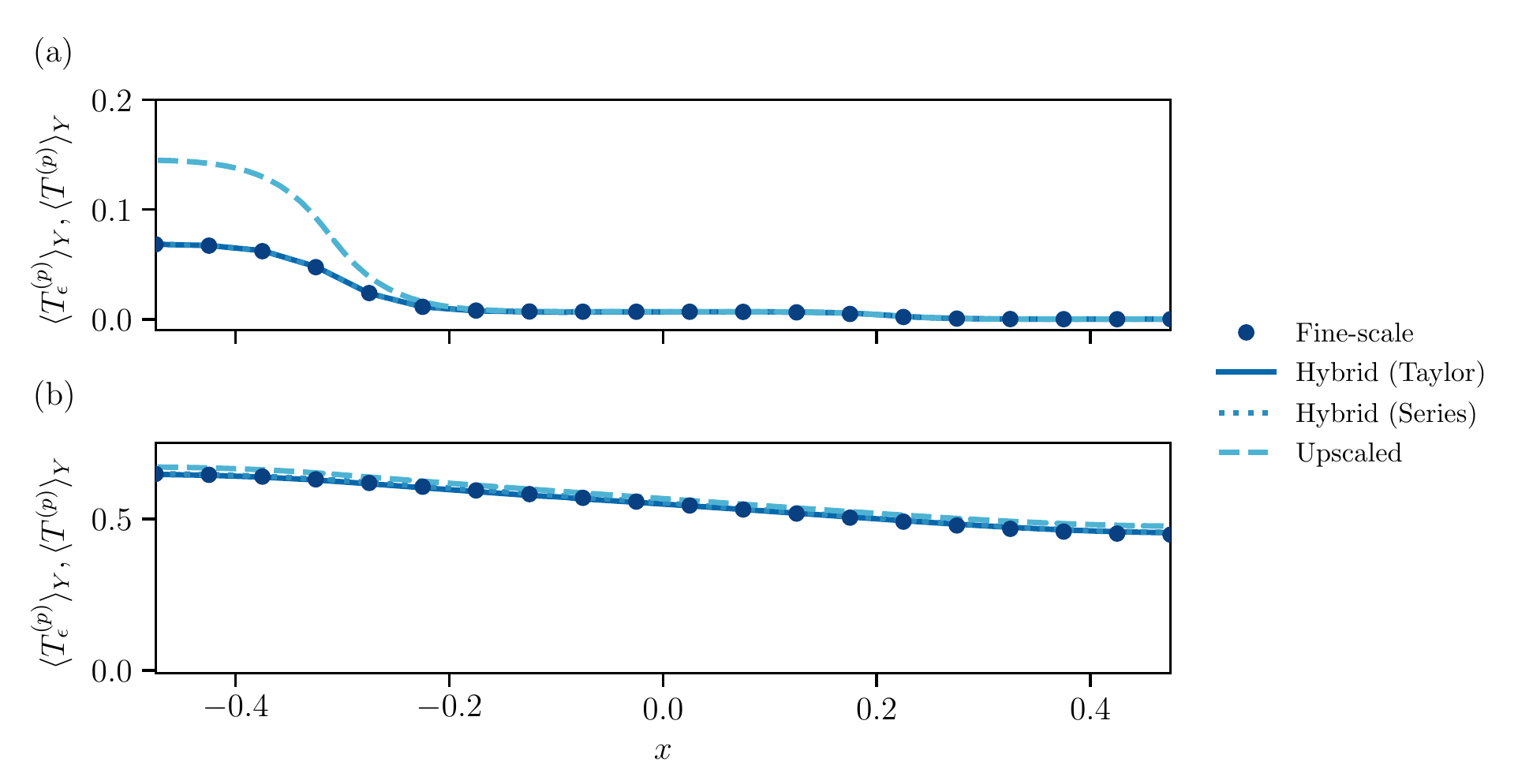}}}

\caption{The centerline ($y=0$) average packing temperature computed at (a) $t=0.02$ and (b) $t=0.20$ with a time step size $\Delta{t}=3.15 \times 10^{-5}$ for the fine scale, hybrid with Taylor and Series expasnion, respectively, and upscaled simulations.}
\label{fig:packing_temp_line_acc}
\end{figure}

\begin{figure}[H]
\centerline{
 {\includegraphics{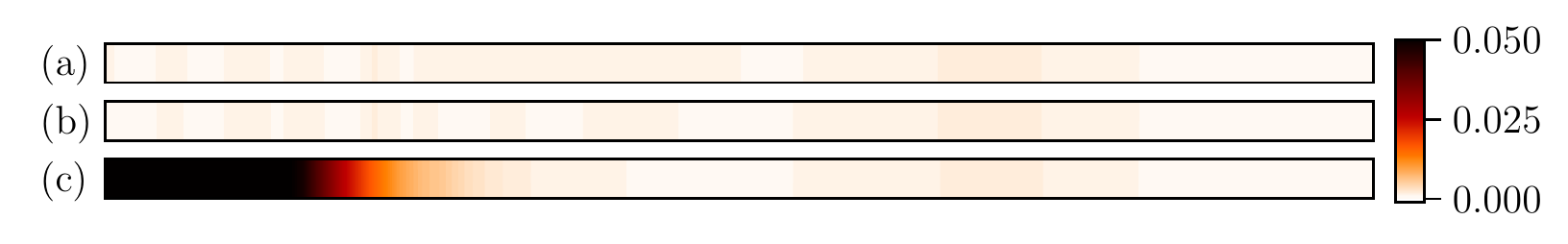}}}

\caption{Absolute error of the average packing temperature computed at $t=0.02$ with a time step size $\Delta{t}=3.15 \times 10^{-5}$ for (a) hybrid with Taylor expansion, (b) hybrid with Series expansion and (c) upscaled simulations.}
\label{fig:packing_temp_err_032_acc}
\end{figure}

\begin{figure}[H]
\centerline{
 {\includegraphics{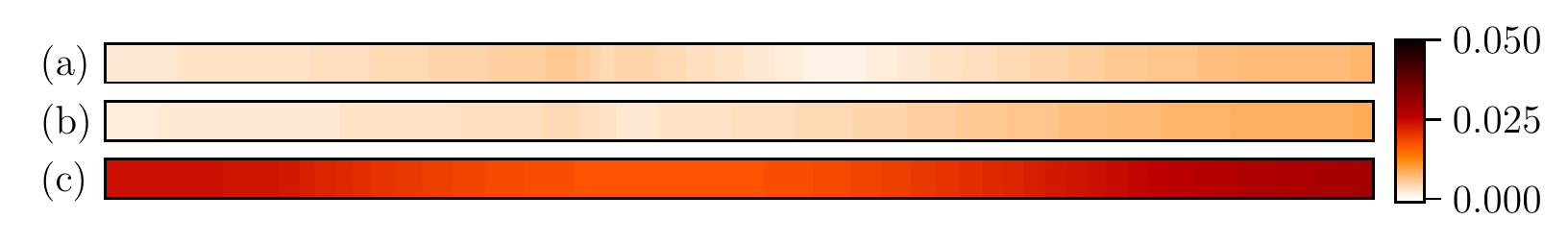}}}

\caption{Absolute error of average packing temperature computed at $t=0.20$ with a time step size $\Delta{t}=3.15 \times 10^{-5}$ for (a) hybrid with Taylor expansion, (b) hybrid with Series expansion and (c) upscaled simulations.}
\label{fig:packing_temp_err_1271_acc}
\end{figure}

\begin{figure}[H]
\centerline{
 {\includegraphics{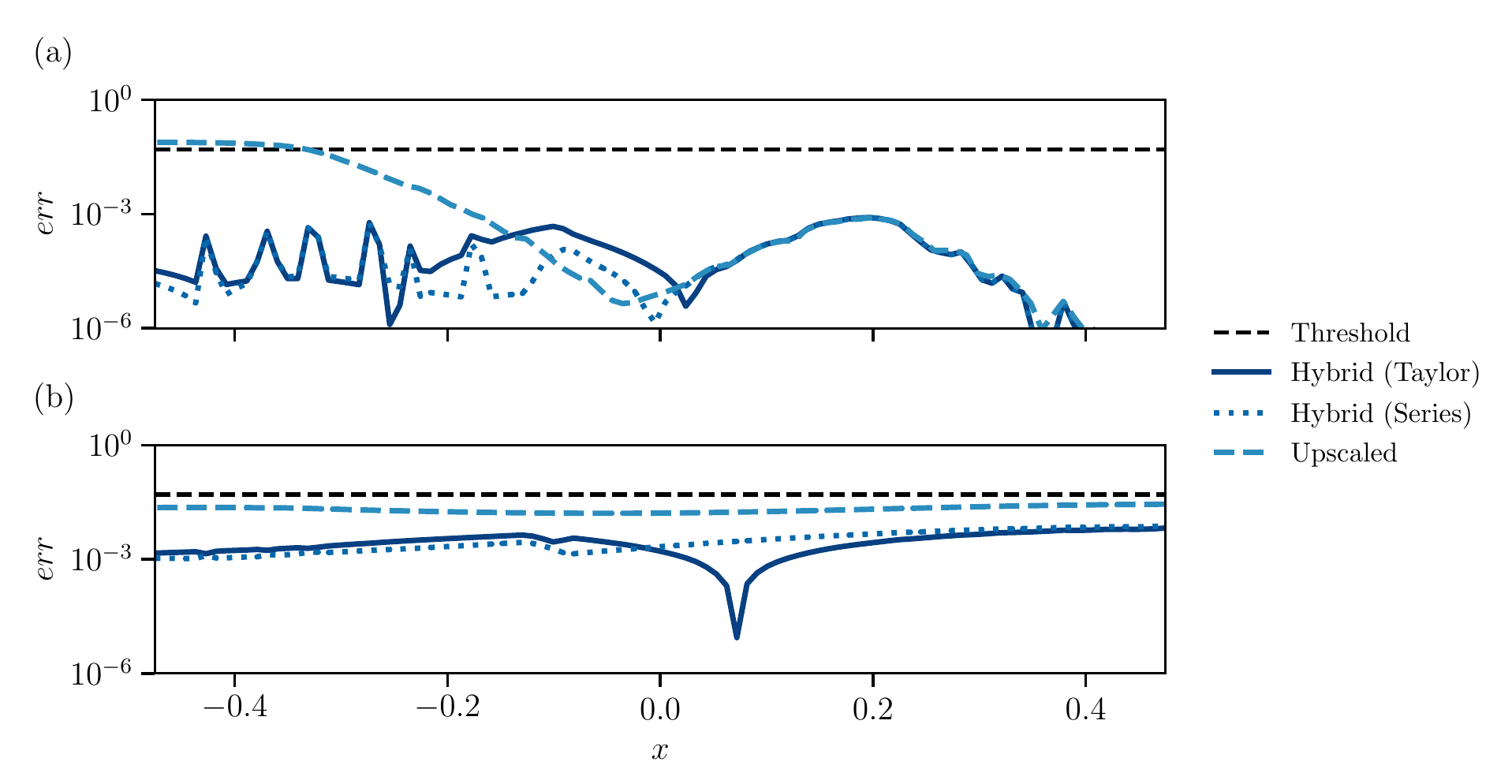}}}
\caption{The centerline ($y=0$) absolute error of the average packing temperature computed at (a) $t=0.02$ and (b) $t=0.20$ with a time step size $\Delta{t}=3.15 \times 10^{-5}$ for hybrid with Taylor and Series expasnion, respectively, and upscaled simulations.}
\label{fig:packing_temp_err_line_acc}
\end{figure}

%% file: coupling.tex
\section{Coupling boundary location on the accuracy of hybrid coupling for packing temperature}
\label{app:pack-temp-loc-appen}

\begin{figure}[H]
\centerline{
 {\includegraphics{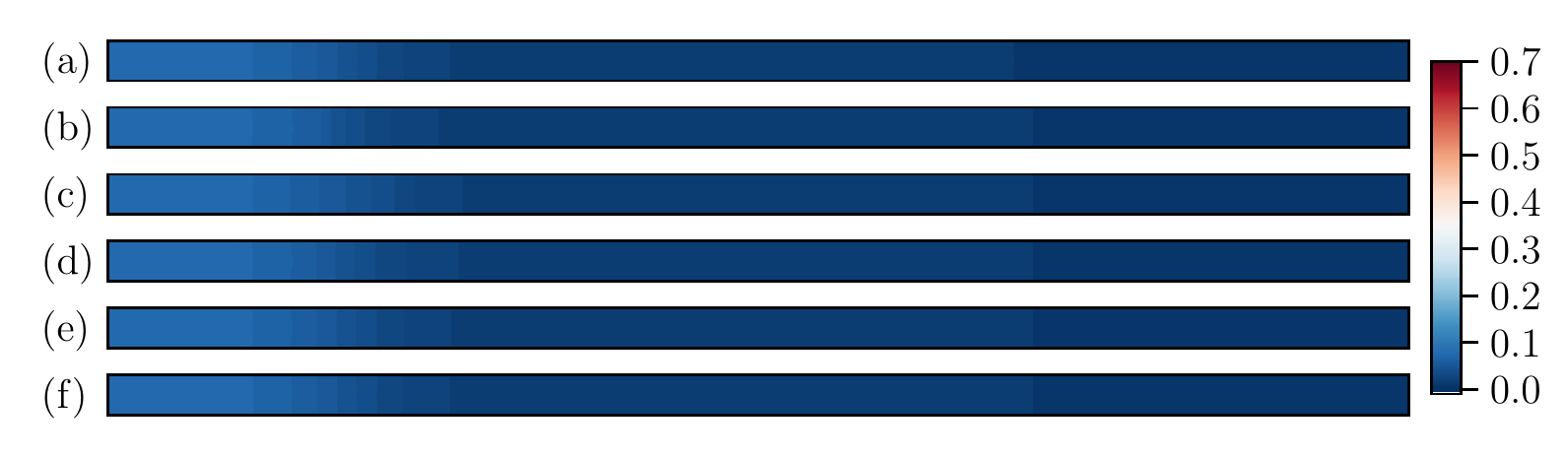}}}

\caption{Average packing temperature computed at $t=0.02$ with a time step size $\Delta{t}=3.15 \times 10^{-5}$ for (a) fine-scale and hybrid simulations with Taylor expansion with $x_{dist}=$ (b) $0.0\epsilon$, (c) $1.0\epsilon$, (d) $1.5\epsilon$, (e) $3.0\epsilon$ and (f) $4.5\epsilon$. }
\label{fig:packing_temp_032_loc_taylor}
\end{figure}

\begin{figure}[H]
\centerline{
 {\includegraphics{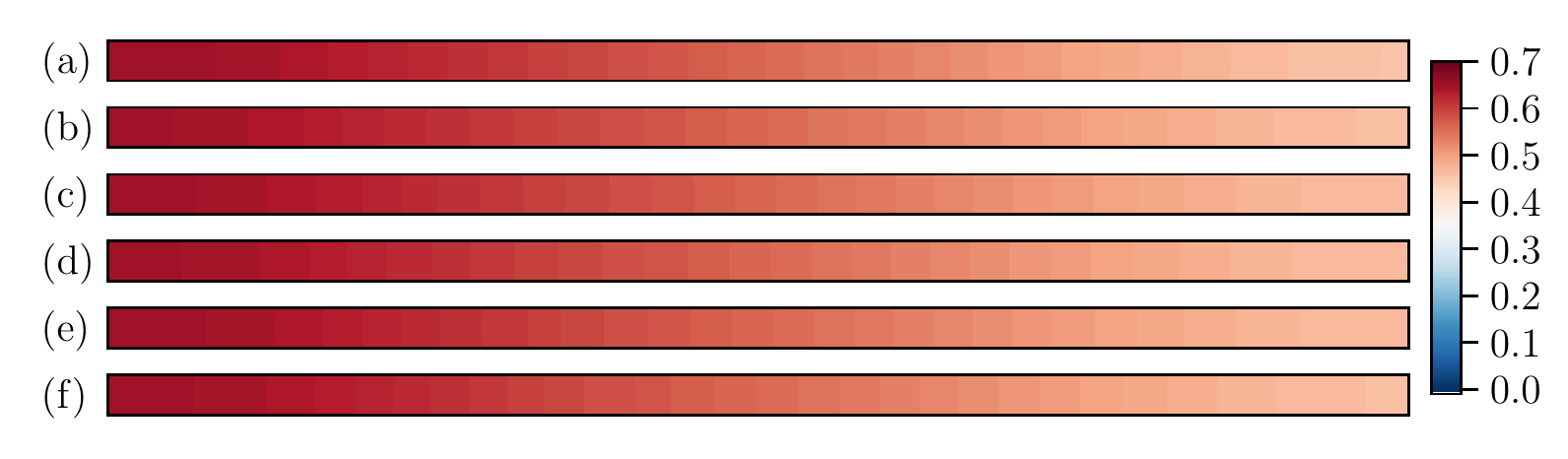}}}

\caption{Average packing temperature computed at $t=0.20$ with a time step size $\Delta{t}=3.15 \times 10^{-5}$ for (a) fine-scale and hybrid simulations with Taylor expansion with $x_{dist}=$ (b) $0.0\epsilon$, (c) $1.0\epsilon$, (d) $1.5\epsilon$, (e) $3.0\epsilon$ and (f) $4.5\epsilon$.}
\label{fig:packing_temp_1271_loc_taylor}
\end{figure}

\begin{figure}[H]
\centerline{
 {\includegraphics{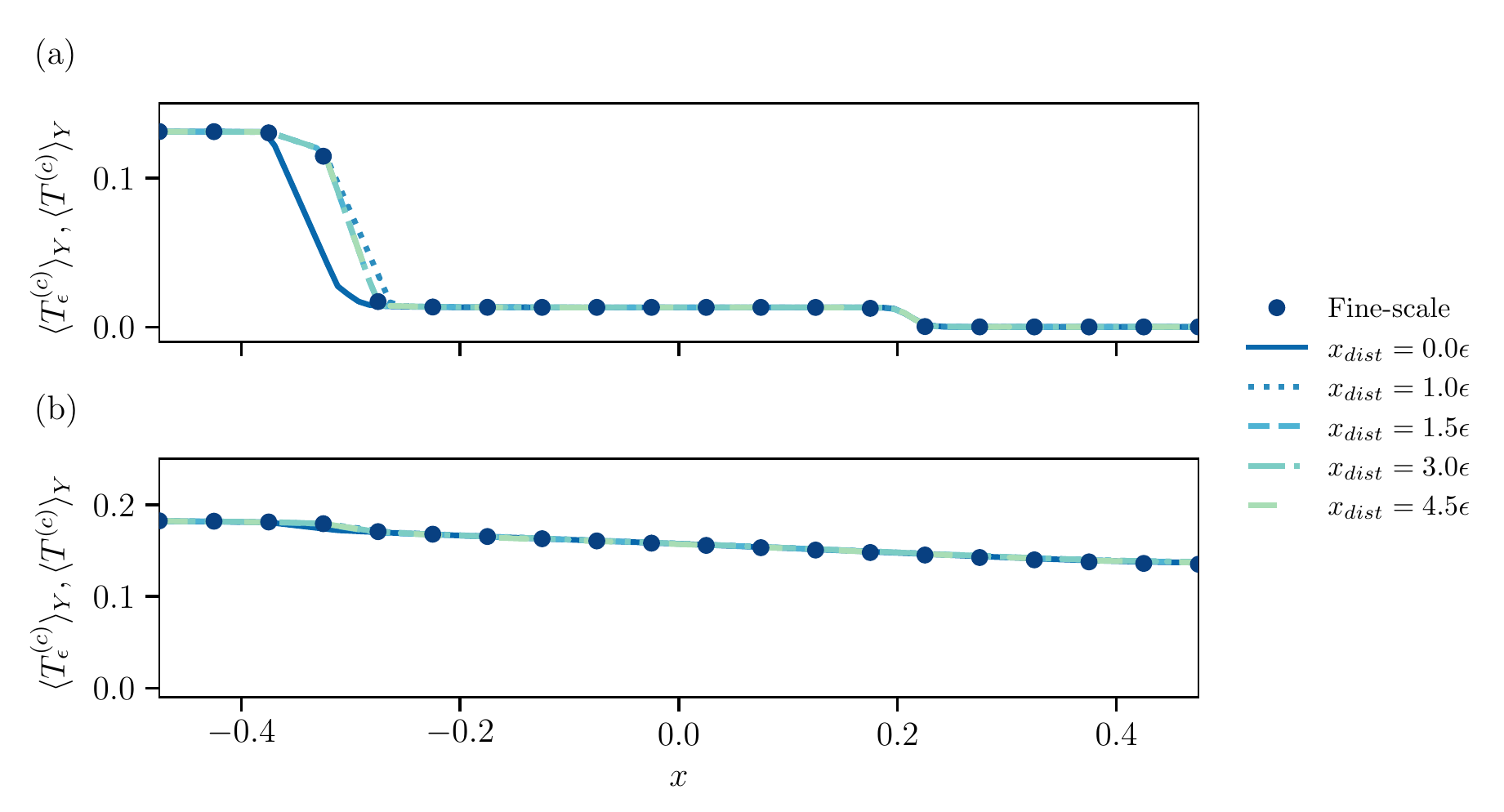}}}

\caption{The centerline ($y=0$) average packing temperature computed at (a) $t=0.02$ and (b) $t=0.20$ with a time step size $\Delta{t}=3.15 \times 10^{-5}$ for fine-scale and hybrid simulations with Taylor expansion with $x_{dist}=0.0\epsilon$, $1.0\epsilon$, $1.5\epsilon$,  $3.0\epsilon$ and $4.5\epsilon$.}
\label{fig:packing_temp_line_loc_taylor}
\end{figure}

\begin{figure}[H]
\centerline{
 {\includegraphics{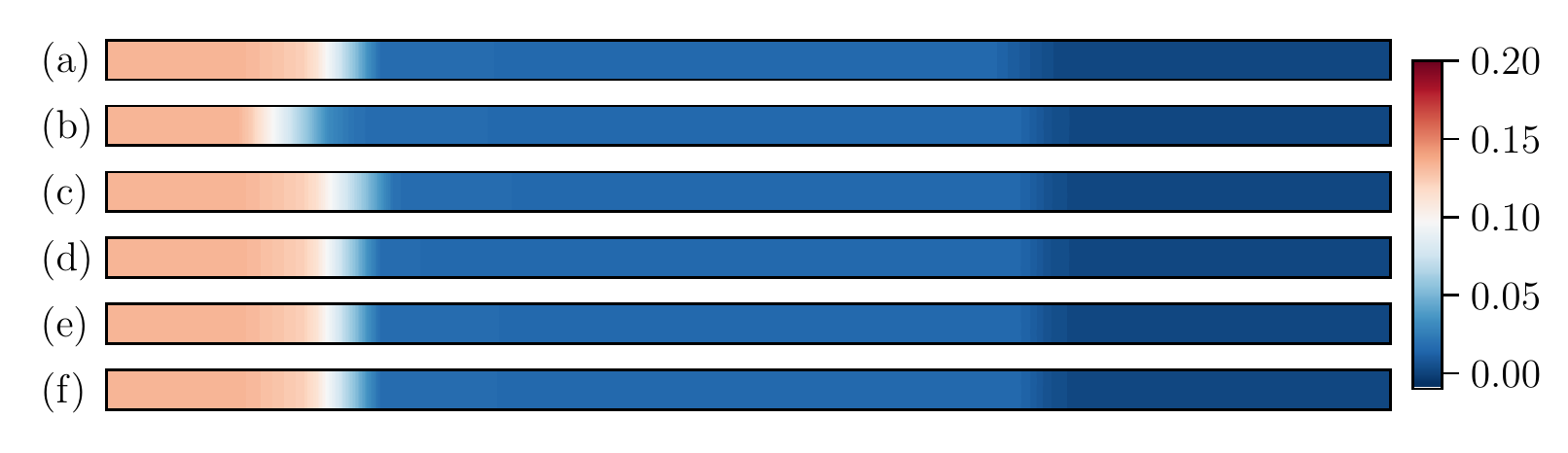}}}

\caption{Average packing temperature computed at $t=0.02$ with a time step size $\Delta{t}=3.15 \times 10^{-5}$ for (a) fine-scale and hybrid simulations with Series expansion with $x_{dist}=$ (b) $0.0\epsilon$, (c) $1.0\epsilon$, (d) $1.5\epsilon$, (e) $3.0\epsilon$ and (f) $4.5\epsilon$.}
\label{fig:packing_temp_032_loc_series}
\end{figure}

\begin{figure}[H]
\centerline{
 {\includegraphics{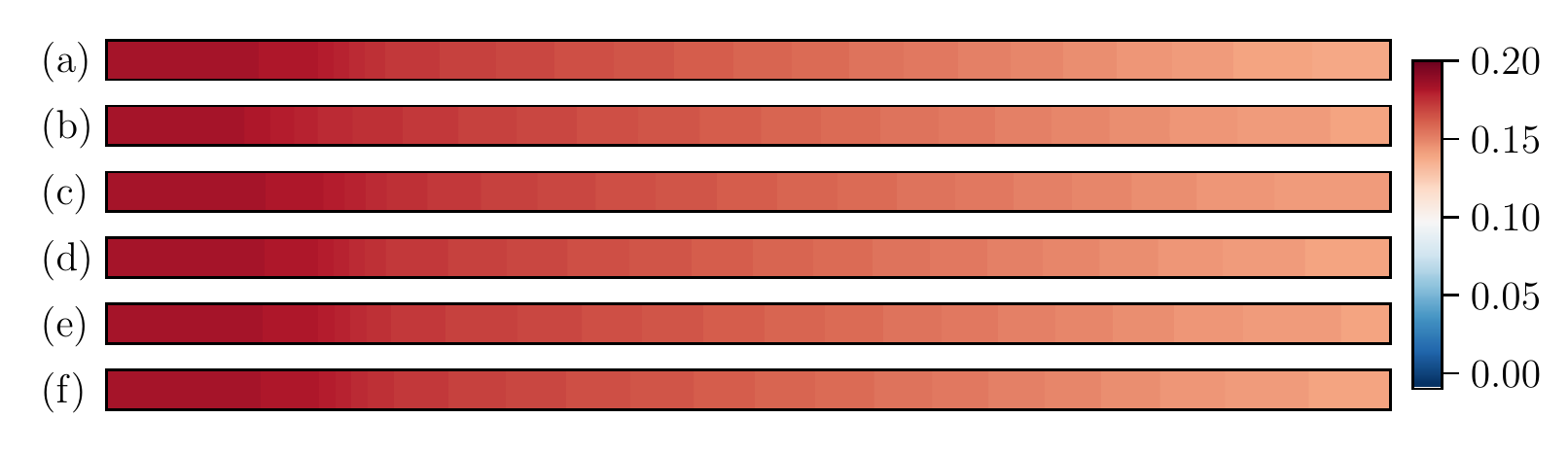}}}

\caption{Average packing temperature computed at $t=0.20$ with a time step size $\Delta{t}=3.15 \times 10^{-5}$ for (a) fine-scale and hybrid simulations with Series expansion with $x_{dist}=$ (b) $0.0\epsilon$, (c) $1.0\epsilon$, (d) $1.5\epsilon$, (e) $3.0\epsilon$ and (f) $4.5\epsilon$.}
\label{fig:packing_temp_1271_loc_series}
\end{figure}

\begin{figure}[H]
\centerline{
 {\includegraphics{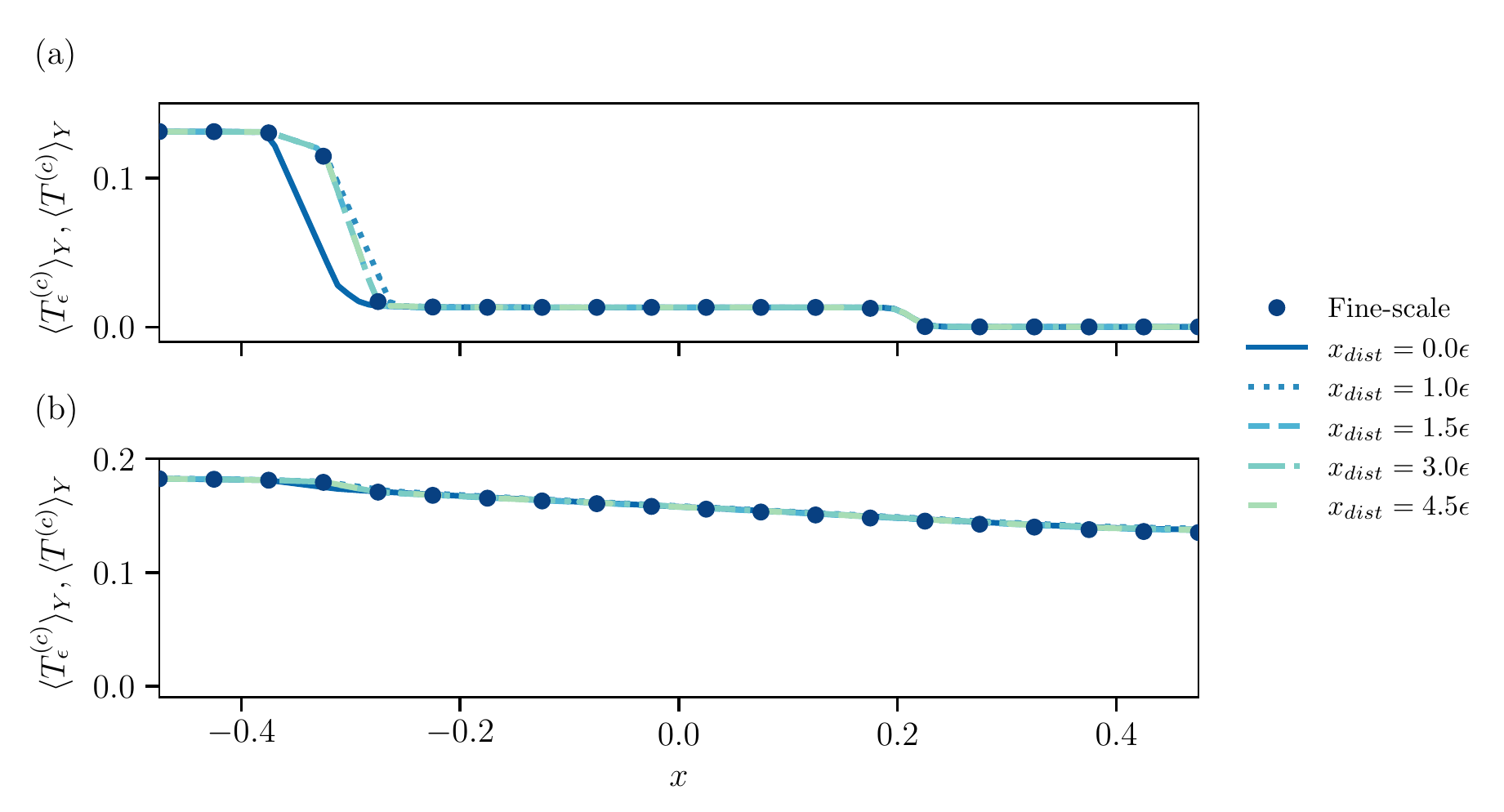}}}

\caption{The centerline ($y=0$) average packing temperature computed at (a) $t=0.02$ and (b) $t=0.20$ with a time step size $\Delta{t}=3.15 \times 10^{-5}$ for fine-scale and hybrid simulations with Series expansion with $x_{dist}=0.0\epsilon$, $1.0\epsilon$, $1.5\epsilon$,  $3.0\epsilon$ and $4.5\epsilon$.}
\label{fig:packing_temp_line_loc_series}
\end{figure}

\begin{figure}[H]
\centerline{
 {\includegraphics{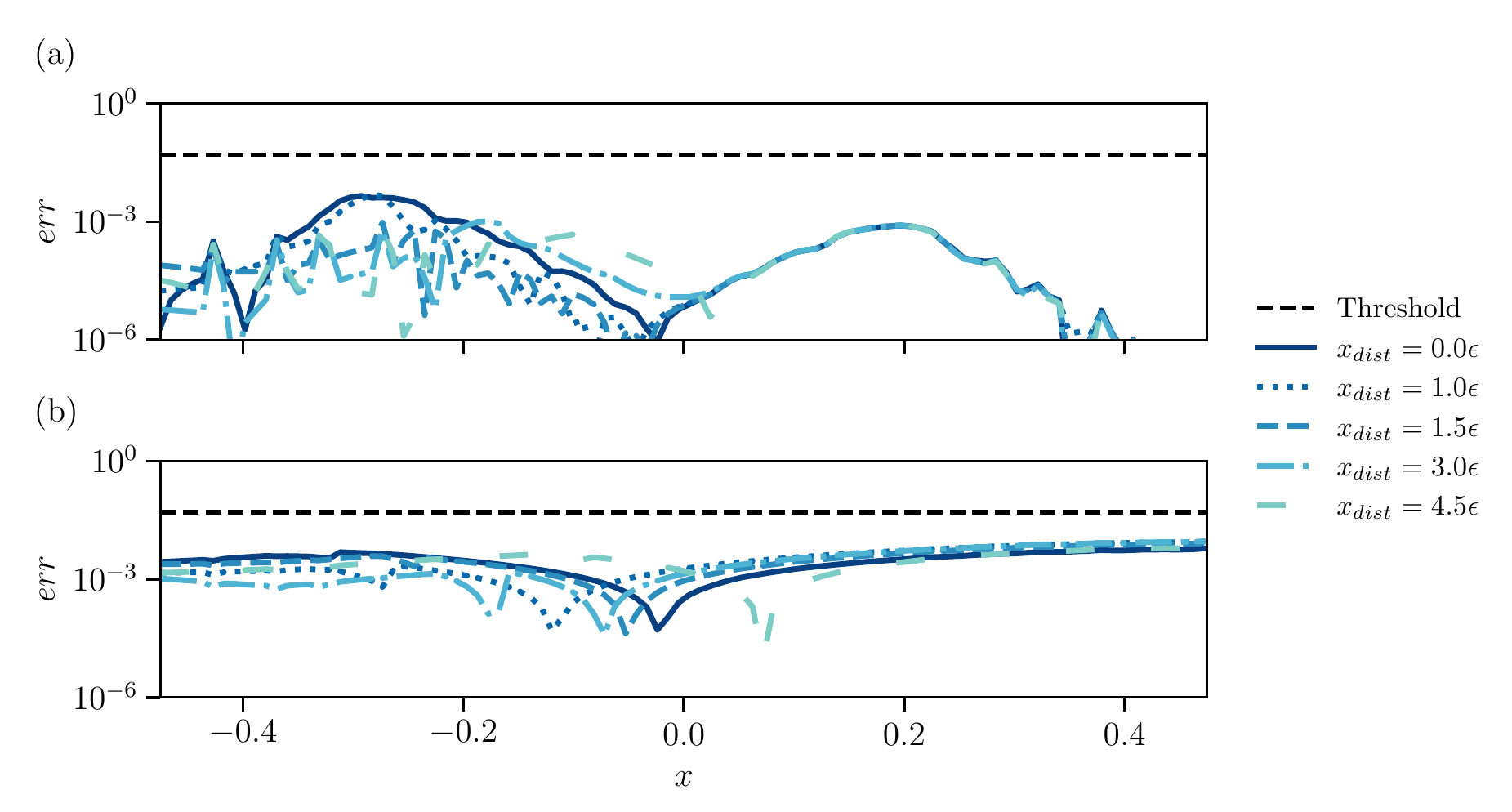}}}

\caption{The centerline ($y=0$) absolute error of the average packing temperature computed at (a) $t=0.02$ and (b) $t=0.20$ with a time step size $\Delta{t}=3.15 \times 10^{-5}$ for hybrid simulations with Taylor expansion with $x_{dist}=0.0\epsilon$, $1.0\epsilon$, $1.5\epsilon$,  $3.0\epsilon$ and $4.5\epsilon$.}
\label{fig:packing_temp_line_loc_taylor_err}
\end{figure}

\begin{figure}[H]
\centerline{
 {\includegraphics{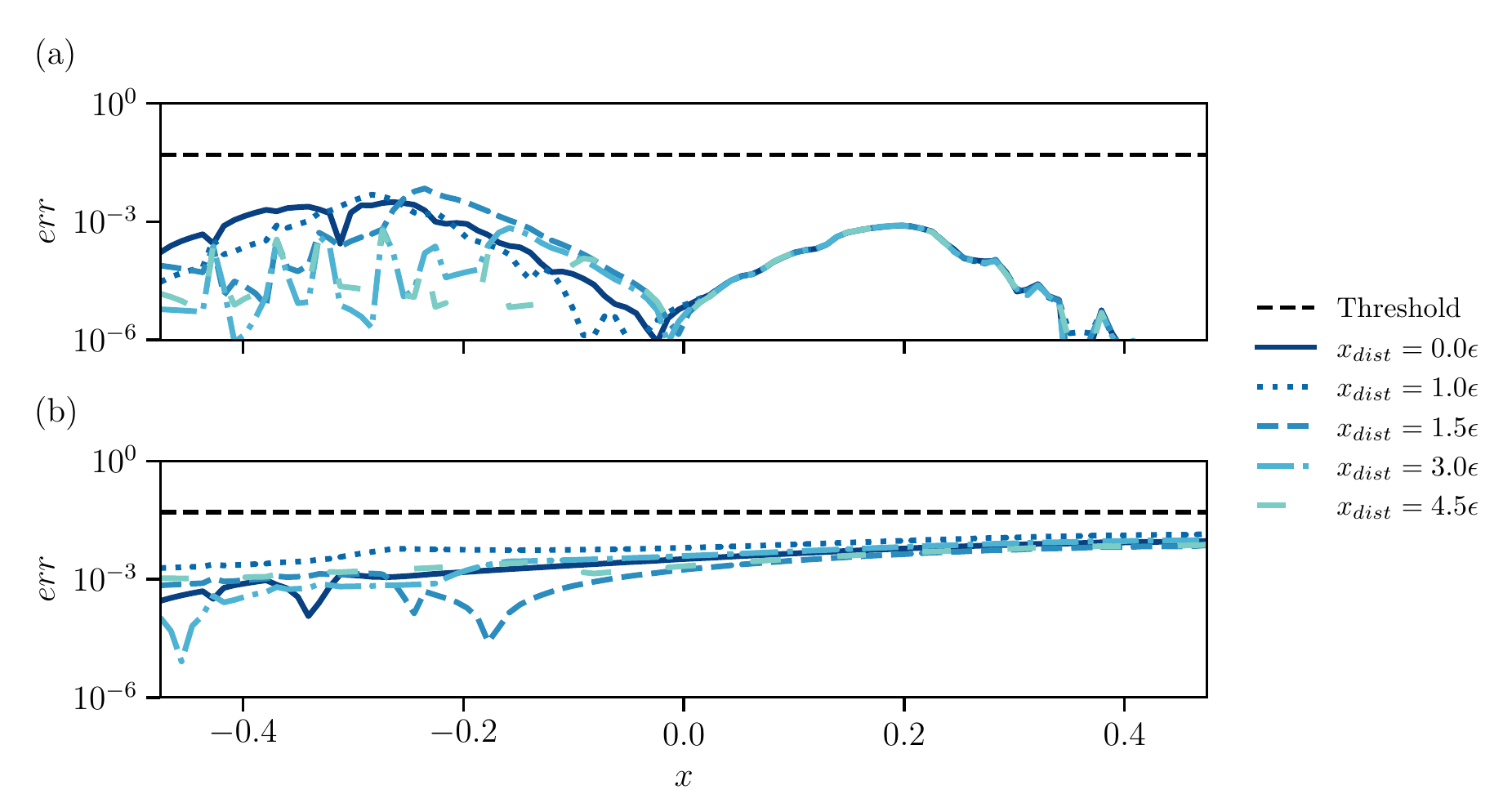}}}

\caption{The centerline ($y=0$) absolute error of the average packing temperature computed at (a) $t=0.02$ and (b) $t=0.20$ with a time step size $\Delta{t}=3.15 \times 10^{-5}$ for hybrid simulations with Series expansion with $x_{dist}=0.0\epsilon$, $1.0\epsilon$, $1.5\epsilon$,  $3.0\epsilon$ and $4.5\epsilon$.}
\label{fig:packing_temp_line_loc_series_err}
\end{figure}